\documentclass{article}
\usepackage{amsfonts}

\begin {document}
\topmargin= -.2in \baselineskip=20pt

\title {Twisted Exponential Sums}
\author {Lei Fu\\
{\small Chern Institute of Mathematics and LPMC, Nankai
University,
Tianjin 300071, P. R. China}\\
{\small leifu@nankai.edu.cn}}
\date{}
\maketitle

\centerline {\bf Abstract}

\bigskip
Let $k$ be a finite field of characteristic $p$, $l$ a prime number
distinct to $p$, $\psi:k\to \overline {\bf Q}_l^\ast$ a nontrivial
additive character, and $\chi:{k^\ast}^n\to \overline{\bf Q}_l^\ast$
a character on ${k^\ast}^n$. Then $\psi$ defines an Artin-Schreier
sheaf ${\cal L}_\psi$ on the affine line ${\bf A}_k^1$, and $\chi$
defines a Kummer sheaf ${\cal K}_\chi$ on the $n$-dimensional torus
${\bf T}_k^n$. Let $f\in k[X_1,X_1^{-1},\ldots, X_n,X_n^{-1}]$ be a
Laurent polynomial. It defines a $k$-morphism $f:{\bf T}_k^n\to {\bf
A}_k^1$. In this paper, we calculate the dimensions and weights of
$H_c^i({\bf T}_{\bar k}^n, {\cal K}_\chi\otimes f^\ast {\cal
L}_\psi)$ under some non-degeneracy conditions on $f$. Our results
can be used to estimate sums of the form
$$\sum_{x_1,\ldots, x_n\in k^\ast} \chi_1(f_1(x_1,\ldots,
x_n))\cdots \chi_m(f_m(x_1,\ldots, x_n))\psi(f(x_1,\ldots,
x_n)),$$ where $\chi_1,\ldots, \chi_m:k^\ast\to {\bf C}^\ast$ are
multiplicative characters, $\psi:k\to {\bf C}^\ast$ is a
nontrivial additive character, and $f_1,\ldots, f_m, f$ are
Laurent polynomials.

\bigskip



\bigskip
\bigskip
\centerline {\bf 0. Introduction}

\bigskip
\bigskip
Let $k$ be a finite field with $q$ elements of characteristic $p$,
let $\chi_1,\ldots, \chi_m:k^\ast\to {\bf C}^\ast$ be nontrivial
multiplicative characters, let $\psi:k\to {\bf C}^\ast$ be a
nontrivial additive character, and let $$f_1(X_1,\ldots,
X_n),\ldots, f_m(X_1,\ldots, X_n), f(X_1,\ldots, X_n)\in
k[X_1,X_1^{-1},\ldots, X_n,X_n^{-1}]$$ be Laurent polynomials. We
make the convention that $\chi_i(0)=0$ $(i=1,\ldots, m)$. In number
theory, we are often lead to study the sum
$$S_1=\sum_{x_1,\ldots, x_n\in k^\ast} \chi_1(f_1(x_1,\ldots,
x_n))\cdots \chi_m(f_m(x_1,\ldots, x_n))\psi(f(x_1,\ldots, x_n)).$$
For this purpose, let's consider another sum
\begin{eqnarray*}
S_2 &=&\sum_{x_1,\ldots, x_{n+m}\in k^\ast}
\chi_1^{-1}(x_{n+1})\cdots \chi_m^{-1}(x_{n+m})\\
&&\qquad\qquad\qquad \psi\left(f(x_1,\ldots,
x_n)+x_{n+1}f_1(x_1,\ldots, x_n)+\cdots +x_{n+m}f_m(x_1,\ldots,
x_n)\right).
\end{eqnarray*}
We have
\begin{eqnarray*}
&&S_2\\&=&\sum_{x_1,\ldots, x_n\in k^\ast}\sum_{x_{n+1}\ldots,
x_{n+m}\in k^\ast}
\left(\chi_1^{-1}(x_{n+1})\psi(x_{n+1}f_1(x_1,\ldots,
x_n))\right)\cdots
\left(\chi_m^{-1}(x_{n+m})\psi(x_{n+m}f_m(x_1,\ldots,
x_n))\right)\\
&&\qquad\qquad\qquad \qquad\qquad\qquad\psi(f(x_1,\ldots, x_n))\\
&=& \sum_{x_1,\ldots, x_n\in k^\ast}\left(\sum_{x_{n+1}\in k^\ast}
\chi_1^{-1}(x_{n+1})\psi(x_{n+1}f_1(x_1,\ldots, x_n))\right)\cdots
\left(\sum_{x_{n+m}\in
k^\ast}\chi_m^{-1}(x_{n+m})\psi(x_{n+m}f_m(x_1,\ldots,
x_n))\right)\\
&&\qquad\qquad\qquad\psi(f(x_1,\ldots, x_n)).
\end{eqnarray*}
For $i=1,\ldots, m$ and $x_1,\ldots, x_n\in k^\ast$, if
$f_i(x_1,\ldots,x_n)=0$, we have
$$\sum_{x_{n+i}\in k^\ast}
\chi_i^{-1}(x_{n+i})\psi(x_{n+i}f_i(x_1,\ldots, x_n))=0;$$ if
$f_i(x_1,\ldots,x_n)\not=0$, we have
\begin{eqnarray*}
\sum_{x_{n+i}\in k^\ast}
\chi_i^{-1}(x_{n+i})\psi(x_{n+i}f_i(x_1,\ldots, x_n))  &=&
\sum_{x\in k^\ast} \chi_i^{-1}\left(\frac{x}{f_i(x_1,\ldots,
x_n)}\right)\psi(x)\\
&=& \chi_i(f_i(x_1,\ldots, x_n))G(\chi_i,\psi),
\end{eqnarray*}
where
$$G(\chi_i,\psi)=\sum_{x\in k^\ast} \chi_i^{-1}(x)\psi(x)$$ is the
Gauss sum. So in any case, we have
$$\sum_{x_{n+i}\in k^\ast}
\chi_i^{-1}(x_{n+i})\psi(x_{n+i}f_i(x_1,\ldots,
x_n))=\chi_i(f_i(x_1,\ldots, x_n))G(\chi_i,\psi).$$ Hence
\begin{eqnarray*}
S_2&=&  \sum_{x_1,\ldots, x_n\in k^\ast}\chi_1(f_1(x_1,\ldots,
x_n))G(\chi_1,\psi)\cdots \chi_m(f_m(x_1,\ldots, x_n))G(\chi_m,\psi)
\psi(f(x_1,\ldots, x_n))\\
&=&G(\chi_1,\psi)\cdots G(\chi_m,\psi)S_1.
\end{eqnarray*}
As the Gauss sums are well-understood, the study of $S_1$ is reduced
to the study of $S_2$.

In this paper, we use $l$-adic cohomology theory to study sums of
the form
$$\sum_{{x_i}\in k^\ast}\chi_1(x_1)\cdots
\chi_n(x_n)\psi(f(x_1,\ldots, x_n)),$$ where $\chi_1,\ldots,
\chi_n$ are multiplicative characters (nontrivial or trivial).
Note that $S_2$ is of this form. Our results complete those in
[DL], where the case of trivial $\chi_i$ is treated. We follow the
approach initiated by Denef and Loeser.

\bigskip
We first associate geometric objects to the above data. The Kummer
covering
$$[q-1]:{\bf T}_k^n\to{\bf T}_k^n, \; x\mapsto x^{q-1}$$ on the torus
${\bf T}_k^n={\rm Spec}\,k[X_1,X_1^{-1},\ldots, X_n,X_n^{-1}]$
defines a ${\bf T}_k^n(k)$-torsor
$$1\to {\bf T}_k^n(k)\to {\bf T}_k^n\stackrel {[q-1]}\to {\bf T}_k^n\to 1,$$
where ${\bf T}_k^n(k)={\rm Hom}_k({\rm Spec}\, k, {\bf T}_k^n)$ is
the group of $k$-rational points in ${\bf T}_k^n$. Fix a prime
number $l$ distinct to $p$. Let $\chi: {\bf T}_k^n(k)={k^\ast}^n\to
\overline {\bf Q}_l^\ast$ be a character. Pushing-forward the above
torsor by $\chi^{-1}$, we get a lisse $\overline {\bf Q}_l$-sheaf
${\cal K}_\chi$ on ${\bf T}_k^n$ of rank 1. We call ${\cal K}_\chi$
the {\it Kummer sheaf} associated to $\chi$. For any rational point
$x$ in ${\bf T}_k^n(k')={\rm Hom}_k({\rm Spec}\, k', {\bf T}_k^n)$
with value in a finite extension $k'$ of $k$, we have
$${\rm Tr}(F_x, ({\cal K}_\chi)_{\bar x})=\chi({\rm Norm}_{k'/k}(x)),$$
where $F_x$ is the geometric Frobenius element at $x$.

The Artin-Schreier covering $${\cal P}:{\bf A}_k^1\to {\bf
A}_k^1,\; x\mapsto x^q-x$$ defines an ${\bf A}_k^1(k)$-torsor
$$0\to {\bf A}_k^1(k) \to {\bf A}_k^1\stackrel {\cal P}\to {\bf A}_k^1\to
0,$$ where ${\bf A}_k^1(k)={\rm Hom}_k({\rm Spec}\, k, {\bf A}_k^1)$
is the group of $k$-rational points in ${\bf A}_k^1$. Let $\psi:
{\bf A}_k^1(k)=k\to \overline {\bf Q}_l^\ast$ be an additive
character. Pushing-forward this torsor by $\psi^{-1}$, we get a
lisse ${\bf Q}_l$-sheaf ${\cal L}_\psi$ of rank $1$ on ${\bf
A}_k^1$, which we call the {\it Artin-Schreier sheaf}. For any
rational point $x$ in ${\bf A}_k^1(k')={\rm Hom}_k({\rm Spec}\, k',
{\bf A}_k^1)$ with value in a finite extension $k'$ of $k$, we have
$${\rm Tr}(F_x, ({\cal L}_\psi)_{\bar x})=\psi({\rm Tr}_{k'/k}(x)),$$
where $F_x$ is the geometric Frobenius element at $x$.

Let
$$f=\sum_{i\in {\bf Z}^n} a_i X^i\in k[X_1,X_1^{-1},\ldots, X_n,
X_n^{-1}]$$  be a Laurent polynomial. The {\it Newton polyhedron}
$\Delta_\infty(f)$ of $f$ at $\infty$ is the convex hull in ${\bf
R}^n$ of the set $\{i\in {\bf Z}^n|a_i\not= 0\}\cup \{0\}.$ We say
$f$ is {\it non-degenerate} with respect to $\Delta_\infty(f)$ if
for any face $\tau$ of $\Delta_\infty(f)$ not containing $0$, the
locus of
$$\frac{\partial f_\tau}{\partial X_1}=\cdots=
\frac{\partial f_\tau}{\partial X_n}=0$$ in ${\bf T}_k^n$ is
empty, where $$f_\tau=\sum_{i\in \tau} a_i X^i.$$ This is
equivalent to saying that the morphism
$$f_\tau:{\bf T}_k^n={\rm Spec}\,A[X_1,X_1^{-1},\ldots,
X_n, X_n^{-1}] \to {\bf A}_k^1={\rm Spec}\, k[T]$$ defined by the
$k$-algebra homomorphism $$k[T]\to k[X_1,X_1^{-1},\ldots, X_n,
X_n^{-1}],\; T\mapsto f_\tau$$ is smooth.

\bigskip

The first main result of this paper is the following theorem.

\bigskip
\noindent {\bf Theorem 0.1.} Let $f:{\bf T}_k^n\to {\bf A}_k^1$ be a
$k$-morphism defined by a Laurent polynomial $f\in
k[X_1,X_1^{-1},\ldots, X_n,X_n^{-1}]$ that is non-degenerate with
respect to $\Delta_\infty(f)$ and let ${\cal K}_\chi$ be a Kummer
sheaf on ${\bf T}_k^n$. Suppose ${\rm dim}(\Delta_\infty(f))=n$.
Then

(i) $H_c^i({\bf T}_{\bar k}^n, {\cal K}_\chi\otimes f^\ast{\cal
L}_\psi)=0$ for $i\not=n$.

(ii) ${\rm dim} (H_c^n ({\bf T}_{\bar k}^n, {\cal K}_\chi\otimes
f^\ast{\cal L}_\psi))=n! {\rm vol}(\Delta_\infty(f)).$

(iii) If $0$ is an interior point of $\Delta_\infty(f)$, then
$H_c^n({\bf T}_{\bar k}^n, {\cal K}_\chi\otimes f^\ast{\cal
L}_\psi)$ is pure of weight $n$.

\bigskip
Here the conclusion of (iii) means that for any eigenvalue $\lambda$
of the geometric Frobenius element  $F$ in ${\rm Gal}(\bar k/k)$
acting on $H_c^n({\bf T}_{\bar k}^n, {\cal K}_\chi\otimes
f^\ast{\cal L}_\psi)$, $\lambda$ is an algebraic number, and all the
galois conjugates of $\lambda$ have archimedean absolute value
$q^{\frac{n}{2}}.$

Note that we have
$$[q-1]_\ast [q-1]^\ast f^\ast {\cal L}_\psi
\cong \bigoplus_\chi ({\cal K}_\chi\otimes f^\ast {\cal
L}_\psi),$$ where $[q-1]:{\bf T}_k^n\to {\bf T}_k^n$ is the Kummer
covering, and in the direct sum on the right-hand side, $\chi$
goes over the set of all characters $\chi: {\bf T}_k^n(k)\to
\overline {\bf Q}_l$. The composition $f\circ [q-1]$ is defined by
the Laurent polynomial
$$f'(X_1,\ldots, X_n)=f(X_1^{q-1},\ldots, X_n^{q-1}).$$ Note that
$f'$ is also non-degenerate with respect to its Newton polyhedron
at $\infty$. We have
\begin{eqnarray*}
H_c^i({\bf T}_{\bar k}^n, {f'}^\ast {\cal L}_\psi)&\cong&
H_c^i({\bf T}_{\bar k}^n, [q-1]_\ast [q-1]^\ast f^\ast {\cal
L}_\psi) \\
&\cong& \bigoplus_\chi H_c^i({\bf T}_{\bar k}^n, {\cal
K}_\chi\otimes f^\ast {\cal L}_\psi).
\end{eqnarray*}
So $H_c^i({\bf T}_{\bar k}^n, {\cal K}_\chi\otimes f^\ast {\cal
L}_\psi)$ are direct factors of $H_c^i({\bf T}_{\bar k}^n,
{f'}^\ast {\cal L}_\psi)$. Hence Theorem 0.1 (i) and (iii) follow
directly from the main theorem 1.3 in [DL] applied to $f'$. Using
[I1] 2.1, one can show $$\chi_c({\bf T}_{\bar k}^n, {\cal
K}_\chi\otimes f^\ast{\cal L}_\psi)= \chi_c({\bf T}_{\bar k}^n,
f^\ast{\cal L}_\psi),$$ where $\chi_c$ denotes the Euler
characteristic for the cohomology with compact support. Hence
Theorem 0.1 (ii) can also be deduced from [DL] 1.3.

In this paper, we give a proof of Theorem 0.1 independent of the
main theorem of [DL]. On the other hand, our second main Theorem
0.4 on the weights of $H_c^n({\bf T}_{\bar k}^n, {\cal
K}_\chi\otimes f^\ast {\cal L}_\psi)$ doesn't seem to follow from
the corresponding theorem in [DL].

\bigskip
\noindent {\bf Corollary 0.2.} Let $f\in k[X_1,X_1^{-1},\ldots,
X_n,X_n^{-1}]$ be a Laurent polynomial that is non-degenerate with
respect to $\Delta_\infty(f)$ and suppose ${\rm
dim}(\Delta_\infty(f))=n$. Then for any multiplicative characters
$\chi_1,\ldots, \chi_n:k^\ast\to {\bf C}^\ast$ and any nontrivial
additive character $\psi:k\to {\bf C}^\ast$, we have
$$|\sum_{x_i\in k^\ast}
\chi_1(x_1)\cdots\chi_n(x_n)\psi(f(x_1,\ldots,x_n))|\leq n!{\rm
vol}(\Delta_\infty(f))q^{\frac{n}{2}}.$$

\bigskip
\noindent {\bf Proof.} Note that the values of $\chi_i$ and $\psi$
are algebraic integers. In particular, we may consider them to
have values in $\overline {\bf Q}_l$. Let $\chi: {k^\ast}^n\to
\overline {\bf Q}_l^\ast$ be the character defined by
$$\chi(x)=\chi_1(x_1)\ldots\chi_n(x_n)$$ for any $x=(x_1,\ldots,
x_n)\in {k^\ast}^n$. We then have
$$\sum_{x_i\in k^\ast}
\chi_1(x_1)\cdots\chi_n(x_n)\psi(f(x_1,\ldots,x_n))=\sum_{x\in
{\bf T}_k^n(k)}{\rm Tr}(F_x,({\cal K}_\chi\otimes f^\ast{\cal
L}_\psi)_{\bar x}).$$ By the Grothendieck trace formula ([SGA
$4\frac{1}{2}$], [Rapport] Th\'eor\`eme 3.2), we have
$$\sum_{x\in
{\bf T}_k^n(k)}{\rm Tr}(F_x,({\cal K}_\chi\otimes f^\ast{\cal
L}_\psi)_{\bar x})=\sum_{i=0}^{2n}(-1)^i {\rm Tr}(F,H_c^i({\bf
T}_{\bar k}^n, {\cal K}_\chi\otimes f^\ast{\cal L}_\psi)).$$ By
[D] 3.3.1, for any eigenvalue $\lambda$ of $F$ acting on
$H_c^n({\bf T}_{\bar k}^n, {\cal K}_\chi\otimes f^\ast{\cal
L}_\psi)$, $\lambda$ is an algebraic number and all its galois
conjugates have archimedean absolute value $\le q^{\frac{n}{2}}$.
Combined with Theorem 0.1 (i) and (ii), we get
$$|\sum_{i=0}^{2n} (-1)^i{\rm Tr}(F,H_c^i({\bf T}_{\bar k}^n,
{\cal K}_\chi\otimes f^\ast{\cal L}_\psi))|\le n!{\rm
vol}(\Delta_\infty(f))q^{\frac{n}{2}}.$$ So we have
$$|\sum_{x_i\in k^\ast}
\chi_1(x_1)\cdots\chi_n(x_n)\psi(f(x_1,\ldots,x_n))|\leq n!{\rm
vol}(\Delta_\infty(f))q^{\frac{n}{2}}.$$

\bigskip
Combining Corollary 0.2 with the discussion at the beginning, and
using the fact that Gauss sums have absolute value
$q^{\frac{1}{2}}$, we get the following.

\bigskip
\noindent {\bf Corollary 0.3.}  Let $f, f_1,\ldots, f_m \in
k[X_1,X_1^{-1},\ldots, X_n,X_n^{-1}]$ be Laurent polynomials,
$\chi_1,\ldots, \chi_m:k^\ast\to {\bf C}^\ast$ nontrivial
multiplicative characters, and $\psi:k\to {\bf C}^\ast$ a
nontrivial additive character. Suppose the Laurent polynomial
$$F(x_1,\ldots, x_n,x_{n+1},\ldots, x_{n+m})
=f(x_1,\ldots, x_n)+x_{n+1}f_1(x_1,\ldots, x_n) +\ldots + x_{n+m}
f_m(x_1,\ldots,x_n)$$ is non-degenerate with respect to
$\Delta_\infty(F)$, and ${\rm dim}(\Delta_\infty(F))=m+n$. Then we
have
$$|\sum_{x_i\in k^\ast}
\chi_1(f_1(x_1,\ldots, x_n))\cdots\chi_m(f_m(x_1,\ldots,
x_n))\psi(f(x_1,\ldots,x_n))|\leq (n+m)!{\rm
vol}(\Delta_\infty(F))q^{n/2}.$$

\bigskip
Under the assumption of Theorem 0.1, let
$$E({\bf T}_k^n, f,\chi)=\sum_{w\in {\bf Z}} e_wT^w,$$
where $e_w$ is the number of eigenvalues counted with multiplicities
of the geometric Frobenius element  $F$ in ${\rm Gal}(\bar k/k)$
acting on $H_c^n({\bf T}_{\bar k}^n, {\cal K}_\chi\otimes
f^\ast{\cal L}_\psi)$ with weight $w$. Our next goal is to determine
$E({\bf T}_k^n, f,\chi)$.

\bigskip
For any convex polyhedral cone $\delta$ in ${\bf R}^n$ with $0$
being a face, define the convex polytope ${\rm poly}(\delta)$ to be
the intersection of $\delta$ with a hyperplane in ${\bf R}^n$ which
does not contain $0$ and intersects each one dimensional face of
$\delta$. Note that ${\rm poly}(\delta)$ is defined only up to
combinatorial equivalence. For any convex polytope $\Delta$ in ${\bf
R}^n$ and any face $\tau$ of $\Delta$, define ${\rm
cone}_\Delta(\tau)$ to be the cone generated by $u'-u$
($u'\in\Delta$, $u\in \tau$), and define ${\rm
cone}_\Delta^\circ(\tau)$  to be the image of ${\rm
cone}_\Delta(\tau)$ in ${\bf R}^n/{\rm span}(\tau-\tau)$. Note that
$0$ is a face of ${\rm cone}_\Delta^\circ(\tau)$. We define
polynomials $\alpha(\delta)$ and $\beta(\Delta)$ in one variable $T$
inductively by the following formulas:
\begin{eqnarray*}
\alpha(\{0\})&=&1,\\
\beta(\Delta)&=&(T^2-1)^{{\rm dim}(\Delta)}+\sum_{\tau \hbox { face
of } \Delta,\;
\tau\not=\Delta}(T^2-1)^{{\rm dim}(\tau)}\alpha({\rm cone}_{\Delta}^\circ(\tau)),\\
\alpha(\delta)&=&{\rm trunc}_{\leq {\rm
dim}(\delta)-1}((1-T^2)\beta({\rm poly}(\delta))),
\end{eqnarray*}
where ${\rm trunc}_{\leq d}(\cdot)$ denotes taking the degree $\leq
d$ part of a polynomial. These polynomials are first introduced by
Stanley [S]. Note that $\alpha(\delta)$ and $\beta(\Delta)$ only
involve even powers of $T$, and they depend only on the
combinatorial types of $\delta$ and $\Delta$. If $\delta$ is a
simplicial cone, that is, if $\delta$ is generated by linearly
independent vectors, then using induction on ${\rm dim}(\delta)$,
one can verify $\alpha(\delta)=1$.

\bigskip
Let $\chi:{\bf T}_k^n(k)\to \overline {\bf Q}_l^\ast$ be a
character. For a rational convex polytope $\Delta$ in ${\bf R}^n$ of
dimension $n$, let $T$ be the set of faces $\tau$ of $\Delta$ so
that $\tau\not=\Delta$, $0\in\tau$, and ${\cal K}_\chi\cong
p_\tau^\ast{\cal K}_{\tau}$ for a Kummer sheaf ${\cal K}_{\tau}$ on
the torus ${\bf T}_\tau={\rm Spec}\, k[{\bf Z}^n\cap {\rm
span}(\tau-\tau)]$, where
$$p_\tau:{\bf T}_k^n={\rm Spec}\, k[{\bf Z}^n]\to {\bf T}_\tau=
{\rm Spec}\, k[{\bf Z}^n\cap {\rm span}(\tau-\tau)]$$ is the
morphism defined by the canonical homomorphism
$$k[{\bf Z}^n\cap {\rm span}(\tau-\tau)]\hookrightarrow  k[{\bf
Z}^n].$$ Define $$e(\Delta,\chi)= n!{\rm vol}(\Delta)+\sum_{\tau\in
T} (-1)^{n-{\rm dim}(\tau)}({\rm dim}(\tau))!{\rm
vol}(\tau)\alpha({\rm cone}_{\Delta}^\circ(\tau))(1)$$ and define a
polynomial $E(\Delta,\chi)$ inductively by
$$E(\Delta,\chi)=e(\Delta,\chi)T^n-\sum_{\tau\in T} (-1)^{n-{\rm
dim}(\tau)}E({\tau},\chi_\tau)\alpha({\rm
cone}_{\Delta}^\circ(\tau)).$$ Our second main result is the
following.

\bigskip
\noindent {\bf Theorem 0.4.} Let $f:{\bf T}_k^n\to {\bf A}_k^1$ be a
$k$-morphism defined by a Laurent polynomial $f\in
k[X_1,X_1^{-1},\ldots, X_n,X_n^{-1}]$ that is non-degenerate with
respect to $\Delta_\infty(f)$ and let ${\cal K}_\chi$ be a Kummer
sheaf on ${\bf T}_k^n$. Suppose ${\rm dim}(\Delta_\infty(f))=n$. Let
$$E({\bf T}_k^n, f,\chi)=\sum_{w\in {\bf Z}} e_wT^w,$$
where $e_w$ is the number of eigenvalues counted with
multiplicities of the geometric Frobenius element  $F$ in ${\rm
Gal}(\bar k/k)$ acting on $H_c^n({\bf T}_{\bar k}^n, {\cal
K}_\chi\otimes f^\ast{\cal L}_\psi)$ with weight $w$. Then $E({\bf
T}_k^n, f,\chi)$ is a polynomial of degree $\leq n$, and
\begin{eqnarray*}
E({\bf T}_k^n, f,\chi)&=&E(\Delta_\infty(f),\chi),\\
e_n&=&e(\Delta_\infty(f),\chi).
\end{eqnarray*}

\bigskip
\noindent {\bf Corollary 0.5.} Let $f:{\bf T}_k^n\to {\bf A}_k^1$
be a $k$-morphism defined by a Laurent polynomial $f\in
k[X_1,X_1^{-1},\ldots, X_n,X_n^{-1}]$ that is non-degenerate with
respect to $\Delta_\infty(f)$ and let ${\cal K}_\chi$ be a Kummer
sheaf on ${\bf T}_k^n$. Suppose ${\rm dim}(\Delta_\infty(f))=n$
and suppose for any face $\tau$ of $\Delta_\infty(f)$ of
codimension one containing $0$, ${\cal K}_\chi$ is not the inverse
image of any Kummer sheaf on ${\bf T}_\tau={\rm Spec}\, k[{\bf
Z}^n\cap {\rm span}(\tau-\tau)]$ under the morphism
$$p_\tau:{\bf T}_k^n={\rm Spec}\, k[{\bf Z}^n]\to {\bf T}_\tau=
{\rm Spec}\, k[{\bf Z}^n\cap {\rm span}(\tau-\tau)].$$ Then
$H_c^n({\bf T}_{\bar k}^n, {\cal K}_\chi\otimes f^\ast{\cal
L}_\psi)$ is pure of weight $n$.

\bigskip
\noindent {\bf Remark 0.6.} Note that this corollary together with
Theorem 0.1 is Theorem 4.2 in [AS], except that Adolphson and
Sperber prove the theorem for almost all $p$.

\bigskip
\noindent {\bf Proof of Corollary 0.5.} Under our assumption, the
set $T$ of faces $\tau$ of $\Delta_\infty(f)$ so that
$\tau\not=\Delta_\infty(f)$, $0\in\tau$, and ${\cal K}_\chi\cong
p_\tau^\ast{\cal K}_{\tau}$ for a Kummer sheaf ${\cal K}_{\tau}$ on
the torus ${\bf T}_\tau={\rm Spec}\, k[{\bf Z}^n\cap {\rm
span}(\tau-\tau)]$ is empty. Therefore
$$E(\Delta_\infty(f),\chi)=e(\Delta_\infty(f),\chi)T^n.$$
By Theorem 0.4, this implies $$E({\bf T}_k^n, f,\chi)=
e(\Delta_\infty(f),\chi)T^n,$$ and hence $H_c^n({\bf T}_{\bar
k}^n, {\cal K}_\chi\otimes f^\ast{\cal L}_\psi)$ is pure of weight
$n$.

\bigskip
For a rational convex polytope $\Delta$ in ${\bf R}^n$ with
dimension $n$, recall that $T$ is the set of faces $\tau$ of
$\Delta$ so that $\tau\not=\Delta$, $0\in\tau$, and ${\cal
K}_\chi\cong p_\tau^\ast{\cal K}_{\tau}$ for a Kummer sheaf ${\cal
K}_{\tau}$ on the torus ${\bf T}_\tau={\rm Spec}\, k[{\bf Z}^n\cap
{\rm span}(\tau-\tau)]$. Define
$$V_i(\Delta, \chi)=\sum_{\tau\in T, \;{\rm dim}(\tau)=i} {\rm
vol}(\tau)$$ for $0\leq i\leq n-1$ and define
$$V_n(\Delta,\chi)={\rm vol}(\Delta).$$
Let $\tau_0$ be the smallest face of $\Delta$ containing $0$. We
say $\Delta$ is {\it simplicial at the origin} if there are
exactly $n-{\rm dim}(\tau_0)$ faces of $\Delta$ that have
dimension $n-1$ and contain $\tau_0$. If $\Delta$ is simplicial at
the origin and $\tau$ is a face containing $0$, then the number of
faces of $\Delta$ that have dimension $k$ and contain $\tau$ is
$\left(
\begin{array}{c}
n-{\rm dim}(\tau)\\
n-k
\end{array}
\right).$

The following corollary is Theorem 4.8 in [AS].

\bigskip
\noindent {\bf Corollary 0.7.} Let $f:{\bf T}_k^n\to {\bf A}_k^1$
be a $k$-morphism defined by a Laurent polynomial $f\in
k[X_1,X_1^{-1},\ldots, X_n,X_n^{-1}]$ that is non-degenerate with
respect to $\Delta_\infty(f)$ and let ${\cal K}_\chi$ be a Kummer
sheaf on ${\bf T}_k^n$. Suppose ${\rm dim}(\Delta_\infty(f))=n$
and suppose $\Delta_\infty(f)$ is simplicial at the origin. Let
$e_w$ be the number of eigenvalues counted with multiplicities of
the geometric Frobenius element $F$ in ${\rm Gal}(\bar k/k)$
acting on $H_c^n({\bf T}_{\bar k}^n, {\cal K}_\chi\otimes
f^\ast{\cal L}_\psi)$ with weight $w$. Then we have
$$e_w=\sum_{i=0}^w (-1)^{w-i}i! \left(
\begin{array}{c}
n-i\\
n-w
\end{array}
\right)V_i(\Delta_\infty(f),\chi).$$

\bigskip
\noindent {\bf Proof.} The hypothesis that $\Delta_\infty(f)$ is
simplicial at the origin implies that $$\alpha({\rm
cone}_{\Delta_\infty(f)}^\circ(\tau))=1$$ for any face $\tau$ of
$\Delta_\infty(f)$ containing $0$. By Theorem 0.4, we have
\begin{eqnarray*}
e_n&=& e(\Delta_\infty(f),\chi)\\
&=&  n!{\rm vol}(\Delta_\infty(f))+\sum_{\tau\in T} (-1)^{n-{\rm
dim}(\tau)}({\rm dim}(\tau))!{\rm vol}(\tau) \\
&=& n!V_n(\Delta_\infty(f),\chi)+\sum_{i=0}^{n-1} (-1)^{n-i} i!
V_i(\Delta_\infty(f),\chi)\\
&=& \sum_{i=0}^{n} (-1)^{n-i} i! V_i(\Delta_\infty(f),\chi).
\end{eqnarray*}
This proves our assertion for $w=n$. We use induction on $n$.
Under our assumption, we have
$$E(\Delta_\infty(f),\chi)=e(\Delta_\infty(f),\chi)T^n-\sum_{\tau\in T} (-1)^{n-{\rm
dim}(\tau)}E({\tau},\chi_\tau).$$ For $w\leq n-1$, taking the
coefficients of $T^w$ on both sides of the above equality and
applying Theorem 0.4 and the induction hypothesis to the pairs
$({\tau},\chi_\tau)$ $(\tau\in T)$, we get
$$e_w=-\sum_{\tau\in T,\; w\leq {\rm dim}(\tau)\leq n-1} (-1)^{n-{\rm
dim}(\tau)}\sum_{i=0}^w (-1)^{w-i}i! \left(
\begin{array}{c} {\rm
dim}(\tau)-i\\ {\rm dim}(\tau)-w
\end{array}
\right)\sum_{\tau'\prec \tau,\; {\rm dim}(\tau')=i,\;\tau'\in T}
{\rm vol}(\tau').$$ So we have
$$e_w=\sum_{i=0}^w(-1)^{w-i+n+1}i!\sum_{\tau'\in T,\; {\rm
dim}(\tau')=i}\; \sum_{\tau'\prec \tau,\; w\leq {\rm
dim}(\tau)\leq n-1} (-1)^{{\rm dim}(\tau)}\left(\begin{array}{c}
{\rm dim}(\tau)-i\\ {\rm dim}(\tau)-w\end{array} \right){\rm
vol}(\tau').$$ By our assumption, for each $\tau'$ containing $0$
of dimension $i$, the number of faces of $\Delta_\infty(f)$ that
have dimension $k$ and contain $\tau'$ is $\left(
\begin{array}{c}
n-i\\
n-k
\end{array}
\right).$ So we have
$$e_w=\sum_{i=0}^w(-1)^{w-i+n+1}i!\sum_{\tau'\in T,\; {\rm
dim}(\tau')=i} \sum_{k=w}^{n-1} (-1)^k \left(\begin{array}{c} k-i\\
k-w\end{array} \right) \left(
\begin{array}{c}
n-i\\
n-k
\end{array}
\right){\rm vol}(\tau').$$ We have
\begin{eqnarray*}
&&\sum_{k=w}^{n-1} (-1)^k \left(\begin{array}{c} k-i\\
k-w\end{array} \right) \left(
\begin{array}{c}
n-i\\
n-k
\end{array}
\right)\\
&=& \sum_{k=w}^{n-1} (-1)^k \left(\begin{array}{c} n-w\\
k-w\end{array} \right) \left(
\begin{array}{c}
n-i\\
n-w
\end{array}
\right)\\
&=&(-1)^w \left(
\begin{array}{c}
n-i\\
n-w
\end{array}
\right)\sum_{j=0}^{n-1-w} (-1)^j \left(\begin{array}{c} n-w\\
j\end{array} \right)\\
&=&(-1)^{n-1}\left(
\begin{array}{c}
n-i\\
n-w
\end{array}
\right).
\end{eqnarray*}
So we have
\begin{eqnarray*}
e_w&=&\sum_{i=0}^w(-1)^{w-i+n+1}i!\sum_{\tau'\in T,\; {\rm
dim}(\tau')=i}(-1)^{n-1}\left(
\begin{array}{c}
n-i\\
n-w
\end{array}
\right){\rm vol}(\tau').\\
&=&\sum_{i=0}^w (-1)^{w-i}i! \left(
\begin{array}{c}
n-i\\
n-w
\end{array}
\right)V_i(\Delta_\infty(f),\chi).
\end{eqnarray*}
This proves our assertion.

\bigskip
The paper is organized as follows. In \S 1, we summarize basic
results on toric schemes. These results can be found in [F]. This
section is mainly for the purpose of fixing notations. A toric
scheme contains an open dense torus. In \S 2, we study extensions of
Kummer sheaves on tori to toric schemes. In \S 3, we study the
compactifications by toric schemes of a morphism on a torus defined
by a Laurent polynomial. In \S 4, we prove our main results.

\bigskip
\noindent{\bf Acknowledgement.} The research is supported by the
NSFC (10525107).

\bigskip
\bigskip
\centerline {\bf 1. Toric Schemes}

\bigskip
\bigskip
In this paper, a {\it lattice} $N$ is a free abelian group of finite
rank. Let $M={\rm Hom}_{\bf Z}(N,{\bf Z})$ be the dual lattice of
$N$, let $V=N\otimes_{\bf Z} {\bf R}$ be the real vector space
generated by $N$, and let $V^\ast={\rm Hom}_{\bf R}(V,{\bf R})$ be
the dual vector space of $V$. We have a canonical identification
$M\otimes_{\bf Z} {\bf R}\cong V^\ast$. A {\it convex polyhedral
cone} $\sigma$ in $V$ is a subset of the form
$$\sigma=\{r_1v_1+\cdots +r_k v_k | r_i\geq 0\},$$
where $v_1,\ldots, v_k$ is a finite family of elements in $V$,
which is called a family of {\it generators} of $\sigma$. We say
$\sigma$ is {\it rational} with respect to the lattice $N$ if
$v_1,\ldots, v_k$ can be chosen to lie in $N$. Define the {\it
dual} $\check\sigma$ of $\sigma$ to be
$$\check \sigma =\{u\in V^\ast |\langle u, v\rangle \geq 0 \hbox {
for all } v\in \sigma\}.$$ If $\sigma$ is a rational convex
polyhedral cone in $V$ with respect to $N$, then $\check\sigma$ is
a rational convex polyhedral cone in $V^\ast$ with respect to $M$
([F] 1.2 (9)). Under this condition, the semigroup $M\cap \check
\sigma$ is finitely generated ([F] 1.2, Proposition 1). Let $A$ be
a commutative ring. Then the ring $A[M\cap \check\sigma]$ is a
finitely generated $A$-algebra. For any $u\in M\cap\check \sigma$,
denote the corresponding element in $A[M\cap \check \sigma]$ by
$\chi^u$. We have
$$\chi^{u_1+u_2}=\chi^{u_1}\chi^{u_2}$$
in $A[M\cap \check \sigma]$ for any $u_1,u_2\in M\cap \check
\sigma$. A {\it face} $\tau$ of $\sigma$ is a subset of the form
$$\tau=\sigma \cap u^\perp=\{v\in \sigma| \langle u,v\rangle =0\}$$
for some $u\in \check \sigma$. It is also a convex polyhedral cone
([F] 1.2 (2)). We use the notation $\tau\prec \sigma$ or
$\sigma\succ\tau$ to denote $\tau$ being a face of $\sigma$. When
$\sigma$ is rational, so is $\tau$, and we may then choose $u\in
M\cap \check\sigma$. We then have
$$M\cap \check \tau=M\cap \check \sigma + {\bf Z}_{\geq 0} (-u)$$
by [F] 1.2, Proposition 2, and hence $$A[M\cap \check \tau]=A[M\cap
\check \sigma]_{\chi^u}.$$ Define
$$U_\sigma={\rm Spec}\, A[M\cap \check \sigma].$$
Then the canonical homomorphism $$A[M\cap \check\sigma]
\hookrightarrow A[M\cap\check \tau]$$ defines an open immersion
$$U_\tau\hookrightarrow U_\sigma.$$
A {\it fan} $\Sigma$ is a finite family of rational convex
polyhedral cones in $V$ satisfying the following properties:

(a) $0$ is face of each cone in $\Sigma$.

(b) A face of a cone in $\Sigma$ is also in $\Sigma$.

(c) If $\sigma,\sigma'\in \Sigma$, then $\sigma\cap \sigma'$ is a
face of both $\sigma$ and $\sigma'$.

\noindent For any $\sigma,\sigma'\in \Sigma$, by the discussion
above, $U_{\sigma\cap\sigma'}$ can be considered as an open
subscheme of both $U_\sigma$ and $U_{\sigma'}$. Gluing
$U_{\sigma}$ and $U_{\sigma'}$ along $U_{\sigma\cap\sigma'}$, we
get a scheme $X_A(\Sigma)$ of finite type over $A$, which we call
the {\it toric scheme} associated to the fan $\Sigma$. The scheme
$X_A(\Sigma)$ is separated over $A$ ([F] 1.4, Lemma). If
$\bigcup_{\sigma\in \Sigma} \sigma=V$, then it is proper over $A$
([F] 2.4, Proposition). A rational convex polyhedral cone $\sigma$
is called {\it regular} if it can be generated by part of a basis
of $N$. Under this condition, $U_\sigma$ is smooth over $A$ ([F]
2.1). A fan $\Sigma$ is called {\it regular} if all the cones in
$\Sigma$ are regular. Under this condition, $X_A(\Sigma)$ is
smooth over $A$.

Each $U_\sigma$ $(\sigma\in \Sigma)$ can be regarded as an open
subscheme of $X_A(\Sigma)$. Taking $\sigma=0$, we see the torus
$$U_0={\rm Spec}\, A[M]$$
is an open subscheme of $X_A(\Sigma)$. One can show this torus is
dense in $X_A(\Sigma)$. Let $n$ be the rank of $N$. Then the torus
${\rm Spec}\, A[M]$ has relative dimension $n$, and we denote it
by ${\bf T}_A^n$. For any $\sigma\in \Sigma$, the $A$-algebra
homomorphism
$$A[M\cap \check\sigma]\to A[M]\otimes_A A[M\cap \check \sigma],\;
\chi^u\mapsto \chi^u\otimes \chi^u \; (u\in M\cap \check \sigma)$$
defines an action
$${\bf T}_A^n\times _A U_\sigma \to U_\sigma.$$ These actions for
$\sigma\in \Sigma$ can be glued together to give an action
$${\bf T}_A^n\times_A X_A(\Sigma)\to X_A(\Sigma)$$ which extends the
action of ${\bf T}_A^n$ on itself.

For any $\tau\in \Sigma$, let
$$N_\tau=N\cap \tau + (-N\cap \tau)=N\cap {\rm span}(\tau)$$ be
the group generated by $N\cap \tau$, and let $N(\tau)=N/N_\tau$.
Then $N(\tau)$ is torsion free, and hence a lattice. Let $M(\tau)$
be its dual lattice. We have a canonical isomorphism
$$M(\tau)\cong M\cap \tau^\perp.$$ For each $\sigma\in \Sigma$ with
$\tau\prec \sigma$, let $\bar \sigma$ be the image of $\sigma$ in
$N(\tau)\otimes_{\bf Z} {\bf R}\cong V/{\rm span}(\tau)$. Note
that $\sigma$ is completely determined by $\bar \sigma$, that is,
for any $\sigma,\sigma'\in \Sigma$ containing $\tau$, we have
$\sigma=\sigma'$ if and only if $\bar\sigma=\overline {\sigma'}$.
The family
$${\rm star}(\tau)=\{\bar\sigma|\sigma\in\Sigma,\;
\tau\prec\sigma\}$$ is a fan in $N(\tau)\otimes_{\bf Z}{\bf R}$. Let
$$V(\tau)=X_A({\rm star}(\tau)).$$ It is obtained by gluing
$$U_\sigma(\tau)={\rm Spec}\,A[M(\tau)\cap
\check{\bar\sigma}]={\rm Spec}\,A[M\cap\check \sigma \cap
\tau^\perp]$$ together ($\sigma\in {\rm Star}(\tau)$). Let
$$O_\tau= U_\tau(\tau)={\rm Spec}\; A[M\cap\tau^\perp]$$ be the
open dense torus in $V(\tau)$. Its relative dimension is ${\rm
dim}(V)-{\rm dim}(\tau).$ For any $\sigma\in\Sigma$ with
$\tau\prec\sigma$, consider the map
\begin{eqnarray*}
M\cap\check\sigma&\to& A[M\cap\check\sigma\cap\tau^\perp], \cr
u&\mapsto&
\left\{
\begin{array}{cc}
\chi^u& \hbox { if } u\in\check\sigma \cap\tau^\perp,\cr 0& \hbox
{ if } u\not \in\check\sigma \cap\tau^\perp.
\end{array}\right.
\end{eqnarray*}
It is a semigroup homomorphism, where the semigroup law on
$A[M\cap\check\sigma\cap\tau^\perp]$ is multiplication. (To prove
this, we use the fact that $\check\sigma\cap\tau^\perp$ is a face
of $\check\sigma$ ([F] 1.2 (10)), and a sum of elements in a cone
lies in a face if and only if each summand lies in the face.) It
induces an $A$-algebra epimorphism
$$A[M\cap\check \sigma]\to A[M\cap\check \sigma\cap \tau^\perp]$$
and hence a closed immersion
$$U_\sigma(\tau)\to U_\sigma.$$ For any $\sigma,\sigma'\in \Sigma$ with
$\tau\prec\sigma'\prec\sigma$, the diagram
$$\begin{array}{ccc}
U_{\sigma'}(\tau)&\to& U_{\sigma'}\\
\downarrow&&\downarrow\\
U_\sigma(\tau)&\to& U_\sigma
\end{array}$$
commutes and is Cartesian. So we can glue these closed immersions
together to get a closed immersion
$$V(\tau)\to \bigcup_{\tau\prec\sigma} U_\sigma.$$ One can show
for those $\sigma\in \Sigma$ not containing $\tau$, the image of
$V(\tau)$ in $\bigcup_{\tau\prec\sigma} U_\sigma$ is disjoint from
$U_\sigma$. So composing the above closed immersion with the open
immersion $\bigcup_{\tau\prec\sigma} U_\sigma\to X_A(\Sigma)$, we
get a closed immersion
$$V(\tau)\to X_A(\Sigma).$$
We regard the open dense torus $O_\tau$ of $V(\tau)$ as a subscheme
of $X_A(\Sigma)$ through this closed immersion. One can verify
$O_\tau$ $(\tau\in\Sigma)$ are disjoint in $X_A(\Sigma)$. Actually
$O_\tau$ are the orbits of the torus action ${\bf T}_A^n$ on
$X_A(\Sigma)$. Moreover, we have ([F] 3.1, Proposition)
\begin{eqnarray*}
U_\sigma&=&\coprod_{\gamma\prec\sigma} O_\gamma,\\
V(\tau)&=& \coprod_{\tau\prec\gamma} O_\gamma.
\end{eqnarray*}

\bigskip
Let $\delta$ be a rational convex polyhedral cone of dimension
${\rm dim}(V)$ in the dual space $V^\ast$. We have
$$\check\delta\cap(-\check\delta)=\delta^\perp=0.$$
So by [F] 1.2 (10), $0$ is a face of $\check\delta$, and hence the
family $\Sigma(\delta)$ of faces of $\check \delta$ is a fan in $V$,
and we have
$$X_A(\Sigma(\delta))=U_{\check \delta}={\rm Spec}\,A[M\cap
\delta].$$ By [F] 1.2 (10), the map $\tau\mapsto \check
\delta\cap\tau^\perp$ sets up a one-to-one correspondence between
the family of faces of $\delta$ and the family of faces of
$\check\delta$. We claim
$$\check\delta\cap \tau^\perp=({\rm cone}_\delta(\tau))^\vee,$$
where ${\rm cone}_\delta(\tau)$ is the cone in $V^\ast$ generated
by $u'-u$ $(u'\in \delta,\; u\in \tau)$. Indeed, if $v\in ({\rm
cone}_\delta(\tau))^\vee$, then for any $u'\in\delta, \; u\in
\tau$, we have
$$\langle u',v\rangle\geq \langle u,v\rangle.$$ Taking
$u,u'\in\tau$, we see $\langle\cdot,v\rangle |_\tau$ is constant.
As $0\in \tau$, we have $\langle\cdot,v\rangle |_\tau=0$, that is,
$v\in\tau^\perp$. The above inequality then implies that $\langle
u',v\rangle\geq 0$ for all $u'\in \delta$, that is, $v\in
\check\delta$. So we have $v\in \check\delta\cap \tau^\perp$.
Hence $({\rm cone}_\delta(\tau))^\vee\subset\check \delta\cap
\tau^\perp$. It is not hard to see $\check\delta\cap
\tau^\perp\subset({\rm cone}_\delta(\tau))^\vee.$ So
$\check\delta\cap \tau^\perp=({\rm cone}_\delta(\tau))^\vee.$ By
[F] 1.2 (10), we have
\begin{eqnarray*}
{\rm dim}\left(({\rm cone}_\delta(\tau))\bigcap (-{\rm
cone}_\delta(\tau))\right)&=&{\rm dim}(V)-{\rm dim}(({\rm
cone}_\delta(\tau))^{\vee})\\&=&{\rm dim}(V)-{\rm
dim}(\check\delta\cap\tau^\perp)\\&=&{\rm dim}(\tau).
\end{eqnarray*}
As
$${\rm span}(\tau)=\tau-\tau\subset ({\rm cone}_\delta(\tau))\bigcap (-{\rm
cone}_\delta(\tau))$$ and $({\rm cone}_\delta(\tau))\bigcap (-{\rm
cone}_\delta(\tau))$ is a linear space, we have
$$({\rm cone}_\delta(\tau))\bigcap
(-{\rm cone}_\delta(\tau))=\tau-\tau.$$ We summarize these results
as follows.

\bigskip
\noindent {\bf Proposition 1.1.} Let $\delta$ be a rational convex
polyhedral cone of dimension ${\rm dim}(V)$ in $V^\ast$. For any
face $\tau$ of $\delta$, let ${\rm cone}_{\delta}(\tau)$ be the
cone in $V^\ast$ generated by $u'-u$ $(u'\in \delta, u\in \tau)$.
We have
$$({\rm cone}_\delta(\tau))^\vee=\check\delta\cap \tau^\perp.$$
The family
$$\Sigma(\delta)=\{({\rm cone}_\delta(\tau))^\vee| \tau\prec
\delta\}$$ is a fan in $V=N\otimes_{\bf Z} {\bf R}$, and coincides
with the fan consisting of faces of $\check\delta$.  We have
$$X_A(\Sigma(\delta))=U_{\check \delta}={\rm Spec}\, A[M\cap
\delta].$$ Moreover, we have
$${\rm dim}({\rm cone}_\delta(\tau))^\vee={\rm dim}(V)-{\rm
dim}(\tau)$$ and
$$(-{\rm cone}_\delta(\tau))\bigcap({\rm cone}_\delta(\tau))=\tau-\tau.$$

\bigskip
A convex {\it polytope} $P$ in $V=N\otimes_{\bf Z}{\bf R}$ is a
subset of $V$ which can be written as the convex hull of a finite
family of element in $V$. If this finite family can be chosen to lie
in $N\otimes_{\bf Z} {\bf Q}$, we say $P$ is {\it rational}. If $P$
is a rational convex polytope in $V$ such that $0$ is an interior
point of $P$, then the set $\Sigma$ of cones over faces on the
boundary of $P$ is a fan, and the toric scheme $X_A(\Sigma)$ is
proper over $A$.

Let $\Delta$ be a rational convex polytope in the dual space
$V^\ast=M\otimes_{\bf Z}{\bf R}$. First consider the case where
$0$ is an interior point of $\Delta$. A {\it face} $\tau$ of
$\Delta$ is a subset of the form
$$\tau=\{u\in \Delta |\langle u, v\rangle =r\},$$
where $r$ is a real number and $v\in V$ is a vector such that
$\langle u, v\rangle \geq r$ for all $u\in \Delta$. Since $0$ is
an interior point of $\Delta$, we have either $v=0$ or $r<0$. When
$v=0$, the face $\tau$ is just $\Delta$. When $r<0$, we can always
choose $v$ so that $r=-1$. So a proper face of $\Delta$ is of the
form
$$\tau =\{u\in \Delta |\langle u, v\rangle =-1\},$$
where $v\in V$ is a vector such that $\langle u, v\rangle \geq -1$
for all $u\in \Delta$.

Define the {\it polar set} $\Delta^\circ$ of $\Delta$ to be
$$\Delta ^\circ =\{v\in V |\langle u,v\rangle \geq -1 \hbox { for
all } u\in \Delta\}.$$ It is a rational convex polytope in $V$
([F] 1.5, Proposition), and $0$ is in its interior. Denote the fan
in $V$ of cones over faces on the boundary of $\Delta^\circ$ by
$\Sigma(\Delta)$.

For any face $\tau$ of $\Delta$, let
$$\tau^\ast=\{v\in \Delta^\circ | \langle u,v\rangle =-1 \hbox {
for all } u\in \tau\}.$$ Then ([F] 1.5, Proposition) $\tau\mapsto
\tau^\ast$ sets up a one-to-one correspondence between faces of
$\Delta$ and faces of $\Delta^\circ$, and
$${\rm dim} (\tau) + {\rm dim} (\tau^\ast) ={\rm dim} (V)-1.$$
The cone $C(\tau^\ast)$ in $V$ over $\tau^\ast$ is
\begin{eqnarray*}
C(\tau^\ast)&=&\{tv| t\geq 0, v\in V, \langle u, v\rangle \geq -1
\hbox { for all } u\in \Delta, \langle u, v\rangle= -1 \hbox { for
all }
u\in \tau\}\\
&=& \{v|v\in V, \langle u',v\rangle \geq \langle u,v\rangle \hbox {
for
all } u'\in \Delta, u\in \tau\}\\
&=& ({\rm cone}_\Delta(\tau))^{\vee},
\end{eqnarray*}
where ${\rm cone}_\Delta(\tau)$ is the cone in $V^\ast$ generated
by $u'-u$ ($u'\in \Delta, u\in \tau$). So the fan $\Sigma(\Delta)$
is the set
$$\Sigma(\Delta)=\{({\rm cone}_\Delta(\tau))^{\vee}|\tau\prec \Delta\}.$$
Note that we have
\begin{eqnarray*}
{\rm dim}(({\rm cone}_\Delta(\tau))^{\vee})&=&{\rm dim}
(C(\tau^\ast))\\
&=& {\rm dim}(\tau^\ast) +1\\
&=& {\rm dim}(V)-{\rm dim}(\tau).
\end{eqnarray*}
By [F] 1.2 (10), we have
$${\rm dim}\left(({\rm cone}_\Delta(\tau))\bigcap (-{\rm
cone}_\Delta(\tau))\right)={\rm dim}(V)-{\rm dim}(({\rm
cone}_\Delta(\tau))^{\vee})={\rm dim}(\tau).$$ As
$$\tau-\tau\subset ({\rm cone}_\Delta(\tau))\bigcap (-{\rm
cone}_\Delta(\tau)),$$ we have
$$({\rm cone}_\Delta(\tau))\bigcap
(-{\rm cone}_\Delta(\tau))={\rm span}(\tau-\tau).$$

For any $v\in V$, the {\it first meeting locus} $F_{\Delta}(v)$ is
defined to be the subset of $\Delta$ where the function
$\langle\cdot,v\rangle |_\Delta$ reaches its minimum. If $r$ is
the minimum of $\langle\cdot,v\rangle |_\Delta$, then $\langle
u,v\rangle\geq r$ for all $u\in \Delta$, and $$F_\Delta(v)=\{u\in
\Delta | \langle u,v\rangle=r\}.$$ So $F_\Delta(v)$ is a
(nonempty) face of $\Delta$. For any subset $A$ in $V$, let
$$F_\Delta(A)=\bigcap_{v\in A} F_\Delta(v).$$
It is also face of $\Delta$. (Maybe empty). We have
\begin{eqnarray*}
C(\tau^\ast)&=& \{v|\langle u',v\rangle \geq \langle u,v\rangle
\hbox { for
all } u'\in \Delta, u\in \tau\}\\
&=& \{v|F_\Delta(v)\supset \tau\}.
\end{eqnarray*}
As $\tau^\ast$ is the closure of the complement of its proper
faces, $C(\tau^\ast)$ is the closure of the set
\begin{eqnarray*}
&&C(\tau^\ast)-\bigcup_{{\tau'}^\ast\prec
\tau^\ast,{\tau'}^\ast\not= \tau^\ast} C({\tau'}^\ast) \\
&=& C(\tau^\ast)-\bigcup_{\tau\prec
\tau',\tau\not=\tau'} C({\tau'}^\ast)\\
&=& \{v|F_\Delta(v)\supset \tau\}-\bigcup_{\tau\prec
\tau',\tau\not=\tau'} \{v|F_\Delta(v)\supset
\tau'\}\\
&=&\{v|F_\Delta(v)=\tau\}.
\end{eqnarray*}
Moreover, we have
$$\tau=\bigcap_{v\in C(\tau^\ast)}
F_\Delta(v)=F_\Delta(C(\tau^\ast)).$$ We summarize these results
as follows.

\bigskip
\noindent {\bf Proposition 1.2.} Let $\Delta$ be a rational convex
polytope of dimension ${\rm dim}(V)$ in $V^\ast$. For any face
$\tau$ of $\Delta$, let ${\rm cone}_{\Delta}(\tau)$ be the cone
generated by $u'-u$ $(u'\in \Delta, u\in \tau)$. Then the family
$$\Sigma(\Delta)=\{({\rm cone}_\Delta(\tau))^\vee| \tau\prec
\Delta\}$$ is a fan in $V=N\otimes_{\bf Z} {\bf R}$, and
$X_A(\Sigma(\Delta))$ is proper over $A$. For any $v\in V$, let
$F_\Delta(v)$ be the face of $\Delta$ where the function
$\langle\cdot, v\rangle|_\Delta$ reaches minimum. Then
\begin{eqnarray*}
({\rm cone}_\Delta(\tau))^\vee&=&\{v|F_\Delta(v)\supset \tau\}\\
&=&\overline {\{v|F_\Delta(v)= \tau\}}
\end{eqnarray*}
The map $$\tau\mapsto ({\rm cone}_\Delta(\tau))^\vee$$ sets up a
one-to-one correspondence between faces of $\Delta$ and cones in
the fan $\Sigma(\Delta)$, and the inverse map is
$$\sigma\mapsto F_\Delta(\sigma)=\bigcap _{v\in \sigma}
F_\Delta(v).$$ Moreover, we have
$${\rm dim}({\rm cone}_\Delta(\tau))^\vee={\rm dim}(V)-{\rm
dim}(\tau)$$ and
$$({\rm cone}_\Delta(\tau))\bigcap(-{\rm cone}_\Delta(\tau))={\rm
span}(\tau-\tau).$$

\bigskip
Note that Proposition 1.2 holds for the polytope $\Delta$ if and
only if it holds for a translation of $\Delta$. Making a
translation, we may assume $0$ is an interior point of $\Delta$.
In this case, Proposition 1.2 follows from the previous
discussion.

\bigskip
\bigskip
\centerline {\bf 2. Extensions of Kummer sheaves to toric schemes}

\bigskip
\bigskip
In this section $k$ is a field. Let $m$ be a positive integer
prime to the characteristic of $k$ so that $k$ contains a
primitive $m$-th root of unity. Let
$$\mu_m(k)^n=\{(\zeta_1,\ldots, \zeta_n)|\zeta_i\in k,
\zeta_i^m=1\}.$$ The Kummer covering $$[m]:{\bf T}_k^n\to{\bf
T}_k^n, \; x\mapsto x^m$$ on the torus ${\bf T}_k^n$ defines a
$\mu_m(k)^n$-torsor
$$1\to \mu_m(k)^n\to {\bf T}_k^n\stackrel {[m]}\to {\bf T}_k^n\to 1.$$
Let $\chi: \mu_m(k)^n\to \overline {\bf Q}_l^\ast$ be a character.
Pushing-forward the above torsor by $\chi^{-1}$, we get a lisse
$\overline {\bf Q}_l$-sheaf ${\cal K}_\chi$ on ${\bf T}_k^n$ of
rank 1. We call ${\cal K}_\chi$ the {\it Kummer sheaf} associated
to $\chi$. The following properties of the Kummer sheaf are
standard. We omit their proof.

(a) We have an isomorphism $${\cal K}_\chi|_{\{1\}}\cong \overline
{\bf Q}_l,$$ where $1$ is the identity of the algebraic group ${\bf
T}_k^n$.

(b) Let $s: {\bf T}_k^n\times_k {\bf T}_k^n\to {\bf T}_k^n$ be the
multiplication on the torus, and let ${\rm pr}_1, {\rm pr}_2: {\bf
T}_k^n\times_k {\bf T}_k^n\to {\bf T}_k^n$ be the two projections.
We have an isomorphism $$s^\ast {\cal K}_\chi\cong {\rm pr}_1^\ast
{\cal K}_\chi\otimes {\rm pr}_2^\ast {\cal K}_\chi$$ which is
compatible with the isomorphism
$$(s^\ast {\cal K}_\chi)|_{\{(1,1)\}}\cong ({\rm pr}_1^\ast {\cal
K}_\chi\otimes {\rm pr}_2^\ast {\cal K}_\chi)|_{\{(1,1)\}}$$ defined
by (a). Taking the inverse image of the above isomorphism under the
morphism $${\bf T}_k^n\to {\bf T}_k^n\times_k{\bf T}_k^n,\; x\mapsto
(x,x^{-1}),$$ we deduce an isomorphism
$${\cal K}_{{\chi}^{-1}}\cong {\cal K}_\chi^\vee.$$

(c) Suppose $m'|m$ and suppose there exists a character $\chi':
\mu_{m'}(k)^n\to \overline {\bf Q}_l^\ast$ such that
$$\chi(\zeta)=\chi'(\zeta^{\frac{m}{m'}}).$$ Then we have an
isomorphism
$${\cal K}_\chi\cong {\cal K}_{\chi'}.$$

(d) Let $n_i$ $(i=1,2)$ be positive integers such that
$n=n_1+n_2$, let $\chi_i: \mu_m(k)^{n_i}\to \overline {\bf
Q}_l^\ast$ be characters such that
$$\chi(\zeta_1,\ldots, \zeta_n)=\chi_1(\zeta_1,\ldots,
\zeta_{n_1})\chi_2(\zeta_{n_1+1},\ldots, \zeta_n),$$ and let
$$p_i: {\bf T}_k^n={\bf T}_k^{n_1}\times_k {\bf T}_k^{n_2}\to {\bf
T}_k^{n_i}$$ be the projections. Then we have an isomorphism
$${\cal K}_\chi\cong p_1^\ast {\cal K}_{\chi_1}\otimes p_2^\ast
{\cal K}_{\chi_2}.$$

\bigskip
Now let $N$ be a lattice of rank $n$, let $M={\rm Hom}_{\bf
Z}(N,{\bf Z})$ be its dual lattice, let $V=N\otimes_{\bf Z}{\bf R}$,
and let $V^\ast={\rm Hom}_R(V, R)\cong M\otimes_{\bf Z}{\bf R}$. For
any fan $\Sigma$ in $V$, denote the open immersion of the torus
${\bf T}_k^n={\rm Spec}\,k [M]$ in the toric $k$-scheme
$X_k(\Sigma)$ by $$j: {\bf T}_k^n\hookrightarrow X_k(\Sigma).$$ In
this section, we study the sheaf $j_\ast{\cal K}_\chi$ and the
perverse sheaf $j_{!\ast}(K_{\chi}[n])$ on $X_k(\Sigma)$. (For the
definition of $j_{!\ast}$, see [BBD] 1.4.22 and 1.4.24.)

\bigskip
Let $\delta$ be a rational convex polyhedral cone in $V^\ast$ of
dimension $n$, and let $\Sigma(\delta)$ be the fan defined in
Proposition 1.1. If $0$ is a face of $\delta$, then the map
$$M\cap \delta \to k, \; u\mapsto \left\{
\begin{array}{cc}
1 & \hbox { if } u=0,\\
0 & \hbox { if } u\not =0
\end{array}
\right.$$ is a semigroup homomorphism. It induces an epimorphism
of $k$-algebras $k[M\cap \delta]\to k$ and hence a closed
immersion
$$x_0: {\rm Spec}\, k\to {\rm Spec}\, k[M\cap \delta]=X_k(\Sigma(\delta)).$$
We call $x_0$ the {\it distinguished point} in
$X_k(\Sigma(\delta))$. It is fixed under the torus action on
$X_k(\Sigma(\delta))$.

\bigskip
\noindent {\bf Lemma 2.1.} Let $\delta$ be a rational convex
polyhedral cone in $V^\ast$ of dimension $n$ with $0$ being a face,
let $x_0:{\rm Spec}\, k\to X_k(\Sigma(\delta))$ be the distinguished
point in $X_k(\Sigma(\delta))$, and let $j:{\bf T}_k^n\to
X_k(\Sigma(\delta))$ be the open immersion of the open dense torus
in $X_k(\Sigma(\delta))$. If $\chi:\mu_m(k)^n\to \overline {\bf
Q}_l^\ast$ is a nontrivial character, then $x_0^\ast
(j_{!\ast}({\cal K}_{\chi}[n]))$ is acyclic and $x_0^\ast
(j_\ast{\cal K}_\chi)=0$. Moreover, we have $j_\ast \overline {\bf
Q}_l=\overline {\bf Q}_l$.

\bigskip
\noindent {\bf Proof.} Consider the morphism $$\bar\epsilon: {\bf
T}_k^n\times_k X_k(\Sigma(\delta))\to {\bf T}_k^n\times_k
X_k(\Sigma(\delta)),\; \bar\epsilon(t, x)=(t,tx),$$ where the
second component of $\bar \epsilon$ is defined by the action of
${\bf T}_k^n$ on $X_k^n(\Sigma(\delta))$. Note that $\bar\epsilon$
is an isomorphism. Denote the restriction of $\bar\epsilon$ to
${\bf T}_k^n\times_k{\bf T}_k^n$ by $\epsilon$. Consider the
commutative diagram
$$\begin{array}{ccccc}
{\bf T}_k^n\times_k{\bf T}_k^n&\stackrel{\epsilon}\to& {\bf
T}_k^n\times_k{\bf T}_k^n&\stackrel {{\rm pr}_2}\to &{\bf
T}_k^n \\
\downarrow {\rm id}\times j&&\downarrow {\rm id}\times
j&&\downarrow j \\
{\bf T}_k^n\times_kX_k(\Sigma(\delta))&\stackrel{\bar
\epsilon}\to& {\bf T}_k^n\times_kX_k(\Sigma(\delta))&\stackrel
{\bar {\rm pr}_2}\to &X_k(\Sigma(\delta)),
\end{array}$$
where ${\rm pr}_2$ and $\bar {\rm pr}_2$ are the projections to
the second components. Note that we have $${\rm pr}_2\epsilon =
s,$$ where $s$ is the multiplication on the torus. By the smooth
base change theorem and the fact that $\bar{\rm pr}_2^\ast[n]$ is
an exact functor with respect to the perverse t-structure ([BBD]
page 108-109), we have
$$\bar {\rm pr}_2^\ast (j_{!\ast}({\cal K}_\chi[n]))[n]\cong({\rm id}\times
j)_{!\ast}({\rm pr}_2^\ast {\cal K}_\chi[2n]).$$ So we have
\begin{eqnarray*}
\bar\epsilon^\ast \bar {\rm pr}_2^\ast (j_{!\ast}({\cal
K}_\chi[n]))[n] &\cong& \bar\epsilon^\ast ({\rm id}\times
j)_{!\ast}({\rm pr}_2^\ast{\cal K}_\chi[2n])\\
&\cong& ({\rm id}\times j)_{!\ast}(\epsilon^\ast {\rm pr}_2^\ast
{\cal
K}_\chi[2n])\\
&\cong&({\rm id}\times j)_{!\ast}(s^\ast {\cal K}_\chi[2n])\\
&\cong& ({\rm id}\times j)_{!\ast}({\rm pr}_1^\ast{\cal
K}_\chi\otimes {\rm pr}_2^\ast {\cal K}_\chi[2n]).
\end{eqnarray*}
Using [BBD] 1.4.24, 4.2.7 and 4.2.8, one can show
$$({\rm id}\times j)_{!\ast}({\rm pr}_1^\ast{\cal
K}_\chi\otimes {\rm pr}_2^\ast {\cal K}_\chi[2n])\cong {\rm
pr}_1^\ast ({\cal K}_\chi[n])\otimes \bar {\rm
pr}_2^\ast(j_{!\ast}({\cal K}_\chi[n])).$$ So we have
$$\bar\epsilon^\ast \bar {\rm pr}_2^\ast (j_{!\ast}({\cal
K}_\chi[n]))[n] \cong {\rm pr}_1^\ast ({\cal
K}_\chi[n])\otimes\bar {\rm pr}_2^\ast (j_{!\ast}({\cal
K}_\chi[n])).$$ Let $i$ be the morphism
$${\bf T}_k^n={\bf T}_k^n\times_k {\rm Spec}\, k\stackrel {{\rm id}\times x_0}
\to {\bf T}_k^n\times_k X_k(\Sigma(\delta))$$ and let $\pi:{\bf
T}_k^n\to{\rm Spec}\,k$ be the structure morphism. We have
\begin{eqnarray*}
i^\ast \bar\epsilon^\ast \bar {\rm pr}_2^\ast (j_{!\ast}({\cal
K}_\chi[n]))[n]&\cong&(\bar {\rm pr}_2\bar\epsilon i
)^\ast(j_{!\ast}({\cal K}_\chi[n]))[n]\\&\cong&(x_0\pi)^\ast
(j_{!\ast}({\cal K}_\chi[n]))[n],\\
i^\ast({\rm pr}_1^\ast ({\cal K}_\chi[n])\otimes\bar {\rm pr}_2^\ast
(j_{!\ast}({\cal K}_\chi[n])))&\cong&{\cal K}_\chi[n]\otimes
(x_0\pi)^\ast(j_{!\ast}({\cal K}_\chi[n])).
\end{eqnarray*}
We thus have an isomorphism
$$(x_0\pi)^\ast
(j_{!\ast}({\cal K}_\chi[n]))\cong {\cal K}_\chi\otimes
(x_0\pi)^\ast(j_{!\ast}({\cal K}_\chi[n])).$$ So for all $i$, we
have
$${\cal H}^i((x_0\pi)^\ast
(j_{!\ast}({\cal K}_\chi[n])))\cong {\cal K}_\chi\otimes {\cal H}^i(
(x_0\pi)^\ast(j_{!\ast}({\cal K}_\chi[n]))).$$ If $\chi$ is
nontrivial, this isomorphism implies that ${\cal H}^i((x_0\pi)^\ast
(j_{!\ast}({\cal K}_\chi[n])))=0$. Indeed, we can evaluate $H^0({\bf
T}_{\bar k}^n, \cdot)$ on both sides of the above isomorphism and
use the fact that $H^0({\bf T}_{\bar k}^n, {\cal K}_\chi)=0$. So we
have ${\cal H}^i(x_0^\ast (j_{!\ast}({\cal K}_\chi[n])))=0$ for all
$i$, and hence $x_0^\ast(j_{!\ast}({\cal K}_\chi[n]))$ is acyclic.

Similarly, one can show $x_0^\ast(j_{\ast}{\cal K}_\chi)=0$ if
$\chi$ is nontrivial. To prove $j_\ast\overline {\bf
Q}_l=\overline {\bf Q}_l$, we need to show for any nonempty
connected \'etale $X_k(\Sigma(\delta))$-scheme $U$, we have
$$(j_\ast\overline {\bf Q}_l)(U)=\overline {\bf Q}_l,$$ that is,
$$\overline {\bf
Q}_l(U\times_{X_k(\Sigma(\delta))}{\bf T}_k^n)=\overline {\bf
Q}_l.$$ It suffices to show $U\times_{X_k(\Sigma(\delta))}{\bf
T}_k^n$ is nonempty and connected. Indeed, by [F] 2.1,
$X_k(\Sigma(\delta))$ is normal. So $U$ is normal. As $U$ is
connected, $U$ is integral. $U\times_{X_k(\Sigma(\delta))}{\bf
T}_k^n$ is an open subscheme of $U$. It is nonempty since the
image of $U$ in $X_k(\Sigma(\delta))$ is nonempty open and hence
has nonempty intersection with the open dense torus ${\bf T}_k^n$.
Being a nonempty open subscheme of the integral scheme $U$,
$U\times_{X_k(\Sigma(\delta))}{\bf T}_k^n$ is also integral, and
in particular connected. This finishes the proof of the lemma.

\bigskip
\noindent {\bf Remark 2.2.} Keep the notations in Lemma 2.1. Denote
by $\bar x_0$ the geometric point associated to $x_0$. By [DL]
6.2.3, $x_0^\ast (j_{!\ast}(\overline {\bf Q}_l[n]))$ is pure of
weight $n$, that is, for any $i$ and any eigenvalue $\lambda$ of the
geometric Frobenius element  $F$ in ${\rm Gal}(\bar k/k)$ acting on
$H^i(\bar x_0^\ast (j_{!\ast}(\overline {\bf Q}_l[n])))$, $\lambda$
is an algebraic number, and all its galois conjugates have
archimedean absolute value $q^{\frac{i+n}{2}}$. Moreover, by [DL]
6.2.1, ${\rm dim}(H^i(\bar x_0^\ast (j_{!\ast}(\overline {\bf
Q}_l[n]))))$ coincides with the coefficient of $T^{i+n}$ of the
polynomial $\alpha(\delta)$ defined in the Introduction.

\bigskip
\noindent {\bf Lemma 2.3.} Let $\Sigma$ be a fan in $V$,
$\sigma\in \Sigma$ and $\delta$ the image of $\check \sigma$ under
the canonical homomorphism $V^\ast \to V^\ast/\sigma^\perp$. Note
that $\delta$ is a rational convex polyhedral cone in
$V^\ast/\sigma^\perp$ of dimension ${\rm
dim}(V^\ast/\sigma^\perp)$ and $0$ is a face of $\delta$. Let
$x_0$ be the distinguished point in $X_k(\Sigma(\delta))$, and let
$$j:{\bf T}_k^n\hookrightarrow X_k(\Sigma),\;j':{\bf T}_k^{{\rm
dim}(\sigma)}\hookrightarrow X_k(\Sigma(\delta))={\rm
Spec}\,k[(M/M\cap\sigma^\perp)\cap \delta]$$ be the immersions of
the open dense tori in toric schemes. Denote by
$$p_1:{\bf T}_k^n={\rm Spec}\, k[M]\to O_\sigma=
{\rm Spec}\, k [M\cap \sigma^\perp]$$ the morphism
induced by the canonical homomorphism
$$k[M\cap \sigma^\perp]\hookrightarrow k[M].$$ If ${\cal
K}_\chi=p_1^\ast {\cal K}_{\chi_1}$ for some Kummer sheaf ${\cal
K}_{\chi_1}$ on $O_\sigma$, then
\begin{eqnarray*}
(j_{!\ast}({\cal K}_\chi[n]))|_{O_\sigma}&\cong & ({\cal
K}_{\chi_1}[n-{\rm dim}(\sigma)])\otimes \pi^\ast x_0^\ast
({j'}_{!\ast}(\overline {\bf Q}_l[{\rm dim}(\sigma)])),\\
(j_\ast{\cal K}_\chi)|_{O_\sigma}&\cong& {\cal K}_{\chi_1},
\end{eqnarray*}
where $O_\sigma$ is considered as a subscheme of $X_k(\Sigma)$ as
in \S 1, and $\pi:O_\sigma\to {\rm Spec}\, k$ is the structure
morphism. If ${\cal K}_\chi$ is not the inverse image under $p_1$
of any Kummer sheaf on $O_\sigma$, then $(j_{!\ast}({\cal
K}_\chi[n]))|_{O_\sigma}$ is acyclic and $(j_\ast{\cal
K}_\chi)|_{O_\sigma}=0$.

\bigskip
\noindent {\bf Proof.} Since $M/M\cap\sigma^\perp$ is torsion free
and hence a free abelian group, the short exact sequence
$$0\to M\cap \sigma^\perp\to M\to M/M\cap \sigma^\perp\to 0$$
splits. Fix an isomorphism
$$M\cong M\cap \sigma^\perp\oplus M/M\cap \sigma^\perp.$$ This
isomorphism induces identifications
\begin{eqnarray*}
\check\sigma&\cong &\sigma^\perp\oplus \delta,\\
M\cap \check \sigma&\cong& M\cap \sigma^\perp\oplus (M/M\cap
\sigma^\perp)\cap \delta,\\
k[M\cap \check\sigma]&\cong& k[M\cap\sigma^\perp]\otimes_k
k[(M/M\cap \sigma^\perp)\cap\delta],\\
U_\sigma={\rm Spec}\, k[M\cap \check\sigma]&\cong& {\rm Spec}\,
k[M\cap\sigma^\perp]\times_k {\rm Spec}\,k[(M/M\cap
\sigma^\perp)\cap\delta]=O_\sigma\times_k X_k(\Sigma(\delta)),\\
{\bf T}_k^n={\rm Spec}\, k[M]&\cong& {\rm Spec}\, k[M\cap
\sigma^\perp]\times_k {\rm Spec}\,
k[M/M\cap\sigma^\perp]=O_\sigma\times_k {\bf T}_k^{{\rm
dim}(\sigma)}.
\end{eqnarray*}
We can find Kummer sheaves ${\cal K}_{\chi_1}$ on $O_\sigma$ and
${\cal K}_{\chi_2}$ on ${\bf T}_k^{{\rm dim}(\sigma)}$ so that
${\cal K}_\chi$ is identified with $$p_1^\ast{\cal
K}_{\chi_1}\otimes p_2^\ast{\cal K}_{\chi_2}$$ through the
isomorphism ${\bf T}_k^n\cong O_\sigma\times_k {\bf T}_k^{{\rm
dim}(\sigma)}$, where
$$p_1:O_\sigma\times_k {\bf T}_k^{{\rm dim}(\sigma)}
\to O_\sigma,\, p_2:O_\sigma\times_k {\bf T}_k^{{\rm
dim}(\sigma)}\to {\bf T}_k^{{\rm dim}(\sigma)}$$ are the
projections. Consider the commutative diagram
$$\begin{array}{ccccc}
O_\sigma&\to& U_\sigma& \stackrel{j}\leftarrow& {\bf T}_k^n\\
\parallel&&\downarrow \cong&&\downarrow \cong\\
O_\sigma &\stackrel{i}\to& O_\sigma\times_k
X_k(\Sigma(\delta))&\stackrel {{\rm id}\times j'}\leftarrow
&O_\sigma\times_k {\bf T}_k^{{\rm dim}(\sigma)},
\end{array}$$
where $i$ is the morphism
$$O_\sigma=O_\sigma\times_k{\rm Spec}\, k\stackrel {{\rm id}\times x_0}\to
O_\sigma\times_k X_k(\Sigma(\delta)).$$ Using [BBD] 1.4.24, 4.2.7
and 4.2.8, one can show $$({\rm id}\times
j')_{!\ast}\biggl(p_1^\ast{\cal K}_{\chi_1}[n-{\rm
dim}(\sigma)]\otimes p_2^\ast{\cal K}_{\chi_2}[{\rm
dim}(\sigma)]\biggr) \cong  {\rm pr}_1^\ast{\cal K}_{\chi_1}[n-{\rm
dim}(\sigma)]\otimes {\rm pr}_2^\ast{j'}_{!\ast}({\cal
K}_{\chi_2}[{\rm dim}(\sigma)]),$$ where $${\rm
pr}_1:O_\sigma\times_k X_k(\Sigma(\delta))\to O_\sigma,\; {\rm
pr}_2: O_\sigma\times_k X_k(\Sigma(\delta))\to X_k(\Sigma(\delta))$$
are the projections. So we have
\begin{eqnarray*}
(j_{!\ast}({\cal K}_\chi[n]))|_{O_\sigma}&\cong& i ^\ast({\rm
id}\times j')_{!\ast}\biggl(p_1^\ast{\cal K}_{\chi_1}[n-{\rm
dim}(\sigma)]\otimes p_2^\ast{\cal
K}_{\chi_2}[{\rm dim}(\sigma)]\biggr)\\
&\cong& i^\ast \biggl({\rm pr}_1^\ast{\cal K}_{\chi_1}[n-{\rm
dim}(\sigma)]\otimes {\rm pr}_2^\ast{j'}_{!\ast}({\cal
K}_{\chi_2}[{\rm
dim}(\sigma)])\biggr)\\
&\cong&{\cal K}_{\chi_1}[n-{\rm dim}(\sigma)]\otimes\pi^\ast
x_0^\ast ({j'}_{!\ast}({\cal K}_{\chi_2}[{\rm dim}(\sigma)])),
\end{eqnarray*}
that is, $$(j_{!\ast}({\cal K}_\chi[n]))|_{O_\sigma}\cong{\cal
K}_{\chi_1}[n-{\rm dim}(\sigma)]\otimes\pi^\ast x_0^\ast
({j'}_{!\ast}({\cal K}_{\chi_2}[{\rm dim}(\sigma)])).$$ Similarly,
we have
$$(j_\ast{\cal K}_\chi )|_{O_\sigma}\cong
{\cal K}_{\chi_1}\otimes\pi^\ast x_0^\ast({j'}_\ast{\cal
K}_{\chi_2}).$$ We then apply Lemma 2.1.

\bigskip
\noindent {\bf Proposition 2.4.} Let $\Sigma$ be a fan in $V$,
$\sigma\in \Sigma$,  $$j:{\bf T}_k^n\hookrightarrow
X_k(\Sigma),\;j_1:O_\sigma={\rm Spec}\, k[M\cap\sigma^\perp]\to
V(\sigma)=X_k({\rm star}(\sigma))$$  the immersions of the open
dense tori in toric schemes, and $$p_1:{\bf T}_k^n={\rm Spec}\,
k[M]\to O_\sigma={\rm Spec}\, k [M\cap \sigma^\perp]$$ the morphism
induced by the canonical homomorphism
$$k[M\cap \sigma^\perp]\hookrightarrow k[M].$$ If ${\cal
K}_\chi=p_1^\ast {\cal K}_{\chi_1}$ for a Kummer sheaf ${\cal
K}_{\chi_1}$ on $O_\sigma$, then $$(j_\ast{\cal
K}_\chi)|_{V(\sigma)}\cong {j_1}_\ast{\cal K}_{\chi_1}.$$ If ${\cal
K}_\chi$ is not the inverse image under $p_1$ of any Kummer sheaf on
$O_\sigma$, then $$(j_{\ast}{\cal K}_\chi)|_{V(\sigma)}=0.$$

\bigskip
\noindent {\bf Proof.}  First suppose ${\cal K}_\chi$ is not the
inverse image under $p_1$ of any Kummer sheaf on $O_\sigma$. Then
for any $\tau\in \Sigma$ with $\sigma\prec\tau$, ${\cal K}_\chi$ is
not the inverse image of any Kummer sheaf on $O_\tau$. By Lemma 2.3,
we then have $(j_\ast{\cal K}_\chi)|_{O_\tau}=0$. As
$V(\sigma)=\bigcup_{\sigma\prec \tau} O_\tau$, we have $(j_\ast{\cal
K}_\chi)|_{V(\sigma)}=0.$

Now suppose ${\cal K}_\chi=p_1^\ast{\cal K}_{\chi_1}$ for a Kummer
sheaf ${\cal K}_1$ on $O_\sigma$. As in the proof of Lemma 2.3, we
have a commutative diagram
$$\begin{array}{ccc}
{\bf T}_k^n&\stackrel {j_\sigma} \to & U_\sigma \\
\downarrow\cong &&\downarrow\cong \\
O_\sigma\times_k {\bf T}_k^{{\rm dim}(\sigma)}&\stackrel {{\rm
id}\times j'}\to& O_\sigma\times_k X_k(\Sigma(\delta)),
\end{array}$$
where $\delta$ is the image of $\check\sigma$ in
$V^\ast/\sigma^\perp$, and $j_\sigma$ and $j'$ are the open
immersions of the tori in $U_\sigma$ and in $X_k(\Sigma(\delta))$,
respectively. As ${\cal K}_\chi=p_1^\ast {\cal K}_{\chi_1}$, the
sheaf ${j_\sigma}_\ast{\cal K}_\chi$ can be identified with $({\rm
id}\times j')_\ast (p^\ast {\cal K}_{\chi_1})$, where
$p:O_\sigma\times_k {\bf T}_k^{{\rm dim}(\sigma)}\to O_\sigma$ is
the projection. By [BBD] 4.2.7 (which can be deduced from [SGA
$4\frac{1}{2}$] [Th. finitude] 2.16 and Appendix 2.10), we have
$$({\rm id}\times j')_\ast (p^\ast {\cal
K}_{\chi_1})\cong {\rm pr}_1^\ast {\cal K}_{\chi_1}\otimes {\rm
pr}_2^\ast({j'}_\ast \overline {\bf Q}_l),$$ where $${\rm
pr}_1:O_\sigma\times_k X_k(\Sigma(\delta))\to O_\sigma,\; {\rm
pr}_2:O_\sigma\times_k X_k(\Sigma(\delta))\to X_k(\Sigma(\delta))$$
are the projections. As ${j'}_\ast \overline {\bf Q}_l=\overline
{\bf Q}_l$ by Lemma 2.1, we have
$$({\rm id}\times j')_\ast (p^\ast {\cal
K}_{\chi_1})\cong {\rm pr}_1^\ast {\cal K}_{\chi_1}.$$ It follows
that the canonical homomorphism
$$q_\sigma^\ast {\cal K}_{\chi_1}\to {j_\sigma}_\ast j_\sigma^\ast
q_\sigma^\ast {\cal K}_{\chi_1}$$ is an isomorphism, where
$q_\sigma:U_\sigma\to O_\sigma$ is the morphism defined by the
canonical homomorphism $k[M\cap \sigma^\perp]\hookrightarrow
k[M\cap\check \sigma]$.

For each $\tau\in \Sigma$ with $\sigma\prec\tau$, fix notations by
the following commutative diagram:
$$\begin{array}{ccrcl}
&&O_\sigma&\to& O_\tau\\
&\stackrel{p_1} \nearrow&q_\sigma \uparrow &&\uparrow q_\tau\\
{\bf T}_k^n&\stackrel {j_\sigma}\hookrightarrow
&U_\sigma&\stackrel {j_{\sigma\tau}}\hookrightarrow &U_\tau\\
&&i_\sigma \uparrow&&\uparrow i_\tau \\
&&O_\sigma&\stackrel {l_\tau}\hookrightarrow &U_\tau(\sigma)={\rm
Spec}\, k[M\cap \sigma^\perp\cap \check\tau],
\end{array}$$
and let $j_\tau:{\bf T}_k^n\hookrightarrow U_\tau$ be the canonical
open immersion. If ${\cal K}_\chi$ is the inverse image of a Kummer
sheaf ${\cal K}_\tau$ on $O_\tau$, then the same argument as above
shows that the following canonical homomorphisms are isomorphisms:
$$q_\tau^\ast{\cal K}_{\chi_\tau}\stackrel{\cong}\to
{j_\tau}_\ast j_\tau^\ast q_\tau^\ast {\cal K}_{\chi_\tau}\eqno
(1)$$ $$(q_\tau i_\tau)^\ast {\cal K}_{\chi_\tau}\stackrel
{\cong}\to {l_\tau}_\ast l_\tau^\ast (q_\tau i_\tau)^\ast {\cal
K}_{\chi_\tau}\eqno (2)$$

The isomorphism $$q_\sigma^\ast {\cal K}_{\chi_1}\stackrel \cong
\to {j_\sigma}_\ast j_\sigma^\ast q_\sigma^\ast {\cal
K}_{\chi_1}$$ induces an isomorphism
$${i_\sigma}^\ast q_\sigma^\ast {\cal K}_{\chi_1}\stackrel \cong
\to {i_\sigma}^\ast{j_\sigma}_\ast j_\sigma^\ast q_\sigma^\ast
{\cal K}_{\chi_1}.$$ So we have an isomorphism
$${\cal K}_{\chi_1}\stackrel\cong\to j_1^\ast((j_\ast{\cal
K}_\chi)|_{V(\sigma)})\eqno (3) $$ Applying the adjointness of
$(j_1^\ast, {j_1}_\ast)$ to the inverse of this isomorphism, we get
a homomorphism
$$(j_\ast{\cal
K}_\chi)|_{V(\sigma)}\to {j_1}_\ast {\cal K}_{\chi_1}\eqno (4)$$
Let's show this is an isomorphism. This would prove our
proposition.

We have $V(\sigma)=\bigcup_{\sigma\prec\tau}O_\tau$. It suffices to
prove the restriction of the homomorphism (4) to each $O_\tau$
$(\sigma\prec\tau)$ is an isomorphism. If ${\cal K}_\chi$ is not the
inverse image of any Kummer sheaf on $O_\tau$, the both
$(j_\ast{\cal K}_\chi)|_{O_\tau}$ and $({j_1}_\ast {\cal
K}_{\chi_1})|_{O_\tau}$ vanish by Lemma 2.3, and the restriction to
$O_\tau$ of (4) is trivially an isomorphism. Now suppose ${\cal
K}_\chi$ is the inverse image of a Kummer sheaf ${\cal K}_\tau$ on
$O_\tau$. Let's prove the restriction of (4) to the set
$U_\tau(\sigma)$ (which contains $O_\tau$) is an isomorphism.

Recall that to construct the homomorphism (4), we first construct
an isomorphism (3), and then apply the adjointness of
$(j_1^\ast,{j_1}_\ast)$ to the inverse of (3). The isomorphism (3)
can also be described as follows: From the isomorphism (1), we get
an isomorphism
$$l_\tau^\ast i_\tau^\ast q_\tau^\ast {\cal
K}_{\chi_\tau}\stackrel{\cong}\to l_\tau^\ast i_\tau^\ast
{j_\tau}_\ast j_\tau^\ast q_\tau^\ast {\cal K}_{\chi_\tau}\eqno(5)$$
Note that (3) can be identified with (5). This can be seen by
applying $i_\sigma^\ast$ to the commutative diagram
$$\begin{array}{ccc}
j_{\sigma\tau}^\ast q_\tau^\ast {\cal K}_{\chi_\tau}&\stackrel
\cong\to& j_{\sigma\tau}^\ast {j_\tau}_\ast j_\tau^\ast q_\tau^\ast
{\cal
K}_{\chi_\tau} \\
\cong\downarrow&&\downarrow \cong\\
q_\sigma^\ast{\cal K}_{\chi_1}&\stackrel\cong\to & {j_\sigma}_\ast
j_\sigma^\ast q_\sigma^\ast {\cal K}_{\chi_1}.
\end{array}$$
Applying the adjointness of $(l_\tau^\ast, {l_\tau}_\ast)$ to the
inverse of (5), we get
$$i_\tau^\ast {j_\tau}_\ast j_\tau^\ast q_\tau^\ast {\cal
K}_{\chi_\tau}\to {l_\tau}_\ast l_\tau^\ast i_\tau^\ast q_\tau^\ast
{\cal K}_{\chi_\tau}\eqno (6)$$ Note that (6) can be identified with
the restriction of (4) to $U_\tau(\sigma)$. So it suffices to show
(6) is an isomorphism. But (6) coincides with the composition
$$i_\tau^\ast {j_\tau}_\ast j_\tau^\ast q_\tau^\ast {\cal
K}_{\chi_\tau}\to i_\tau^\ast q_\tau^\ast {\cal
K}_{\chi_\tau}\to{l_\tau}_\ast l_\tau^\ast i_\tau^\ast q_\tau^\ast
{\cal K}_{\chi_\tau},$$ where the first arrow is obtained by
applying $i_\tau^\ast$ to the inverse of (1), and the second arrow
can be identified with (2). As (1) and (2) are isomorphisms, (6) is
also an isomorphism. This finishes the proof of the proposition.

\bigskip
In the following, $A$ is a finitely generated $k$-algebra. For any
fan $\Sigma$ in $V$ and any cone $\sigma$ in $\Sigma$, we denote
${\rm Spec}\, A[M\cap \check\sigma]$ (resp. ${\rm Spec}\, k[M\cap
\check\sigma]$) by $U_\sigma$ (resp.${U_\sigma}_k$), and ${\rm
Spec}\, A[M\cap \sigma^\perp]$ (resp. ${\rm Spec}\, k[M\cap
\sigma^\perp])$ by $O_\sigma$ (resp. ${O_\sigma}_k$). They are
considered as subschemes of $X_A(\Sigma)$ (resp. $X_k(\Sigma)$.)

\bigskip
\noindent {\bf Proposition 2.5.} Let $A$ be a finitely generated
$k$-algebra, let $\Sigma$ be a fan in $V$ and let $Y$ be an
effective Cartier divisor of $X_A(\Sigma)$. Denote by $j:{\bf
T}_k^n\to X_k(\Sigma)$ the immersion of the open dense torus in
$X_k(\Sigma)$, and denote by $(j_{!\ast}({\cal K}_\chi[n]))|_Y$ and
$(j_\ast{\cal K}_\chi)|_Y$ the inverse images under the composition
$$Y\to X_A(\Sigma)\to X_k(\Sigma)$$ of $j_{!\ast}({\cal
K}_\chi[n])$ and $j_\ast{\cal K}_\chi$, respectively. Let $y\in Y$
and let $\sigma$ be the unique cone in $\Sigma$ such that $y\in
O_\sigma$. Suppose the scheme theoretic intersection $Y\cap
O_\sigma$ is smooth over $A$ at $y$, and is not equal to
$O_\sigma$ in any neighborhood of $y$. Then $Y$ is universally
locally acyclic over $A$ at $y$ relative to $(j_{!\ast}({\cal
K}_\chi))|_Y$ and relative to $(j_\ast{\cal K}_\chi)|_Y$.

\bigskip
\noindent {\bf Proof.} As in the proof of Lemma 2.3, $M\cap
\sigma^\perp$ is a direct factor of $M$. Choose a sublattice $M'$ of
$M$ such that
$$M\cong M\cap \sigma^\perp\oplus M'.$$ Let
$\delta$ be the image of $\check\sigma$ in $M'\otimes_{\bf Z}{\bf
R}$. It is a cone in $M'\otimes_{\bf Z}{\bf R}$ of dimension ${\rm
dim}(\sigma)$ with $0$ being a face. We have
\begin{eqnarray*}
\check\sigma&=&\sigma^\perp\oplus \delta,\\
M\cap \check \sigma&=& M\cap \sigma^\perp\oplus M'\cap \delta,\\
A[M\cap \check\sigma]&\cong& A[M\cap\sigma^\perp]\otimes_A
A[M'\cap\delta],\\
U_\sigma={\rm Spec}\, A[M\cap \check\sigma]&\cong& {\rm Spec}\,
A[M\cap\sigma^\perp]\times_A {\rm Spec}\,A[M'
\cap\delta]=O_\sigma\times_A {\rm Spec}\,A[M' \cap\delta].
\end{eqnarray*}
We have a commutative diagram of Cartesian squares
$$\begin{array}{ccccc}
&&Y\cap U_\sigma&\leftarrow & Y\cap O_\sigma\\
&&\downarrow&&\downarrow\\
{\cal O}_\sigma&\leftarrow& U_\sigma&\leftarrow & O_\sigma \\
\downarrow&&\downarrow&&\downarrow \\
{\rm Spec} \,A&\leftarrow&{\rm Spec}\,A[M' \cap\delta]&\leftarrow&
{\rm Spec}\,A,
\end{array}$$
where $U_\sigma\to O_\sigma$ is the projection, $O_\sigma\to
U_\sigma$ is the canonical closed immersion, and  ${\rm Spec}\,
A\to {\rm Spec}\,A[M' \cap\delta]$ is defined by the homomorphism
of $A$-algebras
$$A[M'\cap\delta]\to A,\;
\chi^u\mapsto \left\{
\begin{array}{cc}
1 & \hbox { if } u=0, \\
0 & \hbox { if } u\not =0.
\end{array}\right.$$
Choose an open neighborhood $W$ of $y$ in $U_\sigma$ and $g\in
{\cal O}_{X_A(\Sigma)}(W)$ so that $Y\cap W$ is defined by $g$.
Since $Y\cap O_\sigma$ is smooth over $A$ at $y$ and is not equal
to $O_\sigma$ near $y$, by [SGA 1], Expos\'e II, Th\'eor\`eme 4.10
(iii), the image $dg$ in $(\Omega_{O_\sigma/A}^1)_y\otimes _{{\cal
O}_{O_\sigma,y}} k(y)$ is nonzero. Hence the image of $dg$ in
$(\Omega_{U_\sigma/A[M'\cap\delta]}^1)_y\otimes_{{\cal
O}_{U_\sigma,y}}k(y)$ is nonzero. By [SGA 1], Expos\'e II,
Th\'eor\`eme 4.10 (ii), there exists an open neighborhood $W'$ of
$y$ in $U_\sigma$ and an \'etale $A[M'\cap\delta]$-morphism
$$W'\to {\rm Spec} A[M'\cap \delta][T_1,\ldots, T_r]$$ such that
$Y\cap W'$ is the inverse image of the closed subscheme of ${\rm
Spec} A[M'\cap \delta][T_1,\ldots, T_r]$ defined by the equation
$T_r=0$.
$$\begin{array}{rcl}
Y\cap W'&\to& {\rm Spec}\, A[M'\cap \delta][T_1,\ldots, T_{r-1}]\\
\downarrow&&\downarrow\\
W'&\to&{\rm Spec} A[M'\cap \delta][T_1,\ldots,
T_r] \\
&\searrow&\downarrow\\
&&{\rm Spec}\, A[M'\cap \delta].
\end{array}$$
We have
$${\bf T}_k^n={\rm Spec}\, k[M]\cong {O_\sigma}_k
\times_k {\rm Spec}\, k[M'].$$ Let
$$p_1: {\bf T}_k^n\to {O_\sigma}_k,\;
p_2:{\bf T}_k^n\to{\rm Spec}\, k[M']$$ be the projections. We can
find Kummer sheaves ${\cal K}_{\chi_1}$ on ${O_\sigma}_k$ and
${\cal K}_{\chi_2}$ on ${\rm Spec}\, k[M']$ so that ${\cal
K}_\chi$ is identified with $p_1^\ast{\cal K}_{\chi_1}\otimes
p_2^\ast{\cal K}_{\chi_2}.$ Using [BBD] 1.4.24, 4.2.7 and 4.2.8,
one can verify that on ${U_\sigma}_k$, we have
\begin{eqnarray*}
j_{!\ast}({\cal K}_\chi[n])&\cong& {\rm pr}_1^\ast({\cal
K}_{\chi_1}[n-{\rm dim}(\sigma)])\otimes {\rm
pr}_2^\ast({j'}_{!\ast}({\cal K}_{\chi_2}[{\rm dim}(\sigma)])),\\
j_\ast{\cal K}_\chi&\cong&{\rm pr}_1^\ast{\cal K}_{\chi_1}\otimes
{\rm pr}_2^\ast({j'}_\ast{\cal K}_{\chi_2}),
\end{eqnarray*}
where
\begin{eqnarray*}
{\rm pr}_1: {U_\sigma}_k\cong {O_\sigma}_k\times_k {\rm Spec}\,
k[M'\cap\delta] &\to& {O_\sigma}_k,\\
{\rm pr}_2: {U_\sigma}_k\cong {O_\sigma}_k \times_k {\rm Spec}\,
k[M'\cap\delta]&\to& {\rm Spec}\, k[M'\cap \delta]
\end{eqnarray*}
are the projections, and $j':{\rm Spec}\, k[M']\to {\rm Spec}\,
k[M'\cap \delta]$ is the immersion of the open dense torus in ${\rm
Spec}\, k[M'\cap \delta]=X_k(\Sigma(\delta))$.
$$\begin{array}{ccc}
{\bf T}_k^n&\stackrel {j}\to& {U_\sigma}_k\\
\downarrow \cong&&\downarrow \cong\\
{O_\sigma}_k \times_k {\rm Spec}\,k[M']&\stackrel {{\rm id}\times
j'} \to&{O_\sigma}_k\times_k {\rm Spec}\, k[M'\cap\delta].
\end{array}$$
By [SGA $4\frac{1}{2}$] [Th. finitude] Th\'eor\`eme 2.13 (with
$S={\rm Spec}\,k$), the morphism
$${\rm Spec}\, A[M'\cap \delta][T_1,\ldots, T_{r-1}]\to {\rm Spec}\,A$$ is universally
locally acyclic relative to the inverse images of
${j'}_{!\ast}({\cal K}_{\chi_2}[{\rm dim}(\sigma)])$ and of
${j'}_{\ast}{\cal K}_{\chi_2}$ under the composition
$${\rm Spec}\, A[M'\cap \delta]
[T_1,\ldots, T_{r-1}]\to {\rm Spec}\, k[M'\cap \delta][T_1,\ldots,
T_{r-1}]\to {\rm Spec}\, k[M'\cap \delta].$$ As the morphism
$$Y\cap W'\to {\rm Spec}\, A[M'\cap \delta]
[T_1,\ldots, T_{r-1}]$$ is \'etale, we see $Y\cap W'$ is
universally locally acyclic over $A$ relative to the inverse
images of ${j'}_{!\ast}({\cal K}_{\chi_2}[{\rm dim}(\sigma)])$ and
of ${j'}_{\ast}{\cal K}_{\chi_2}$ under the composition
$$Y\cap W'\to {\rm Spec}\,A[M'\cap \delta]\to {\rm Spec}\,k[M'\cap
\delta].$$ Note that these inverse images coincides with the inverse
images of ${\rm pr}_2^\ast({j'}_{!\ast}({\cal K}_{\chi_2}[{\rm
dim}(\sigma)]))$ and of ${\rm pr}_2^\ast({j'}_\ast{\cal
K}_{\chi_2})$ under the composition
$$Y\cap W'\to W'\to U_\sigma\to {U_\sigma}_k.$$ By Lemma 2.6 below, $Y\cap W'$
is universally locally acyclic over $A$ relative to the inverse
images of  ${\rm pr}_1^\ast({\cal K}_{\chi_1}[n-{\rm
dim}(\sigma)])\otimes {\rm pr}_2^\ast({j'}_{!\ast}({\cal
K}_{\chi_2}[{\rm dim}(\sigma)]))$ and of ${\rm pr}_1^\ast{\cal
K}_{\chi_1}\otimes {\rm pr}_2^\ast({j'}_\ast{\cal K}_{\chi_2})$.
Therefore  $Y\cap W'$ is universally locally acyclic over $A$
relative to the inverse images of $j_{!\ast}({\cal K}_\chi[n])$
and of $j_{\ast}{\cal K}_\chi$. This proves our assertion.

\bigskip
\noindent {\bf Lemma 2.6.} Let $f:X\to S$ be a morphism, ${\cal K}$
an object in the derived category $D^+(X, {\bf Z}/l^\alpha)$ of
\'etale sheaves of $({\bf Z}/l^\alpha)$-modules, and ${\cal F}$ a
flat locally constant sheaf of $({\bf Z}/l^\alpha)$-modules. Suppose
$f$ is locally acyclic relative to ${\cal K}$. Then $f$ is locally
acyclic relative to ${\cal F}\otimes {\cal K}$.

\bigskip
\noindent {\bf Proof.} Let $\bar s$ be an arbitrary geometric point
in $S$, let $\tilde S_{\bar s}$ be the strict henselization of $S$
at $\bar s$. and let $\bar t$ be an arbitrary  geometric point in
$\tilde S_{\bar s}$. Fix notations by the following commutative
diagram of Cartesian squares
$$\begin{array}{ccccc}
X_{\bar t}&\stackrel {j} \to& X\times_S
\tilde S_{\bar s}&\stackrel {i}\leftarrow & X_{\bar s}\\
\downarrow&&\downarrow &&\downarrow \\
\bar t&\to&\tilde S_{\bar s}&\leftarrow &\bar s.
\end{array}$$
Since $f$ is locally acyclic relative to ${\cal K}$, the canonical
morphism
$$i^\ast ({\cal K}|_{X\times_S \tilde S_{\bar s}})\to i^\ast Rj_\ast j^\ast
({\cal K}|_{X\times_S \tilde S_{\bar s}})$$ is an isomorphism.
Since ${\cal F}$ is flat, we have an isomorphism
$$i^\ast({\cal F}|_{X\times_S \tilde S_{\bar s}})\otimes i^\ast
({\cal K}|_{X\times_S \tilde S_{\bar s}}) \stackrel{\cong} \to
i^\ast ({\cal F}|_{X\times_S \tilde S_{\bar s}})\otimes i^\ast
Rj_\ast j^\ast ({\cal K}|_{X\times_S \tilde S_{\bar s}}),$$ that
is
$$i^\ast(({\cal F}\otimes {\cal K})|_{X\times_S \tilde S_{\bar s}})
{\cong} i^\ast ({\cal F}|_{X\times_S \tilde S_{\bar s}}\otimes
Rj_\ast j^\ast ({\cal K}|_{X\times_S \tilde S_{\bar s}})).$$ Since
${\cal F}$ is flat and locally constant, by [SGA 4] Expos\'e XVII
5.2.11.1, we have
\begin{eqnarray*}
{\cal F}|_{X\times_S \tilde S_{\bar s}}\otimes  Rj_\ast j^\ast
({\cal K}|_{X\times_S \tilde S_{\bar s}})&\cong& Rj_\ast (j^\ast
({\cal
F}|_{X\times_S \tilde S_{\bar s}})\otimes j^\ast
({\cal K}|_{X\times_S \tilde S_{\bar s}}))\\
&\cong& Rj_\ast j^\ast (({\cal F}\otimes {\cal K})|_{X\times_S
\tilde S_{\bar s}}).
\end{eqnarray*}
So the canonical morphism
$$i^\ast (({\cal F}\otimes {\cal K})|_{X\times_S \tilde S_{\bar s}})\to i^\ast Rj_\ast
j^\ast (({\cal F}\otimes {\cal K})|_{X\times_S \tilde S_{\bar
s}})$$ is an isomorphism. Therefore ${\cal F}\otimes {\cal K}$ is
locally acyclic.

\bigskip
\noindent {\bf Proposition 2.7.} Let $A$ be a finitely generated
$k$-algebra, let $\Sigma$ be a fan in $V$, let $Y$ be an effective
Cartier divisor of $X_A(\Sigma)$, and let $j_Y:Y\cap {\bf
T}_A^n\to Y$ be the open immersion induced by the open immersion
of the open dense torus in $X_A(\Sigma)$. Suppose $A$ is an
integral domain of dimension $d$ and is smooth over $k$, and
suppose for any $\sigma\in\Sigma$, the scheme theoretic
intersection $Y\cap O_\sigma$ is smooth over $k$ and is not equal
to $O_\sigma$ near any point. (In particular, taking $\sigma=0$,
we see $Y\cap {\bf T}_A^n$ is smooth over $k$ and is not equal to
${\bf T}_A^n$ anywhere). Then we have
\begin{eqnarray*}
{j_Y}_{!\ast}({\cal K}_\chi|_{Y\cap {\bf T}_A^n}[n+d-1])&\cong&
(j_{!\ast}({\cal K}_\chi[n]))|_Y[d-1],\\
{j_Y}_{\ast}({\cal K}_\chi|_{Y\cap {\bf T}_A^n})&\cong&
(j_{\ast}{\cal K}_\chi)|_Y,
\end{eqnarray*}
where $(j_{!\ast}({\cal K}_\chi[n]))|_Y$ and $(j_\ast{\cal
K}_\chi)|_Y$ are the same as in Proposition 2.5, and ${\cal
K}_\chi|_{Y\cap {\bf T}_A^n}$ is the inverse image of ${\cal
K}_\chi$ under the composition $$Y\cap {\bf T}_A^n\to {\bf T}_A^n\to
{\bf T}_k^n.$$

\bigskip
\noindent {\bf Proof.} Let $y$ be an arbitrary point in $Y$. Choose
$\sigma\in\Sigma$ such that $y\in O_\sigma$. Choose a sublattice
$M'$ of $M$ such that
$$M\cong M\cap \sigma^\perp\oplus M',$$ and let
$\delta$ be the image of $\check\sigma$ in $M'\otimes_{\bf Z}{\bf
R}$. Again we have
$$
U_\sigma={\rm Spec}\, A[M\cap \check\sigma]\cong {\rm Spec}\,
A[M\cap\sigma^\perp]\times_A {\rm Spec}\,A[M'
\cap\delta]=O_\sigma\times_A {\rm Spec}\,A[M' \cap\delta],$$ and we
have a commutative diagram of Cartesian squares
$$\begin{array}{ccccc}
&&Y\cap U_\sigma&\leftarrow & Y\cap O_\sigma\\
&&\downarrow&&\downarrow\\
{\cal O}_\sigma&\leftarrow& U_\sigma&\leftarrow & O_\sigma \\
\downarrow&&\downarrow&&\downarrow \\
{\rm Spec} \,A&\leftarrow&{\rm Spec}\,A[M' \cap\delta]&\leftarrow&
{\rm Spec}\,A\\
\downarrow&&\downarrow&&\downarrow \\
{\rm Spec} \,k&\leftarrow&{\rm Spec}\,k[M' \cap\delta]&\leftarrow&
{\rm Spec}\,k.\\
\end{array}$$
Choose an open neighborhood $W$ of $y$ in $U_\sigma$ and $g\in
{\cal O}_{X_A(\Sigma)}(W)$ so that $Y\cap W$ is defined by $g$.
Since $A$ is a smooth $k$-algebra, $O_\sigma={\rm Spec}\, A[M\cap
\sigma^\perp]$ is smooth over $k$. Since $Y\cap O_\sigma$ is
smooth over $k$ at $y$ and is not equal to $O_\sigma$ near $y$, by
[SGA 1] Expos\'e II, Th\'eor\`eme 4.10 (iii), the image of $dg$ in
$(\Omega_{O_\sigma/k}^1)_y\otimes _{{\cal O}_{O_\sigma,y}} k(y)$
is nonzero. Hence the image of $dg$ in
$(\Omega_{U_\sigma/k[M'\cap\delta]}^1)_y\otimes_{{\cal
O}_{U_\sigma,y}}k(y)$ is nonzero. By [SGA 1] Expos\'e II,
Th\'eor\`eme 4.10 (i), there exists an open neighborhood $W'$ of
$y$ in $U_\sigma$ so that $Y\cap W'$ is smooth over
$k[M'\cap\delta]$. Choose Kummer sheaves ${\cal K}_{\chi_1}$ on
${O_\sigma}_k$ and ${\cal K}_{\chi_2}$ on ${\rm Spec}\, k[M']$ so
that ${\cal K}_\chi$ is identified with $p_1^\ast{\cal
K}_{\chi_1}\otimes p_2^\ast{\cal K}_{\chi_2},$ where $$p_1: {\bf
T}_k^n\to {O_\sigma}_k,\; p_2:{\bf T}_k^n\to{\rm Spec}\, k[M']$$
are the projections. As in the proof of Proposition 2.5, on
${U_\sigma}_k$, we have
\begin{eqnarray*}
j_{!\ast}({\cal K}_\chi[n])&\cong& {\rm pr}_1^\ast({\cal
K}_{\chi_1}[n-{\rm dim}(\sigma)])\otimes {\rm
pr}_2^\ast({j'}_{!\ast}({\cal K}_{\chi_2}[{\rm dim}(\sigma)])),\\
j_\ast{\cal K}_\chi&=&{\rm pr}_1^\ast{\cal K}_{\chi_1}\otimes {\rm
pr}_2^\ast({j'}_\ast{\cal K}_{\chi_2}),
\end{eqnarray*}
where
\begin{eqnarray*}
{\rm pr}_1: {U_\sigma}_k\cong {O_\sigma}_k\times_k {\rm Spec}\,
k[M'\cap\delta] &\to& {O_\sigma}_k,\\
{\rm pr}_2: {U_\sigma}_k\cong {O_\sigma}_k \times_k {\rm Spec}\,
k[M'\cap\delta]&\to& {\rm Spec}\, k[M'\cap \delta]
\end{eqnarray*}
are the projections, and $j':{\rm Spec}\, k[M']\to {\rm Spec}\,
k[M'\cap \delta]$ is the immersion of the open dense torus in ${\rm
Spec}\, k[M'\cap \delta]=X_k(\Sigma(\delta))$. We have the following
commutative diagram of Cartesian squares
$$\begin{array}{ccccc}
Y\cap W'\cap {\bf T}_A^n&\stackrel {j_Y}\hookrightarrow &Y\cap W'&& \\
\downarrow&&\downarrow &&\\
{\bf T}_A^n={\rm Spec}\, A[M]&\hookrightarrow & U_\sigma={\rm
Spec}\, A[M\cap \check\sigma]&&\\
\downarrow&&\downarrow&&\\
{\bf T}_k^n={\rm Spec}\, k[M]&\hookrightarrow & {U_\sigma}_k={\rm
Spec}\, k[M\cap \check\sigma]&\stackrel{{\rm pr}_1}\to& {O_{\sigma}}_k\\
p_2\downarrow&&\downarrow {\rm pr}_2&&\downarrow\\
{\rm Spec}\, k[M']&\stackrel{j'}\hookrightarrow & {\rm Spec}\,
k[M'\cap \delta]&\to& {\rm Spec} \, k.
\end{array}$$
Since $Y\cap W'$ is smooth over ${\rm Spec}\, k[M'\cap \delta]$,
by the smooth base change theorem, on $Y\cap W'$, we have $$({\rm
pr}_2^\ast({j'}_\ast{\cal K}_{\chi_2}))|_{Y\cap W'}\cong
{j_Y}_\ast ((p_2^\ast {{\cal K}_{\chi_2}})|_{Y\cap {\bf T}_A^n})$$
So we have
\begin{eqnarray*}
(j_\ast{\cal K}_\chi)|_{Y\cap W'}&\cong& ({\rm pr}_1^\ast{\cal
K}_{\chi_1}\otimes {\rm pr}_2^\ast({j'}_\ast{\cal
K}_{\chi_2}))|_{Y\cap W'}\\&\cong& ({\rm pr}_1^\ast{\cal
K}_{\chi_1})|_{Y\cap W'}\otimes {j_Y}_\ast ((p_2^\ast{{\cal
K}_{\chi_2}})|_{Y\cap {\bf T}_A^n})\\
&\cong& {j_Y}_\ast ({j_Y}^\ast (({\rm pr}_1^\ast{\cal
K}_{\chi_1})|_{Y\cap W'})\otimes (p_2^\ast{{\cal
K}_{\chi_2}})|_{Y\cap {\bf T}_A^n})\\
&\cong& {j_Y}_\ast((p_1^\ast {\cal K}_{\chi_1}\otimes p_2^\ast {\cal
K}_{\chi_2})|_{Y\cap {\bf T}_A^n})\\
&\cong& {j_Y}_\ast({\cal K}_\chi|_{Y\cap {\bf T}_A^n}),
\end{eqnarray*}
where for the third isomorphism, we use [SGA 4] Expos\'e XVII,
5.2.11.1. This proves that $${j_Y}_{\ast}({\cal K}_\chi|_{Y\cap
{\bf T}_A^n})\cong (j_{\ast}{\cal K}_\chi)|_Y$$ holds when
restricted to $Y\cap W'$. As $W'$ is a neighborhood of $y$, and
$y$ can be taken to be any point in $Y$, we have the above
isomorphism everywhere on $Y$. Denote the smooth morphism $Y\cap
W'\to {\rm Spec}\, k[M'\cap \delta]$ by $f$. Its relative
dimension is $n-{\rm dim}(\sigma)+d-1$. Using the smooth base
change theorem and the fact that $f^\ast[n-{\rm dim}(\sigma)+d-1]$
is exact with respect to the perverse t-structure ([BBD] page
108-109), one can prove the assertion about $j_{!\ast}$.

\bigskip
\bigskip
\centerline{\bf 3. Compactification of $f:{\bf T}_k^n\to {\bf
A}_k^1$}

\bigskip
\bigskip
In this section, we take $N$ to be the lattice ${\bf Z}^n$ and
$V={\bf R}^n$. Let
$$G=\sum_{i\in {\bf Z}^n} a_iX^i \in A[X_1,X_1^{-1},\ldots, X_n,
X_n^{-1}]$$ be a Laurent polynomial with coefficients in a
commutative ring $A$. The {\it Newton polyhedron} $\Delta(G)$ of
$G$ is the convex hull in ${\bf R}^n$ of the set $\{i\in {\bf
Z}^n|a_i\not=0\}$. For any face $\tau$ of $\Delta(G)$, set
$$G_\tau=\sum_{i\in\tau} a_iX^i.$$ Assume ${\rm
dim}(\Delta(G))=n$, and let $\Sigma$ be a fan which is a
subdivision of the fan $\Sigma(\Delta(G))$ in Proposition 1.2. Let
$Y$ be the scheme theoretic closure in the toric scheme
$X_A(\Sigma)$ of the locus $G=0$ in the torus ${\bf T}_A^n={\rm
Spec} A[X_1, X_1^{-1},\ldots,X_n, X_n^{-1}]$. Here we regard the
torus ${\bf T}_A^n$ as an open subscheme of $X_A(\Sigma)$. Note
that $X_A(\Sigma)$ and $Y$ are proper over $A$.

\bigskip
\noindent {\bf Proposition 3.1.} Notation as above. Suppose $A$ is
an integral domain. Let $\sigma\in\Sigma$ and let
$\tau=F_{\Delta(G)}(\sigma)$. Then $$(\tau\cap \{i\in {\bf
Z}^n|a_i\not= 0\})\not=\emptyset.$$ Let $P\in \tau\cap \{i\in {\bf
Z}^n|a_i\not= 0\}$. Then we have $X^{-P}G\in A[{\bf Z}^n\cap
\check\sigma]$ and $X^{-P}G_\tau \in A[{\bf
Z}^n\cap\sigma^\perp]$. Moreover, on the open subscheme
$U_\sigma={\rm Spec}\,A[{\bf Z}^n\cap \check \sigma]$ of
$X_A(\Sigma)$, $Y\cap U_\sigma$ is the locus of $X^{-P}G=0$, and
on the torus $O_\sigma={\rm Spec}\, A [{\bf Z}^n\cap
\sigma^\perp]$ (regarded as a subscheme of $X_A(\Sigma)$), $Y\cap
O_\sigma$ is the locus of $X^{-P}G_\tau=0$.

\bigskip
\noindent {\bf Proof.} Since $\Sigma$ is a subdivision of
$\Sigma(\Delta(G))$, we can choose $\sigma'\in \Sigma(\Delta(G))$
such that $\sigma\subset \sigma'$. By the construction in
Proposition 1.2, $F_{\Delta(G)}(\sigma')$ is not empty. As
$\tau=F_{\Delta(G)}(\sigma)\supset F_{\Delta(G)}(\sigma')$, $\tau$
is nonempty. Since $\tau$ is a face of $\Delta(G)$, we must have
$(\tau\cap \{i\in {\bf Z}^n|a_i\not= 0\})\not=\emptyset$ by the
definition of $\Delta(G)$.

For any $v\in \sigma$, we have
$$F_{\Delta(G)}(v)\supset F_{\Delta(G)}(\sigma)\supset\tau.$$ So
$\langle u', v\rangle \geq \langle u, v\rangle $ for any $u'\in
\Delta(G)$ and $u\in \tau$. Hence $v\in ({\rm
cone}_{\Delta(G)}(\tau))^\vee$. So $\sigma\subset({\rm
cone}_{\Delta(G)}(\tau))^\vee$ and hence ${\rm
cone}_{\Delta(G)}(\tau)\subset \check \sigma$. In particular, we
have $X^{-P}G\in A[{\bf Z}^n\cap \check \sigma]$. Since
$$\tau-\tau\subset ({\rm
cone}_{\Delta(G)}(\tau))\cap (-{\rm
cone}_{\Delta(G)}(\tau))\subset
\check\sigma\cap(-\check\sigma)=\sigma^\perp,$$ we have
$X^{-P}G_\tau \in A[{\bf Z}^n\cap\sigma^\perp]$. For any
$Q\in\Delta(G)$ with $Q\not\in\tau=F_{\Delta(G)}(\sigma)$, we have
$Q\not\in F_{\Delta(G)}(v)$ for some $v\in \sigma$. By the
definition of $F_{\Delta(G)}(v)$, we must have $\langle Q,v\rangle
> \langle P,v\rangle$. So $Q-P\not\in \sigma^\perp$. Therefore,
under the epimorphism
$$A[{\bf Z}^n\cup\check\sigma]\to A[{\bf Z}^n\cup\sigma^\perp],\;
\chi^u\mapsto \left\{\begin{array}{cc} \chi^u &\hbox { if } u\in
\sigma^\perp,\\ 0 &\hbox { if
}u\not\in\sigma^\perp,\end{array}\right.$$ $X^{-P}G$ is mapped to
$X^{-P}G_\tau$. As this epimorphism defines the immersion
$O_\sigma \to U_\sigma$, if we can prove $Y\cap U_\sigma$ is the
locus of $X^{-P}G=0$ in $U_\sigma$, then $Y\cap O_\sigma$ is the
locus of $X^{-P}G_\tau=0$ in $O_\sigma$.

Since $Y$ is the scheme theoretic closure in $X_A(\Sigma)$ of the
locus $G=0$ in ${\bf T}_A^n={\rm Spec}\, A[{\bf Z}^n]$, to prove
$Y\cap U_\sigma$ is the locus of $X^{-P}G=0$ in $U_\sigma$, we
need to show the kernel of the composition
$$A[{\bf Z}^n\cup \check \sigma]\to A[{\bf Z}^n]\to A[{\bf
Z}^n]/(G)\eqno (1)$$ is the ideal generated by $X^{-P}G$. Let
$\{v_1,\ldots, v_k\}$ be a minimal family of generators of
$\sigma$, and let $\sigma_i$ $(i=1,\ldots, k)$ be the rays
generated by $v_i$. Note that $\sigma_i$ are one dimensional faces
of $\sigma$. We have $\check \sigma=\bigcap\limits_{i=1}^k \check
\sigma_i$ and hence
$$A[{\bf Z}^n\cap\check \sigma]=\bigcap\limits_{i=1}^k A[{\bf
Z}^n\cap \check \sigma_i].$$ It suffices to show the kernel of the
composition
$$A[{\bf Z}^n\cap \check \sigma_i]\to A[{\bf Z}^n]\to A[{\bf
Z}^n]/(G)\eqno(2)$$ is the ideal generated by $X^{-P}G$ for each
$i$. Indeed, if this is true and if $f$ lies in the kernel of the
composition (1), then $f$ lies in the kernel of the composition
(2) for each $i$. Hence $f$ lies in the ideal of $A[{\bf Z}^n\cap
\check\sigma_i]$ generated by $X^{-P}G$ for each $i$. So
$$\frac{f}{X^{-P}G}\in \bigcap_{i=1}^k A[{\bf Z}^n\cap \check
\sigma_i]=A[{\bf Z}^n\cap \check\sigma],$$ that is, $f$ lies in the
ideal of $A[{\bf Z}^n\cap \check\sigma]$ generated by $X^{-P}G$.

For each $i$, ${\bf Z}^n/({\bf Z}^n\cap {\rm span}(\sigma_i))$ has
no torsion. So ${\bf Z}^n\cap {\rm span}(\sigma_i)$ is a direct
factor of ${\bf Z}^n$. Since $\sigma_i$ is a one-dimensional
rational cone, we can choose a basis ${e_1,\ldots, e_n}$ of ${\bf
Z}^n$ so that $\sigma_i$ is generated by $e_1$. Denote by
$Y_1,\ldots, Y_n$ the coordinates with respect to this basis. We
then have isomorphisms
\begin{eqnarray*}
A[{\bf Z}^n\cap \check \sigma_i]&\cong&
A[Y_1,Y_2,Y_2^{-1},\ldots,Y_n,Y_n^{-1}],\\
A[{\bf Z}^n]&\cong& A[Y_1,Y_1^{-1},Y_2,Y_2^{-1},\ldots,
Y_n,Y_n^{-1}].
\end{eqnarray*}
Through these isomorphisms, the open immersion ${\bf
T}_A^n\hookrightarrow U_{\sigma_i}$  corresponds to the canonical
homomorphism $A[Y_1,Y_2,Y_2^{-1},\ldots,Y_n,Y_n^{-1}]\hookrightarrow
A[Y_1,Y_1^{-1},Y_2,Y_2^{-1},\ldots, Y_n,Y_n^{-1}]$. We need to show
the kernel of the composition
\begin{eqnarray*}
A[Y_1,Y_2,Y_2^{-1},\ldots, Y_n,Y_n^{-1}] &\hookrightarrow &
A[Y_1,Y_1^{-1},Y_2,Y_2^{-1},\ldots, Y_n,Y_n^{-1}]\\&\to&
A[Y_1,Y_1^{-1},Y_2,Y_2^{-1},\ldots, Y_n,Y_n^{-1}]/(G)\\
&=&A[Y_1,Y_1^{-1},Y_2,Y_2^{-1},\ldots, Y_n,Y_n^{-1}]/(X^{-P}G)
\end{eqnarray*}
is the ideal of $A[Y_1,Y_2,Y_2^{-1},\ldots, Y_n,Y_n^{-1}]$ generated
by $X^{-P}G$. Let $B=A[Y_2,Y_2^{-1},\ldots, Y_n,Y_n^{-1}]$. If $f$
lies in the kernel of the above composition, then $f=X^{-P}Gg$ for
some $g\in B[Y_1,Y_1^{-1}]$. Note that both $f$ and $X^{-P}G$ are in
$B[Y_1]$. By the choice of $P$, $X^{-P}G$ has a nonzero constant
term. This implies that $g\in B[Y_1]$. So $f$ lies in the ideal of
$B[Y_1]$ generated by $X^{-P}G$. This finishes the proof of the
proposition.

\bigskip
Let $$f=\sum_{i\in {\bf Z}^n} a_i X^i\in A[X_1,X_1^{-1},\ldots,
X_n, X_n^{-1}]$$  be a Laurent polynomial. Recall that the Newton
polyhedron $\Delta_\infty(f)$ of $f$ at $\infty$ is the convex
hull in ${\bf R}^n$ of the set $\{i\in {\bf Z}^n|a_i\not= 0\}\cup
\{0\}.$ We say $f$ is {\it non-degenerate} with respect to
$\Delta_\infty(f)$ if for any face $\tau$ of $\Delta_\infty(f)$
not containing $0$, the locus of
$$\frac{\partial f_\tau}{\partial X_1}=\cdots=
\frac{\partial f_\tau}{\partial X_n}=0$$ in ${\bf T}_A^n$ is
empty. By [SGA 1] Expos\'e II, Corollaire 4.5, this is equivalent
to saying that the morphism
$$f_\tau:{\bf T}_A^n={\rm Spec}\,A[X_1,X_1^{-1},\ldots,
X_n, X_n^{-1}] \to {\bf A}_A^1={\rm Spec}\, A[T]$$ defined by the
$A$-algebra homomorphism $$A[T]\to A[X_1,X_1^{-1},\ldots, X_n,
X_n^{-1}],\; T\mapsto f_\tau$$ is smooth.

\bigskip
Let $k$ be a field, and let $f=\sum_{i\in {\bf Z}^n} a_i X^i\in
k[X_1,X_1^{-1},\ldots, X_n, X_n^{-1}]$ be a Laurent polynomial with
coefficients in $k$. It defines a $k$-morphism
$$f:{\bf T}_k^n\to {\bf A}_k^1.$$ Let $$A=k[T]$$ and let $$G=f-T.$$
Regard $G$ as a Laurent polynomial over $A$. We have
$$\Delta(G)=\Delta_\infty(f).$$ Suppose ${\rm dim}(\Delta_\infty(f))=n$. Let
$\Sigma$ be a fan that is a subdivision of
$\Sigma(\Delta_\infty(f))$, let $Y$ be the scheme theoretic
closure in $X_A(\Sigma)$ of the locus $G=0$ in ${\bf T}_A^n$, and
let $g:Y\to {\bf A}_k^1$ be the composition $Y\to X_A(\Sigma)\to
{\rm Spec}\, A$. Note that $g$ is proper. The locus of $G=0$ in
${\bf T}_A^n$ is the closed subscheme ${\rm
Spec}\,k[T,X_1,X_1^{-1},\ldots, X_n,X_n^{-1}]/(f-T)$ of ${\bf
T}_A^n={\rm Spec}\,k[T,X_1,X_1^{-1},\ldots, X_n,X_n^{-1}]$. Since
we have an isomorphism
$$k[T,X_1,X_1^{-1},\ldots,
X_n,X_n^{-1}]/(f-T)\stackrel {\cong}\to k[X_1,X_1^{-1},\ldots,
X_n,X_n^{-1}],\; T\mapsto f,$$ the locus of $G=0$ in ${\bf T}_A^n$
can be identified with ${\bf T}_k^n={\rm
Spec}\,k[X_1,X_1^{-1},\ldots, X_n,X_n^{-1}]$ and the restriction
of $g$ to this locus can be identified with $f:{\bf T}_k^n\to {\bf
A}_k^1$. So $g$ is a compactification of $f$. Let $$j_Y:{\bf
T}_k^n\cong Y\cap {\bf T}_A^n\to Y$$ be the open immersion induced
by the immersion of the open dense torus ${\bf T}_A^n$ in
$X_A(\Sigma)$, and let ${\cal K}_\chi$ be a Kummer sheaf on ${\bf
T}_k^n$. In the following, we study the properties of ${j_Y}_\ast
{\cal K}_\chi$ and ${j_Y}_{!\ast}({\cal K}_\chi[n])$. We start
with a lemma.

\bigskip
{\bf Lemma 3.2.} Let $$f=\sum_{i\in {\bf Z}^n} a_i X^i\in
k[X_1,X_1^{-1},\ldots, X_n, X_n^{-1}]$$ be a Laurent polynomial
over a field $k$. Suppose ${\rm dim}(\Delta_\infty(f))=n$. Let
$\Sigma$ be a subdivision of $\Sigma(\Delta_\infty(f))$, $A=k[T]$,
$G=f-T$, $Y$ the scheme theoretic closure in $X_A(\Sigma)$ of the
locus $G=0$ in ${\bf T}_A^n$, and $g:Y\to{\bf A}_k^1$ the
composition $Y\to X_A(\Sigma)\to {\rm Spec}\, A$. Let
$\sigma\in\Sigma$ and let $\tau=F_{\Delta_\infty(f)}(\sigma)$.

(i) $Y\cap O_\sigma$ is not equal to $O_\sigma$ near any point.

(ii) Suppose $0\not\in \tau$. Choose $P\in\tau\cap\{i\in{\bf
Z}^n|a_i\not=0\}$ and let $Z$ be the locus of $X^{-P}f_\tau=0$ in
${O_\sigma}_k={\rm Spec}\, k[{\bf Z}^n\cap \sigma^\perp]$. Then we
have a Cartesian diagram
$$\begin{array}{ccc}
Y\cap O_\sigma &\to & Z \\
g\downarrow&&\downarrow \\
{\bf A}_k^1&\to & {\rm Spec}\, k.
\end{array}$$
If $f$ is non-degenerate, then $Z$ is smooth over $k$, and hence
$Y\cap O_\sigma$ is smooth over ${\bf A}_k^1$, and over $k$.

(iii) Suppose $0\in\tau$. Then we have an isomorphism $$Y\cap
O_\sigma\cong {O_\sigma}_k$$ and $g|_{Y\cap O_\sigma}$ is identified
with the morphism $f_\tau: {O_\sigma}_k\to {\bf A}_k^1$ defined by
the $k$-algebra homomorphism
$$k[T]\to k[{\bf Z}^n\cap
\sigma^\perp], \; T\mapsto f_\tau.$$ In particular, $Y\cap
O_\sigma$ is smooth over $k$. If $f$ is non-degenerate and
$\sigma\in \Sigma(\Delta_\infty(f))$, then ${Y\cap O_\sigma}$ is
smooth over ${\bf A}_k^1$ outside finitely many closed points.

\bigskip
\noindent {\bf Proof.}

(i) Let $P\in\tau\cap(\{i\in{\bf Z}^n|a_i\not=0\}\cup\{0\})$. By
Proposition 3.1, $Y\cap O_\sigma$ is the locus of $X^{-P}G_\tau=0$
in $O_\sigma={\rm Spec}\, A[{\bf Z}^n\cap\sigma^\perp]$. Since
$X^{-P}G_\tau$ is nonzero and $O_\sigma$ is integral, $Y\cap
O_\sigma$ is not equal to $O_\sigma$ near any point.

(ii) Since $0\not\in \tau$, we have $X^{-P}G_\tau=X^{-P}f_\tau$. As
$X^{-P}f_\tau$ does not involve the variable $T$, $Y\cap O_\sigma\to
{\bf A}_k^1$ can be obtained from $Z\to {\rm Spec}\, k$ by base
change.

Suppose $f$ is non-degenerate. Let's prove $Z$ is smooth. Since
${\bf Z}^n/{\bf Z}^n\cap \sigma^\perp$ is torsion free, ${\bf
Z}^n\cap \sigma^\perp$ is a direct factor of ${\bf Z}^n$. Choose a
basis $\{e_1,\ldots, e_n\}$ of ${\bf Z}^n$ so that $\{e_1,\ldots,
e_r\}$ is a basis of ${\bf Z}^n\cap \sigma^\perp$. Let
$(Y_1,\ldots,Y_n)$ be the coordinates with respect to this basis.
Then $X^{-P}f_\tau\in k[{\bf Z}^n\cap \sigma^\perp]$ only depends
on the variables $Y_1,\ldots, Y_r$. Since $f$ is non-degenerate
with respect to $\Delta_\infty(f)$ and $0\not\in\tau$, the locus
of $X^{-P}f_\tau=0$ in ${\bf T}_k^n$ is smooth over $k$. As
$X^{-P}f_\tau$ depends only on $Y_1,\ldots, Y_r$, the locus of
$X^{-P}f_\tau=0$ in ${\bf T}_k^r={\rm Spec}\,k[{\bf Z}^n\cap
\sigma^\perp]$ is smooth over $k$, that is, $Z$ is smooth over
$k$.

(iii) Now suppose $0\in\tau$. We can then take $P=0$. We have
$G_\tau=f_\tau-T$. Since we have an isomorphism
$$A[{\bf Z}^n\cap\sigma^\perp]/(f_\tau-T)\cong k[{\bf
Z}^n\cap\sigma^\perp], \; T\mapsto f_\tau,$$ the locus of $G_\tau=0$
in $O_\sigma={\rm Spec}\, A[{\bf Z}^n\cap \sigma^\perp]$ can be
identified with ${O_\sigma}_k={\rm Spec}\, k[{\bf Z}^n\cap
\sigma^\perp]$, and $g|_{Y\cap O_\sigma}$ can be identified with the
morphism $f_\tau: {O_\sigma}_k\to {\bf A}_k^1$.

Again choose a basis $\{e_1,\ldots, e_n\}$ of ${\bf Z}^n$ so that
$\{e_1,\ldots, e_r\}$ is a basis of ${\bf Z}^n\cap \sigma^\perp$,
where $r={\rm dim}(\sigma^\perp)$. Let $(Y_1,\ldots,Y_n)$ be the
coordinates with respect to this basis. Then $f_\tau$ depend only on
the coordinates $Y_1,\ldots, Y_r$. Note that
$$\Delta_\infty(f_\tau(Y_1,\ldots,Y_r))=\tau.$$
If $f$ is non-degenerate with respect to $\Delta_\infty(f)$, then
$f_\tau(Y_1,\ldots, Y_r)$ is non-degenerate with respect to
$\Delta_\infty(f_\tau(Y_1,\ldots,Y_r))$. Furthermore, if $\sigma\in
\Sigma(\Delta_\infty(f))$, then by Proposition 1.2, we have ${\rm
dim}(\sigma)=n-{\rm dim}(\tau)$, and hence
$${\rm dim}(\Delta_\infty(f_\tau(Y_1,\ldots,Y_r)))={\rm dim}(\tau)=
{\rm dim}(\sigma^\perp)=r.$$ By [DL] 3.5, $f_\tau: {\bf
T}_k^{r}\to {\bf A}_k^1$ is then smooth outside finitely many
closed points. So ${Y\cap O_\sigma}$ is smooth over ${\bf A}_k^1$
outside finitely many closed points under our assumption.

\bigskip
\noindent {\bf Lemma 3.3.} Keep the notations in Lemma 3.2.
Suppose ${\rm dim}(\Delta_\infty(f))=n$,  $f$ is non-degenerate,
and $\Sigma$ is a subdivision of $\Sigma(\Delta_\infty(f))$. Let
$j:{\bf T}_k^n\to X_k(\Sigma)$ be the immersion of the open dense
torus in $X_k(\Sigma)$, and let $(j_{!\ast}({\cal K}_\chi[n]))|_Y$
and $(j_\ast{\cal K}_\chi)|_Y$ the inverse images under the
composition
$$Y\to X_A(\Sigma)\to X_k(\Sigma)$$ of $j_{!\ast}({\cal
K}_\chi[n])$ and $j_\ast{\cal K}_\chi$, respectively. Then we have
\begin{eqnarray*}
{j_Y}_{!\ast}({\cal K}_\chi[n])&\cong&
(j_{!\ast}({\cal K}_\chi[n]))|_Y,\\
{j_Y}_{\ast}{\cal K}_\chi&\cong& (j_{\ast}{\cal K}_\chi)|_Y.
\end{eqnarray*}

\bigskip
\noindent {\bf Proof.} By Lemma 3.2, $Y\cap O_\sigma$ is smooth
over $k$ for any $\sigma\in \Sigma$. Our assertion then follows
from Proposition 2.7.

\bigskip
\noindent {\bf Proposition 3.4.} Keep the notations in Lemma 3.2.
Suppose ${\rm dim}(\Delta_\infty(f))=n$, $f$ is non-degenerate,
and $\Sigma$ is a subdivision of $\Sigma(\Delta_\infty(f))$ with
the property that for any $\sigma \in \Sigma$ with $0\in
F_{\Delta_\infty(f)}(\sigma)$, we have
$\sigma\in\Sigma(\Delta_\infty(f))$.

(i) Outside finitely many closed points in $Y$, $g: Y\to {\bf
A}_k^1$ is universally locally acyclic relative to
${j_Y}_{!\ast}({\cal K}_\chi[n])$, and relative to ${j_Y}_\ast
{\cal K}_\chi$.

(ii) For any $\sigma\in\Sigma$, $g|_{Y\cap V(\sigma)}\to {\bf
A}_k^1$ is universally locally acyclic relative to $({j_Y}_\ast
{\cal K}_\chi)|_{Y\cap V(\sigma)}$ outside finitely many closed
points.

\bigskip
\noindent {\bf Proof.}

(i) By Lemma 3.3, we have
\begin{eqnarray*}
{j_Y}_{!\ast}({\cal K}_\chi[n])&\cong&
(j_{!\ast}({\cal K}_\chi[n]))|_Y,\\
{j_Y}_{\ast}{\cal K}_\chi&\cong& (j_{\ast}{\cal K}_\chi)|_Y,
\end{eqnarray*}
By Lemma 3.2, $Y\cap O_\sigma$ is smooth over $A$ outside finitely
many closed points for any $\sigma\in\Sigma$. So by Proposition 2.5,
$Y$ is universally locally acyclic over $A$ relative to
$(j_{!\ast}({\cal K}_\chi[n]))|_Y$ and relative to $(j_\ast{\cal
K}_\chi)|_Y$ outside finitely many closed points. Thus $g:Y\to {\rm
Spec}\, A={\bf A}_k^1$ is universally locally acyclic relative
${j_Y}_{!\ast}({\cal K}_\chi[n])$ and relative to ${j_Y}_\ast {\cal
K}_\chi$ outside finitely many closed points.

(ii) Let $$p_1:{\bf T}_k^n={\rm Spec}\,k[{\bf Z}^n]\to
{O_\sigma}_k={\rm Spec}\, k[{\bf Z}^n\cap \sigma^\perp]$$ be the
projection. If ${\cal K}_\chi$ is not the inverse image of a Kummer
sheaf on ${O_\sigma}_k$, then by Proposition 2.4, $j_\ast{\cal
K}_\chi$ vanishes on ${V(\sigma)}_k={\rm Spec}\, k[{\bf Z}^n\cap
\check\sigma]$. As ${j_Y}_\ast{\cal K}_\chi=(j_\ast{\cal
K}_\chi)|_Y$, ${j_Y}_\ast{\cal K}_\chi$ vanishes on $Y\cap
V(\sigma)$. Hence $g|_{Y\cap V(\sigma)}\to {\bf A}_k^1$ is
universally locally acyclic relative to $({j_Y}_\ast {\cal
K}_\chi)|_{Y\cap V(\sigma)}$ in this case.

Now suppose ${\cal K}_\chi=p_1^\ast{\cal K}_{\chi_1}$ for some
Kummer sheaf ${\cal K}_{\chi_1}$ on ${O_\sigma}_k$. Let
$\Sigma_1={\rm star}(\sigma)$ and let $j_1:{O_{\sigma}}_k\to
X_k(\Sigma_1)={V(\sigma)}_k$ be the immersion of the open dense
torus. By Proposition 2.4, we have $$(j_\ast{\cal
K}_\chi)|_{V(\sigma)_k}\cong {j_1}_\ast{\cal K}_{\chi_1}.$$ So we
have
\begin{eqnarray*}
({j_Y}_{\ast}{\cal K}_\chi)|_{Y\cap V(\sigma)}&\cong&
(j_{\ast}{\cal K}_\chi)|_{Y\cap V(\sigma)}\\
&\cong& ({j_1}_\ast{\cal K}_{\chi_1})|_{Y\cap V(\sigma)}.
\end{eqnarray*}
For each $\tau\in\Sigma_1$, by Lemma 3.2, $Y\cap O_\tau$ is smooth
over $A$ outside finitely many points. So by Proposition 2.5
applied to the toric scheme $X_A(\Sigma_1)=V(\sigma)$ and the
Cartier divisor $Y\cap V(\sigma)$, $Y\cap V(\sigma)$ is
universally locally acyclic over $A$ relative to $({j_1}_\ast{\cal
K}_{\chi_1})|_{Y\cap V(\sigma)}$ outside finitely many closed
points. Hence ${Y\cap V(\sigma)}$ is universally locally acyclic
over $A$ relative to $({j_Y}_\ast {\cal K}_\chi)|_{Y\cap
V(\sigma)}$ outside finitely many closed points.

\bigskip
\noindent {\bf Lemma 3.5.} Let $Y$ be a scheme of finite type over
a field $k$, ${\cal F}$ a $\overline {\bf Q}_l$-sheaf on $Y$, $J$
a finite set, and $Y_j$ a closed subscheme of $Y$ for each $j\in
J$. For any $I\subset J$, let
$$Y_I=\bigcap_{j\in I}Y_j,\; Y_I^\circ=Y_I-\bigcup_{j\in J-I} Y_j.$$
We have $Y_\emptyset=Y$, and we denote
$Y_\emptyset^\circ=Y-\bigcup_{j\in J}Y_j$ by $Y^\circ$.

(i) Let $g:Y\to {\bf A}_k^1$ be a $k$-morphism. Suppose
$R^i(g|_{Y_I^\circ})_!{\cal F}$ are tame at $\infty$ for all
$I\subset J$ and all $i$. Then $R^i(g|_{Y_I})_!{\cal F}$ are tame at
$\infty$ for all $I\subset J$ and all $i$. In particular,
$R^ig_!{\cal F}$ are tame at $\infty$ for all $i$.

(ii) Conversely, suppose $R^i(g|_{Y_I})_!{\cal F}$ are tame at
$\infty$ for all $I\subset J$ and all $i$. Then
$R^i(g|_{Y_I^\circ})_!{\cal F}$ are tame at $\infty$ for all
$I\subset J$ and all $i$. In particular, $R^i(g|_{Y^\circ})_!{\cal
F}$ are tame at $\infty$ for all $i$.

(iii) Let $\phi$ be an integer valued function on the power set of
$J$ with the property $\phi(I')\leq \phi(I)$ for any $I\subset
I'\subset J$. Suppose $H_c^i(Y_I^\circ\otimes_k\overline k,{\cal
F})=0$ for all $I\subset J$ and all $i>\phi(I)$. Then
$H_c^i(Y_I\otimes_k\overline k,{\cal F})=0$ for all $I\subset J$
and all $i>\phi(I)$. In particular, $H_c^i(Y\otimes_k\overline
k,{\cal F})=0$ for all $i>\phi(\emptyset)$.

(iv) Let $\phi$ be as in (iii). Suppose furthermore that for any
$I\subset I'\subset J$ with $Y_{I'}\not=\emptyset, Y_{I}$, we have
$\phi(I')<\phi(I)$. If $H_c^i(Y_I\otimes_k\overline k,{\cal F})=0$
for all $I\subset J$ and all $i>\phi(I)$, then
$H_c^i(Y_I^\circ\otimes_k\overline k,{\cal F})=0$ for all
$I\subset J$ and all $i>\phi(I)$. In particular,
$H_c^i(Y^\circ\otimes_k\overline k,{\cal F})=0$ for all
$i>\phi(\emptyset)$.

\bigskip
\noindent {\bf Proof.} We will prove (i) and (iv). The proof of
(ii) and (iii) is similar.

(i) It suffices to prove $R^ig_!{\cal F}$ are tame at $\infty$ for
all $i$. Indeed, applying this result to $g|_{Y_I}$ and the closed
subschemes $Y_I\cap Y_j$ $(j\in J-I)$, we see $R^i(g|_{Y_I})_!{\cal
F}$ are tame at $\infty$ for all $i$ and all $I\subset J$. We use
induction on the number of elements of $J$. The case where
$J=\emptyset$ is trivial. If $J=\{1\}$, we have $Y_\emptyset^\circ
=Y-Y_1$ and $Y_J^\circ=Y_1$. By our assumption,
$R^i(g|_{Y-Y_1})_!{\cal F}$ and $R^i(g|_{Y_1})_!{\cal F}$ are tame
at $\infty$. We have a long exact sequence
$$\cdots\to R^i(g|_{Y-Y_1})_!{\cal F}\to
R^ig_!{\cal F}\to R^i(g|_{Y_1})_!{\cal F}\to\cdots.$$ It follows
that $R^ig_!{\cal F}$ are tame at $\infty$. Suppose $J=\{1,\ldots,
m\}$ $(m\geq 2)$ and suppose our assertion holds for those $J$ with
less than $m$ elements. Applying the induction hypothesis to $Y-Y_1$
and the closed subschemes $Y_j-Y_1$ $(j\in J-\{1\})$, we see
$R^i(g|_{Y-Y_1})_!{\cal F}$ are tame at $\infty$. Applying the
induction hypothesis to $Y_1$ and the closed subschemes $Y_1\cap
Y_j$ $(j\in J-\{1\})$, we see $R^i(g|_{Y_1})_!{\cal F}$ are tame at
$\infty$. It follows from the previous long exact sequence that
$R^ig_!{\cal F}$ are tame at $\infty$.

(iv) As in (i), it suffices to show
$H_c^i(Y^\circ\otimes_k\overline k,{\cal F})=0$ for
$i>\phi(\emptyset)$. We use induction on the number of element of
$J$. The case where $J=\emptyset$ is trivial. If $J=\{1\}$, we
have $Y_\emptyset=Y$, $Y_J=Y_1$ and $Y^\circ=Y-Y_1$. The cases
where $Y_J=\emptyset$ or $Y_J=Y_\emptyset$ are trivial. Suppose
$Y_J\not=\emptyset,Y_\emptyset$. Then $\phi(J)<\phi(\emptyset)$.
We have a long exact sequence
$$\cdots\to H_c^{i-1}(Y_1\otimes_k\overline k, {\cal F})\to
H_c^i((Y-Y_1)\otimes_k\overline k,{\cal F})\to
H_c^i(Y\otimes_k\overline k,{\cal F})\to\cdots.$$ By our
assumption, we have $H_c^{i-1}(Y_1\otimes_k\overline k, {\cal
F})=0$ for $i-1>\phi(J)$ and $H_c^i(Y\otimes_k\overline k,{\cal
F})=0$ for $i>\phi(\emptyset)$. It follows that
$H_c^i((Y-Y_1)\otimes_k\overline k,{\cal F})=0$ for
$i>\phi(\emptyset)$. Let $J=\{1,\ldots, m\}$ $(m\geq 2)$ and
suppose our assertion holds for those $J$ with less than $m$
elements. Applying the induction hypothesis to $Y$ and the closed
subschemes $Y_j$ $(j\in J-\{1\})$, we see
$H^i((Y-\bigcup_{j\not=1}Y_j)\otimes_k\overline k, {\cal F})=0$
for all $i\geq \phi(\emptyset)$. Applying the induction hypothesis
to $Y_1$ and the closed subschemes $Y_1\cap Y_j$ $(j\in J-\{1\})$,
we see $H^{i-1}((Y_1-\bigcup_{j\not=1}Y_j)\otimes_k\overline k,
{\cal F})=0$ for all $i-1\geq \phi(\{1\})$. We have a long exact
sequence
$$\cdots\to H^{i-1}((Y_1-\bigcup_{j\not=1}Y_j)\otimes_k\overline k,
{\cal F})\to  H^i((Y-\bigcup_j Y_j)\otimes_k\overline k, {\cal
F})\to H^i((Y-\bigcup_{j\not=1}Y_j)\otimes_k\overline k, {\cal
F})\to\cdots.$$ If $Y_1\not=\emptyset, Y$, then
$\phi(\{1\})<\phi(\emptyset)$. It follows that $H^i((Y-\bigcup_j
Y_j)\otimes_k\overline k, {\cal F})=0$ for all
$i>\phi(\emptyset)$. If $Y_1=\emptyset$, then $H^i((Y-\bigcup_j
Y_j)\otimes_k\overline k, {\cal F})=H^i((Y-\bigcup_{j\not=1}
Y_j)\otimes_k\overline k, {\cal F})=0$ for $i\geq
\phi(\emptyset)$. If $Y_1=Y$, then $Y-\bigcup_j Y_j=\emptyset$ and
$H^i((Y-\bigcup_j Y_j)\otimes_k\overline k, {\cal F})=0$ for all
$i$.

\bigskip
\noindent {\bf Lemma 3.6.} Let $k$ be a field. Suppose $f:{\bf
T}_k^n\to {\bf A}_k^1$ is a $k$-morphism defined by a Laurent
polynomial $f\in k[X_1,X_1^{-1},\ldots, X_n,X_n^{-1}]$ that is
non-degenerate with respect to $\Delta_\infty(f)$. Then
$R^if_!{\cal K}_\chi$ are tame at $\infty$ for all $i$.

\bigskip
\noindent {\bf Proof.} For the Kummer covering
$$[m]:{\bf T}_k^n\to {\bf T}_k^n, \; x\mapsto x^m,$$ we have
$$[m]_\ast\overline {\bf Q}_l\cong\bigoplus_{\chi}{\cal K}_\chi,$$
where $\chi:\mu_m(k)^n\to \overline {\bf Q}_l^\ast$ goes over the
set of characters of $\mu_m(k)^n$. The composition $f\circ [m]$ is
defined by the Laurent polynomial $f(x_1^m, \ldots, x_n^m)$. It is
not hard to see that this Laurent polynomial is non-degenerate
with respect to its Newton polyhedron at $\infty$. So by [DL] 4.2,
$R^i(f\circ [m])_!\overline {\bf Q}_l$ are tame at $\infty$ for
all $i$. We have
\begin{eqnarray*}
R^i(f\circ [m])_!\overline {\bf Q}_l&\cong&
R^if_!([m]_\ast\overline {\bf Q}_l) \\
&\cong & \bigoplus_{\chi} R^if_! {\cal K}_\chi.
\end{eqnarray*}
So $R^if_!{\cal K}_\chi$ are tame at $\infty$ for all $i$ and all
$\chi$.

\bigskip
\noindent {\bf Lemma 3.7.} Keep the notations in Lemma 3.2.
Suppose $f$ is non-degenerate, ${\rm dim}(\Delta_\infty(f))=n$,
$\Sigma$ is a subdivision of $\Sigma(\Delta_\infty(f))$, and
$\sigma\in\Sigma$ so that $0\not \in
F_{\Delta_\infty(f)}(\sigma)$. Then $R^i(g|_{Y\cap
O_\sigma})_!({j_Y}_{!\ast}({\cal K}_\chi[n]))$ and $R^i(g|_{Y\cap
O_\sigma})_!({j_Y}_{\ast}{\cal K}_\chi)$ are constant sheaves.

\bigskip
\noindent {\bf Proof.} By Lemma 3.2 (ii), we have a commutative
diagram of Cartesian squares
$$\begin{array}{ccc}
Y\cap O_\sigma&\to& Z\\
\downarrow &&\downarrow \\
X_A(\Sigma)&\to & X_k(\Sigma) \\
\downarrow &&\downarrow \\
{\bf A}_k^1&\to &{\rm Spec}\, k,
\end{array}$$
where $Z$ is a closed subscheme of ${O_\sigma}_k$ and we regard
${O_\sigma}_k$ as a subscheme of $X_k(\Sigma)$. Let $j:{\bf
T}_k^n\to X_k(\Sigma)$ be the immersion of the open dense torus in
$X_k(\Sigma)$. By Lemma 3.3, we have
\begin{eqnarray*}
{j_Y}_{!\ast}({\cal K}_\chi[n])&\cong&
(j_{!\ast}({\cal K}_\chi[n]))|_Y,\\
{j_Y}_{\ast}({\cal K}_\chi)&\cong& (j_{\ast}{\cal K}_\chi)|_Y.
\end{eqnarray*}
Fix notations by the following diagram
$$\begin{array}{rcl}
Y\cap O_\sigma&\stackrel {\rho'}\to &Z\\
g \downarrow &&\downarrow \lambda\\
{\bf A}_k^1&\stackrel{\rho} \to& {\rm Spec}\, k,
\end{array}$$
By the proper base change theorem, we have
\begin{eqnarray*}
R^i(g|_{Y\cap O_\sigma})_!({j_Y}_\ast {\cal K}_\chi)&\cong&
R^i(g|_{Y\cap
O_\sigma})_!((j_{\ast}{\cal K}_\chi)|_{Y\cap O_\sigma}) \\
&\cong& R^i(g|_{Y\cap
O_\sigma})_!{\rho'}^\ast((j_{\ast}{\cal K}_\chi)|_Z) \\
&\cong& \rho^\ast R^i \lambda_! ((j_{\ast}{\cal K}_\chi)|_Z),
\end{eqnarray*}
that is,
$$R^i(g|_{Y\cap O_\sigma})_!({j_Y}_\ast {\cal K}_\chi)\cong
\rho^\ast R^i \lambda_! ((j_{\ast}{\cal K}_\chi)|_Z).$$ Hence
$R^i(g|_{Y\cap O_\sigma})_!({j_Y}_\ast {\cal K}_\chi)$ are constant
for all $i$. Similarly, $R^i(g|_{Y\cap
O_\sigma})_!({j_Y}_{!\ast}({\cal K}_\chi[n]))$ are also constant for
all $i$.

\bigskip
\noindent {\bf Proposition  3.8.} Keep the notations in Lemma 3.2.
Suppose ${\rm dim}(\Delta_\infty(f))=n$,  $f$ is non-degenerate, and
$\Sigma$ is a subdivision of $\Sigma(\Delta_\infty(f))$. For any
$\sigma\in \Sigma$, $R^i(g|_{Y\cap V(\sigma)})_\ast
({j_Y}_{!\ast}({\cal K}_\chi[n]))$ and $R^i(g|_{Y\cap
V(\sigma)})_\ast ({j_Y}_{\ast}{\cal K}_\chi)$ are tame at $\infty$
for all $i$. In particular, taking $\sigma=0$, we see $R^ig_\ast
({j_Y}_{!\ast}({\cal K}_\chi[n]))$ and $R^ig_\ast({j_Y}_{\ast}{\cal
K}_\chi)$ are tame at $\infty$ for all $i$.

\bigskip
\noindent {\bf Proof.} Let $J$ be the set of one dimensional cones
in $\Sigma$, and for each $\sigma\in J$, let $Y_\sigma=Y\cap
V(\sigma)$. For any subset $I=\{\sigma_1,\ldots, \sigma_r\}$ of
$J$, we have
\begin{eqnarray*}
Y_I&=&\bigcap_{i=1}^r Y_{\sigma_i}\\
&=&Y\bigcap \biggl(\bigcap_{i=1}^r V(\sigma_i)\biggr)\\
&=& Y\bigcap \biggl(\bigcap_{i=1}^r\coprod_{\sigma_i\prec \gamma}
O_\gamma\biggr)\\
&=& Y\bigcap \biggl(\coprod_{\sigma_1,\ldots,
\sigma_r\prec\gamma}O_\gamma\biggr).
\end{eqnarray*}
So if $\sigma$ is the smallest cone in $\Sigma$ containing
$\sigma_1,\ldots, \sigma_r$, then we have
$$Y_I=Y\cap V(\sigma).$$
If there is no cone in $\Sigma$ containing all $\sigma_1,\ldots,
\sigma_r$, then $Y_I=\emptyset$. Note that since $\sigma_1,\ldots,
\sigma_r$ are one-dimensional cones in $\Sigma$, the smallest cone
in $\Sigma$ containing $\sigma_1,\ldots, \sigma_r$ is exactly the
cone generated by $\sigma_1,\ldots, \sigma_r$ if this cone lies in
$\Sigma$. Suppose $J-I=\{\sigma_{r+1},\ldots,\sigma_s\}$. Then we
have
\begin{eqnarray*}
Y_I^\circ&=&Y_I-\bigcup_{i=r+1}^s Y_{\sigma_j}\\
&=& Y\bigcap\biggl(\coprod_{\sigma_1,\ldots, \sigma_r\prec
\gamma}O_\gamma-\bigcup_{i=r+1}^s
\coprod_{\sigma_j\prec\gamma}O_\gamma\biggr)\\
&=& Y\bigcap\biggl(\coprod_{\begin{array}{c}
\sigma_1,\ldots, \sigma_r\prec\gamma\\
\sigma_{r+1},\ldots, \sigma_s\not\prec\gamma\end{array}}O_\gamma\biggr)\\
&=& Y\cap O_\sigma,
\end{eqnarray*}
where $\sigma$ is again the cone in $\Sigma$ generated by
$\sigma_1,\ldots, \sigma_r$. (If the cone generated by
$\sigma_1,\ldots, \sigma_r$ does not lie in $\Sigma$, then
$Y_I^\circ=\emptyset$.) By Lemma 3.5 (i), to prove our assertion, it
suffices to show that for any $\sigma\in \Sigma$, $R^i(g|_{Y\cap
O_\sigma})_\ast ({j_Y}_{!\ast}({\cal K}_\chi[n]))$ and
$R^i(g|_{Y\cap O_\sigma})_\ast ({j_Y}_{\ast}{\cal K}_\chi)$ are tame
at $\infty$ for all $i$. If $0\not\in F_{\Delta_\infty(f)}(\sigma)$,
this follows from Lemma 3.7. Suppose $0\in
F_{\Delta_\infty(f)}(\sigma)$. By Lemma 3.2 (iii), $Y\cap O_\sigma$
can be identified with ${O_\sigma}_k$ and $g|_{Y\cap O_\sigma}:Y\cap
O_\sigma\to {\bf A}_k^1$ can be identified with
$f_\tau:{O_\sigma}_k\to {\bf A}_k^1$, where $\tau=
F_{\Delta_\infty(f)}(\sigma)$. Using Lemma 3.3, one can check
$({j_Y}_{!\ast}({\cal K}_\chi[n]))|_{Y\cap O_\sigma}$ and
$({j_Y}_{!\ast}{\cal K}_\chi)|_{Y\cap O_\sigma}$ are identified with
$(j_{!\ast}({\cal K}_\chi[n]))|_{{O_\sigma}_k}$ and $(j_\ast{\cal
K}_\chi)|_{{O_\sigma}_k}$, respectively. So it suffice to verify
$R^i{f_\tau}_!((j_{!\ast}({\cal K}_\chi[n]))|_{{O_\sigma}_k})$ and
$R^i{f_\tau}_!((j_\ast{\cal K}_\chi)|_{{O_\sigma}_k})$ are tame at
$\infty$ for all $i$.  This follows from Lemma 3.6 and the
description of $(j_{!\ast}({\cal K}_\chi[n]))|_{{O_\sigma}_k}$ and
$(j_\ast{\cal K}_\chi)|_{{O_\sigma}_k}$ in Lemma 2.3. This finishes
the proof of the proposition.

\bigskip
A Laurent polynomial
$$G(X_1,\ldots, X_n)=\sum_{i\in {\bf Z}^n} a_iX^i\in k[X_1, X_1^{-1},
\ldots, X_n,X_n^{-1}]$$ is called {\it $0$-non-degenerate} with
respect to its Newton polyhedron $\Delta(G)$ if for any face
$\tau$ of $\Delta(G)$, the locus of $G_\tau=\sum_{i\in
\tau}a_iX^i=0$ in ${\bf T}_k^n$ is smooth over $k$.

\bigskip
\noindent {\bf Proposition 3.9.} Let $G(X_1,\ldots, X_n)\in k[X_1,
X_1^{-1}, \ldots, X_n,X_n^{-1}]$ be a Laurent polynomial that is
0-non-degenerate with respect to $\Delta(G)$. Then
$$\chi_c({\bf T}_{\bar k}^n\cap G^{-1}(0), {\cal
K}_\chi)=(-1)^{n-1}n!{\rm vol}(\Delta(G)),$$ where
$$\chi_c({\bf T}_{\bar k}^n\cap G^{-1}(0), {\cal
K}_\chi)=\sum_i(-1)^i {\rm dim}H_c^i({\bf T}_{\bar k}^n\cap
G^{-1}(0), {\cal K}_\chi)$$ is the Euler characteristic.

\bigskip
\noindent {\bf Proof.} Let $r={\rm dim}(\Delta(G))$. If $r<n$,
then after a suitable change of coordinates on the torus, we may
assume $G$ depends on $r$ variables. Without loss of generality,
assume
$$G(X_1,\ldots, X_n)=H(X_1,\ldots, X_r)$$ for some Laurent
polynomial $H(X_1,\ldots, X_r)$. Then we have $${\bf T}_{\bar
k}^n\cap G^{-1}(0)=({\bf T}_{\bar k}^r\cap H^{-1}(0))\times_{\bar
k} {\bf T}_{\bar k}^{n-r}.$$ Let $$p_1:{\bf T}_k^n={\bf
T}_k^r\times {\bf T}_k^{n-r}\to {\bf T}_k^r,\; p_2:{\bf
T}_k^n={\bf T}_k^r\times {\bf T}_k^{n-r}\to {\bf T}_k^{n-r}$$ be
the projections. We can find Kummer sheaves ${\cal K}_{\chi_1}$ on
${\bf T}_k^r$ and ${\cal K}_{\chi_2}$ on ${\bf T}_k^{n-r}$ so that
$${\cal K}_\chi\cong p_1^\ast{\cal K}_{\chi_1}\otimes
p_2^\ast{\cal K}_{\chi_2}.$$ By the K\"{u}nneth formula, we have
\begin{eqnarray*}
\chi_c({\bf T}_{\bar k}^n\cap G^{-1}(0), {\cal K}_\chi)&=& \chi_c
(({\bf T}_{\bar k}^r\cap H^{-1}(0))\times_{\bar k} {\bf T}_{\bar
k}^{n-r},p_1^\ast{\cal K}_{\chi_1}\otimes p_2^\ast{\cal K}_{\chi_2})\\
&=& \chi_c ({\bf T}_{\bar k}^r\cap H^{-1}(0), {\cal K}_{\chi_1})
\chi_c ({\bf T}_{\bar k}^{n-r},{\cal K}_{\chi_2}).
\end{eqnarray*}
We have $$\chi_c ({\bf T}_{\bar k}^{n-r},{\cal K}_{\chi_2})=0.$$
(To see this, we use the K\"unneth formula to reduce to the case
where the dimension of the torus is 1. We then use the
Grothendieck-Ogg-Shafarevich formula.) It follows that
$$\chi_c({\bf T}_{\bar k}^n\cap G^{-1}(0), {\cal
K}_\chi)=0=(-1)^{n-1}n!{\rm vol}(\Delta(G)).$$

Now suppose ${\rm dim}\Delta(G)=n$. Let $\Sigma$ be a regular fan
that is a subdivision of $\Sigma(\Delta(G))$. Let $Y$ be the scheme
theoretic closure in $X_k(\Sigma)$ of the locus $G=0$ in ${\bf
T}_k^n$. By Proposition 3.1 and the assumption that $G$ is
0-non-degenerate with respect to $\Delta(G)$, $Y\cap O_\sigma$ is
smooth over $k$ for any $\sigma\in \Sigma$. By [DL] 2.3, $Y$ is
smooth over $k$. So $Y$ is a smooth compactification of the locus
$G=0$ in ${\bf T}_k^n$. Note that $({\cal K}_\chi)|_{{\bf T}_k^n\cap
G^{-1}(0)}$ is tamely ramified along $Y-({\bf T}_k^n\cap G^{-1}(0))$
in the sense of [I1] 2.6. (Indeed, the inverse image of ${\cal
K}_\chi$ under the Kummer covering
$$[m]:{\bf T}_k^n\to {\bf T}_k^n,\; x\mapsto x^m$$
is constant.) So by [I1] 2.7, we have
$$\chi_c({\bf T}_{\bar k}^n\cap G^{-1}(0), {\cal
K}_\chi)=\chi_c({\bf T}_{\bar k}^n\cap G^{-1}(0), \overline {\bf
Q}_l).$$ (To apply [I1] 2.7, we only require $Y$ is a normal
compactification of the locus $G=0$ in ${\bf T}_k^n$. If we take
$\Sigma$ to be just an arbitrary subdivision of
$\Sigma(\Delta(G))$, then using [DL] 2.3, one can verify $Y$ is
normal. So $\Sigma$ being regular is not necessary for our
purpose.) By [DL] 2.7, we have
$$\chi_c({\bf T}_{\bar k}^n\cap G^{-1}(0), \overline {\bf
Q}_l)=(-1)^{n-1}n!{\rm vol}(\Delta(G)).$$ So we have
$$\chi_c({\bf T}_{\bar k}^n\cap G^{-1}(0), {\cal
K}_\chi)=(-1)^{n-1}n!{\rm vol}(\Delta(G)).$$

\bigskip
\bigskip
\centerline{\bf 4. Cohomology and weights}

\bigskip
\bigskip
In this section, $k$ is a finite field of characteristic $p$ with
$q$ elements, $\psi:(k,+)\to \overline{\bf Q}_l^\ast$ a nontrivial
additive character, and ${\cal L}_\psi$ the lisse sheaf of rank 1 on
${\bf A}_k^1$ obtained by pushing-forward the ${\bf
A}_k^1(k)$-torsor
$$0\to {\bf A}_k^1(k) \to {\bf A}_k^1\stackrel {\cal P}\to {\bf A}_k^1\to
0$$ by $\psi^{-1}$, where $${\cal P}:{\bf A}_k^1\to {\bf A}_k^1,\;
x\mapsto x^q-x$$ is the Artin-Schreier covering.

\bigskip
The following lemma is essentially Propositions 3.1 and 7.1 in [DL].
We include its proof for completeness.

\bigskip
\noindent {\bf Lemma 4.1.} Let $Y$ be a scheme of finite type over
$k$, let $g:Y\to {\bf A}_k^1$ a proper $k$-morphism, and let
${\cal K}$ be an object in the derived category
$D_c^b(Y,\overline{\bf Q}_l)$ of $\overline {\bf Q}_l$-sheaves on
$Y$ defined in [D] 1.1.2. Suppose $R^i g_\ast {\cal K}$ are tame
at $\infty$ for all $i$, and suppose $g$ is locally acyclic
relative to ${\cal K}$ outside finitely many closed points.

(i) The canonical homomorphisms
$$H_c^i (Y\otimes_k\bar k,
{\cal K}\otimes g^\ast{\cal L}_\psi)\to H^i (Y\otimes_k\bar k, {\cal
K}\otimes g^\ast{\cal L}_\psi)$$ are isomorphisms for all $i$.

(ii) If ${\cal K}$ is perverse, then $$H_c^i (Y\otimes_k\bar k,
{\cal K}\otimes g^\ast{\cal L}_\psi)=0$$ for all $i>0$.

(iii) Suppose $Y$ is pure of dimension $n$ and ${\cal K}$ is a
$\overline {\bf Q}_l$-sheaf. Then $$H_c^i(Y\otimes_k\bar k, {\cal
K}\otimes g^\ast{\cal L}_\psi)=0$$ for all $i>n$.

\bigskip
\noindent {\bf Proof.}

(i) Let $\iota: {\bf A}_k^1\hookrightarrow {\bf P}_k^1$ be the
canonical open immersion. Since $R^jg_\ast{\cal K}$ are tame at
$\infty$ and ${\cal L}_\psi$ is totally wild at $\infty$, we have
$$(\iota_\ast (R^jg_\ast{\cal K}\otimes {\cal L}_\psi))_{\overline \infty}=0.$$
By [SGA $4\frac{1}{2}$] [Sommes trig.] Proposition 1.19 and Exemple
1.19.1, the canonical homomorphisms
$$H^i_c({\bf A}_{\bar k}^1, R^jg_\ast{\cal K}\otimes{\cal L}_\psi)\to
H^i({\bf A}_{\bar k}^1, R^jg_\ast{\cal K}\otimes{\cal L}_\psi)$$
are isomorphisms. We have spectral sequences
\begin{eqnarray*}
E_2^{ij}=H^i_c({\bf A}_{\bar k}^1, R^jg_\ast{\cal K}\otimes{\cal
L}_\psi)\cong H^i_c({\bf A}_{\bar k}^1, R^jg_\ast({\cal K}\otimes
g^\ast {\cal L}_\psi))&\Rightarrow& H_c^{i+j}(Y\otimes_k\bar k,
{\cal K}\otimes g^\ast{\cal L}_\psi),\\
E_2^{ij}=H^i({\bf A}_{\bar k}^1, R^jg_\ast{\cal K}\otimes{\cal
L}_\psi)\cong H^i({\bf A}_{\bar k}^1, R^jg_\ast({\cal K}\otimes
g^\ast {\cal L}_\psi))&\Rightarrow& H^{i+j}(Y\otimes_k\bar k,
{\cal K}\otimes g^\ast{\cal L}_\psi),
\end{eqnarray*}
where to get the first spectral sequence, we use the assumption that
$g$ is proper. It follows that the canonical homomorphisms
$$H_c^i (Y\otimes_k\bar k,
{\cal K}\otimes g^\ast{\cal L}_\psi)\to H^i (Y\otimes_k\bar k,
{\cal K}\otimes g^\ast{\cal L}_\psi)$$ are isomorphisms.

(ii) Let $s$ be an arbitrary closed point in ${\bf A}_k^1$, $S$
the henselization of ${\bf A}_k^1$ at $s$, and $\eta$ the generic
point of $S$. Since $g$ is proper, we have a long exact sequence
$$\begin{array}{ccccccccc} \cdots&\to& H^j(X_{\bar s}, {\cal K})&\to&
H^j(X_{\bar\eta},{\cal K})&\to&
H^j(X_{\bar s}, R\Phi ({\cal K}))&\to&\cdots,\\
&& \cong \uparrow &&\cong\uparrow &&&&\\
&& (R^jg_\ast {\cal K})_{\bar s} && (R^jg_\ast {\cal
K})_{\bar\eta}&&&&
\end{array}$$
where $R\Phi$ is the vanishing cycle functor. Since $g$ is locally
acyclic relative to ${\cal K}$ outside finitely many closed points,
$R\Phi ({\cal K})$ is supported on finitely many closed points. By
[I2] Corollaire 4.6, $R\Phi ({\cal K})[-1]$ is perverse. It follows
that
$$H^j(X_{\bar s}, R\Phi ({\cal K}))=\bigoplus_{x\in {\rm supp}(R\Phi ({\cal K}))}
(R^{j+1}\Phi ({\cal K})[-1])_{\bar x}=0$$ for all $j\geq 0$. The
above long exact sequence then shows that $R^jg_\ast ({\cal K})$
are lisse on ${\bf A}_k^1$ for all $j\geq 1$, and is an extension
of a lisse sheaf by a punctual sheaf for $j=0$. As $R^jg_\ast
{\cal K}$ are tame at $\infty$, all these lisse sheaves are
constant on ${\bf A}_{\bar k}^1$. It follows that $H^i_c({\bf
A}_{\bar k}^1, R^jg_\ast{\cal K}\otimes{\cal L}_\psi)=0$ if $j\geq
1$, and if $j=0$ and $i\geq 1$. As in the proof of (i), the
canonical homomorphisms
$$H^i_c({\bf A}_{\bar k}^1, R^jg_\ast{\cal K}\otimes{\cal L}_\psi)\to
H^i({\bf A}_{\bar k}^1, R^jg_\ast{\cal K}\otimes{\cal L}_\psi)$$
are isomorphisms for all $i$. Combined with the Weak Lefschetz
theorem, we get $H^i_c({\bf A}_{\bar k}^1, R^jg_\ast{\cal
K}\otimes{\cal L}_\psi)=0$ for all $i\geq 2$. It follows that
$H^i_c({\bf A}_{\bar k}^1, R^jg_\ast{\cal K}\otimes{\cal
L}_\psi)=0$ for all $i+j>0$. We have a spectral sequence
$$E_2^{ij}=H^i_c({\bf
A}_{\bar k}^1, R^jg_\ast{\cal K}\otimes{\cal L}_\psi)\Rightarrow
H_c^{i+j}(Y\otimes_k\bar k, {\cal K}\otimes g^\ast{\cal
L}_\psi).$$ So we have $H^i_c(Y\otimes_k\bar k, {\cal K}\otimes
g^\ast{\cal L}_\psi)=0$ for all $i> 0$.

(iii) The proof is similar to that of (ii). Instead of using the
perversity of $R\Phi ({\cal K})[-1]$ in (i), we use the fact that
$R^j\Phi ({\cal K})=0$ for $j>n-1$ ([SGA 7] Expos\'e I, 4.2). We
leave the details to the reader.

\bigskip
We are now ready to prove the main theorem of this paper.

\bigskip
\noindent {\bf Theorem 4.2.} Let $f:{\bf T}_k^n\to {\bf A}_k^1$ be a
$k$-morphism defined by a Laurent polynomial $f\in
k[X_1,X_1^{-1},\ldots, X_n,X_n^{-1}]$ that is non-degenerate with
respect to $\Delta_\infty(f)$ and let ${\cal K}_\chi$ be a Kummer
sheaf on ${\bf T}_k^n$. Suppose ${\rm dim}(\Delta_\infty(f))=n$.
Then

(i) $H_c^i({\bf T}_{\bar k}^n, {\cal K}_\chi\otimes f^\ast{\cal
L}_\psi)=0$ for $i\not=n$.

(ii) ${\rm dim} (H_c^n ({\bf T}_{\bar k}^n, {\cal K}_\chi\otimes
f^\ast{\cal L}_\psi))=n! {\rm vol}(\Delta_\infty(f)).$

(iii)  Let
$$E({\bf T}_k^n, f,\chi)=\sum_{w\in {\bf Z}} e_wT^w,$$
where $e_w$ is the number of eigenvalues with weight $w$ counted
with multiplicities of the geometric Frobenius element  $F$ in
${\rm Gal}(\bar k/k)$ acting on $H_c^i({\bf T}_{\bar k}^n, {\cal
K}_\chi\otimes f^\ast{\cal L}_\psi)$. Then $E({\bf T}_k^n,
f,\chi)$ is a polynomial of degree $\leq n$, and
\begin{eqnarray*}
E({\bf T}_k^n, f,\chi)&=&E(\Delta_\infty(f),\chi),\\
e_n&=&e(\Delta_\infty(f),\chi),
\end{eqnarray*}
where $E(\Delta_\infty(f),\chi)$ and $e(\Delta_\infty(f),\chi)$ are
defined in the introduction.

(iv) If $0$ is an interior point of $\Delta_\infty(f)$, then
$H_c^n({\bf T}_{\bar k}^n, {\cal K}_\chi\otimes f^\ast{\cal
L}_\psi)$ is pure of weight $n$.

\bigskip
\noindent {\bf Proof.}

(i) Let $\Sigma$ be a subdivision of $\Sigma(\Delta_\infty(f))$ such
that for any $\sigma\in \Sigma$ with $0\in
F_{\Delta_\infty(f)}(\sigma)$, we have
$\sigma\in\Sigma(\Delta_\infty(f))$. Let $A=k[T]$, let $Y$ be the
scheme theoretic closure in $X_A(\Sigma)$ of the locus $f-T=0$ in
${\bf T}_A^n$, and let $g: Y\to {\bf A}_k^1$ be the composition
$Y\to X_A(\Sigma)\to {\rm Spec}\, A$. Then $g$ is proper and
$g|_{Y\cap {\bf T}_A^n}$ can be identified with $f:{\bf T}_k^n\to
{\bf A}_k^1$. Let $j_Y:{\bf T}_k^n\cong Y\cap {\bf T}_A^n\to Y$ be
the open immersion induced by the immersion of the open dense torus
${\bf T}_A^n$ in $X_A(\Sigma)$. By Proposition 3.4 (ii), for any
$\sigma\in\Sigma$, $g:Y\cap V(\sigma)\to {\bf A}_k^1$ is locally
acyclic relative to $({j_Y}_\ast{\cal K}_\chi)|_{Y\cap V(\sigma)}$
outside finitely many closed points. By Proposition 3.8,
$R^i(g|_{Y\cap V(\sigma)})_\ast ({j_Y}_{\ast}{\cal K}_\chi)$ are
tame at $\infty$ for all $i$. By Lemma 4.1 (iii), we have
$$H_c^i((Y\cap V(\sigma))\otimes_k \bar k, {j_Y}_\ast{\cal K}_\chi
\otimes g^\ast{\cal L}_\psi)=0$$ for all $i>{\rm dim}(Y\cap
V(\sigma))=n-{\rm dim}(\sigma)$. Let $J$ be the set of nonzero
cones in $\Sigma$. For each $\sigma\in J$, let $Y_\sigma=Y\cap
V(\sigma)$. Then for any $\sigma_1,\ldots, \sigma_r\in J$, we have
\begin{eqnarray*}
Y_{\sigma_1}\cap\cdots\cap Y_{\sigma_r}&=&
Y\bigcap \biggl(\bigcap_{i=1}^r V(\sigma_i)\biggr)\\
&=& Y\bigcap \biggl(\bigcap_{i=1}^r\coprod_{\sigma_i\prec \gamma}
O_\gamma\biggr)\\
&=& Y\bigcap \biggl(\coprod_{\sigma_1,\ldots,
\sigma_r\prec\gamma}O_\gamma\biggr).
\end{eqnarray*}
So if $\sigma$ is the smallest cone in $\Sigma$ containing
$\sigma_1,\ldots, \sigma_r$, then we have
$$Y_{\sigma_1}\cap\cdots\cap Y_{\sigma_r}=Y\cap V(\sigma).$$
If there is no cone in $\Sigma$ containing $\sigma_1,\ldots,
\sigma_r$, then $Y_{\sigma_1}\cap\cdots\cap Y_{\sigma_r}$ is empty.
The condition of Lemma 3.5 (iv) holds with $$\phi(\{\sigma_1,\ldots,
\sigma_r\})={\rm dim}(Y_{\sigma_1}\cap\cdots\cap Y_{\sigma_r}).$$ So
we have
$$H_c^i\left(\biggl(Y-\bigcup_{\sigma\in J}Y_\sigma\biggr)\otimes_k \bar k,{j_Y}_\ast{\cal K}_\chi
\otimes g^\ast{\cal L}_\psi\right)=0$$ for all $i>{\rm dim}(Y)=n$,
that is,
$$H_c^i({\bf T}_{\bar k}^n, {\cal K}_\chi
\otimes f^\ast{\cal L}_\psi)=0$$ for all $i>n$. By the Weak
Lefschetz theorem and Poincar\'e duality, we have
$$H_c^i({\bf T}_{\bar k}^n, {\cal K}_\chi
\otimes f^\ast{\cal L}_\psi)=0$$ for all $i<n$. This proves (i).

(ii) By (i), we have
$${\rm dim} (H_c^n({\bf T}_{\bar k}^n,{\cal K}_\chi\otimes f^\ast{\cal
L}_\psi))=(-1)^n\chi_c({\bf T}_{\bar k}^n, {\cal K}_\chi\otimes
f^\ast{\cal L}_\psi).$$ So it suffices to show
$$\chi_c({\bf T}_{\bar k}^n, {\cal K}_\chi\otimes f^\ast{\cal L}_\psi)=(-1)^nn! {\rm
vol}(\Delta_\infty(f)).$$ We have
\begin{eqnarray*}
\chi_c({\bf T}_{\bar k}^n, {\cal K}_\chi\otimes f^\ast {\cal
L}_\psi)&=&\chi_c({\bf A}_{\bar k}^1, Rf_!({\cal K}_\chi\otimes
f^\ast {\cal L}_\psi))\\
&=& \chi_c({\bf A}_{\bar k}^1, Rf_!{\cal K}_\chi\otimes {\cal
L}_\psi).
\end{eqnarray*}
Since ${\cal L}_\psi$ is a lisse sheaf of rank 1 on ${\bf A}_k^1$,
for any closed point $x$ in ${\bf A}_{\bar k}^1$, we have
$${\rm sw}_x(Rf_!{\cal K}_\chi\otimes {\cal L}_\psi)
={\rm sw}_x(Rf_!{\cal K}_\chi).$$ Applying the
Grothendieck-Ogg-Shafarevich formula ([SGA 5] Expos\'e X, 7.1) to
both $\chi_c({\bf A}_{\bar k}^1, Rf_!{\cal K}_\chi\otimes {\cal
L}_\psi)$ and to $\chi_c({\bf A}_{\bar k}^1, Rf_!{\cal K}_\chi)$, we
see that
$$\chi_c({\bf A}_{\bar k}^1, Rf_!{\cal K}_\chi\otimes {\cal
L}_\psi)=\chi_c({\bf A}_{\bar k}^1, Rf_!{\cal K}_\chi) +{\rm
sw}_\infty(Rf_!{\cal K}_\chi)-{\rm sw}_\infty(Rf_!{\cal
K}_\chi\otimes {\cal L}_\psi).$$ By Lemma 3.6, $R^if_!{\cal
K}_\chi$ are tame at $\infty$ for all $i$. So $${\rm
sw}_\infty(Rf_!{\cal K}_\chi)=0.$$ Moreover, ${\cal L}_\psi$ has
Swan conductor $1$ at $\infty$ and has rank $1$. So we have
$${\rm sw}_\infty(Rf_!{\cal K}_\chi\otimes {\cal L}_\psi)={\rm
rank}(Rf_!{\cal K}_\chi)=\chi_c(f^{-1}(\bar\eta), {\cal
K}_\chi),$$ where $\eta$ is the generic point of ${\bf A}_{k}^1$.
Therefore, we have
\begin{eqnarray*}
\chi_c({\bf A}_{\bar k}^1, Rf_!{\cal K}_\chi\otimes {\cal
L}_\psi)&=& \chi_c({\bf A}_{\bar k}^1, Rf_!{\cal
K}_\chi)-\chi_c(f^{-1}(\bar\eta), {\cal K}_\chi).
\end{eqnarray*}
Note that
$$\chi_c({\bf A}_{\bar k}^1, Rf_!{\cal
K}_\chi)=\chi_c({\bf T}_{\bar k}^n, {\cal K}_\chi)=0.$$  So we
have
$$\chi_c({\bf A}_{\bar k}^1, Rf_!{\cal K}_\chi\otimes {\cal
L}_\psi)= -\chi_c(f^{-1}(\bar \eta), {\cal K}_\chi)$$ Therefore
$$\chi_c({\bf A}_{\bar k}^1, Rf_!{\cal K}_\chi\otimes {\cal
L}_\psi)= -\chi_c(f^{-1}(a), {\cal K}_\chi)$$ for sufficiently
general geometric point $a$ in ${\bf A}_{\bar k}^1$. We claim that
for sufficiently general $a$, the Laurent polynomial $G=f-a$ is
$0$-non-degenerate with respect to $\Delta(G)=\Delta_\infty(f)$.
Hence by Proposition 3.9, we have
\begin{eqnarray*}
\chi_c(f^{-1}(a), {\cal K}_\chi)&=& \chi_c({\bf T}_{\bar k}^n\cap
G^{-1}(0), {\cal K}_\chi)\\
&=&(-1)^{n-1}n!{\rm vol}(\Delta(G))\\
&=& (-1)^{n-1}n!{\rm vol}(\Delta_\infty(f)).
\end{eqnarray*}
Therefore $$\chi_c({\bf T}_{\bar k}^n, {\cal K}_\chi\otimes
f^\ast{\cal L}_\psi)=(-1)^nn! {\rm vol}(\Delta_\infty(f)).$$

Let's prove our claim. If $\tau$ is a face of $\Delta(G)$ that
does not contain $0$, then  $G_\tau=f_\tau$ (for any $a$). Since
$f$ is non-degenerate with respect to $\Delta_\infty(f)$, the
system of equations
$$\frac{\partial f_\tau}{\partial X_1}=\cdots=
\frac{\partial f_\tau}{\partial X_n}=0$$ has no solution in ${\bf
T}_{\bar k}^n$. So the locus of $G_\tau=0$ in ${\bf T}_k^n$ is
smooth over $k$ (for any $a$). Now suppose $\tau$ is a face of
$\Delta(G)$ containing $0$. Then $\Delta_\infty(f_\tau)=\tau$ and
$G_\tau=f_\tau-a$. Note that $f_\tau:{\bf T}_k^n\to {\bf A}_k^1$ is
non-degenerate with respect to $\Delta_\infty(f_\tau)$. By [DL] 4.3,
there exists a finite set $S$ of closed points in ${\bf A}_k^1$ such
that $f_\tau$ is smooth outside $f_\tau^{-1}(S)$. Choose $a$ to be
outside $S$. Then the locus of $G_\tau=f_\tau-a=0$ in ${\bf T}_k^n$
is smooth over $k$. So for sufficiently general $a$, $G$ is
$0$-non-degenerate with respect to $\Delta(G)$. This finishes the
proof of (ii).

(iii) By Proposition 3.4 (i), $g:Y\to {\bf A}_k^1$ is locally
acyclic relative to ${j_Y}_{!\ast}({\cal K}_\chi[n])$ outside
finitely many closed points. By Proposition 3.8, $R^ig_{\ast}
({j_Y}_{!\ast}({\cal K}_\chi[n]))$ are tame at $\infty$ for all $i$.
By Lemma 4.1 (ii), we have
$$H_c^i(Y\otimes_k \bar k, {j_Y}_{!\ast}({\cal K}_\chi[n])
\otimes g^\ast{\cal L}_\psi)=0$$ for all $i>0$. Replacing ${\chi}$
by ${\chi}^{-1}$ and $\psi$ by $\psi^{-1}$, we see
$$H_c^i(Y\otimes_k \bar k, {j_Y}_{!\ast}({\cal K}_{\chi^{-1}}[n])
\otimes g^\ast{\cal L}_{\psi^{-1}})=0$$ for all $i>0$. By Lemma
4.1 (i), this implies that
$$H^i(Y\otimes_k \bar k, {j_Y}_{!\ast}({\cal K}_{\chi^{-1}}[n])
\otimes g^\ast{\cal L}_{\psi^{-1}})=0$$ for all $i>0$. By
Poincar\'e duality, we have an isomorphism
$$H_c^i(Y\otimes_k \bar k, {j_Y}_{!\ast}({\cal K}_\chi[n])\otimes
g^\ast{\cal L}_\psi)\cong {\rm Hom}\biggl(H^{-i}\bigl(Y\otimes_k
\bar k, {j_Y}_{!\ast}({\cal K}_{\chi^{-1}}[n])\otimes g^\ast{\cal
L}_{\psi^{-1}}\bigr), \overline {\bf Q}_l(-n)\biggr).$$ Here we use
the fact that the Verdier dual of ${j_Y}_{!\ast}({\cal
K}_\chi[n])\otimes g^\ast{\cal L}_\psi$ is ${j_Y}_{!\ast}({\cal
K}_{\chi^{-1}}[n])\otimes g^\ast{\cal L}_{\psi^{-1}}(n)$. It follows
that
$$H^i_c(Y\otimes_k \bar k, {j_Y}_{!\ast}({\cal K}_{\chi}[n])
\otimes g^\ast{\cal L}_{\psi})=0$$ also for all $i<0$. Moreover,
the main theorem of [D] ([D] 3.3.1 and 6.2.3) implies that
$H^0_c(Y\otimes_k \bar k, {j_Y}_{!\ast}({\cal K}_{\chi}[n])
\otimes g^\ast{\cal L}_{\psi})$ is pure of weight $n$. For any
mixed object $K$ in the derived category $D_c^b(s, \overline {\bf
Q}_l)$, where $s={\rm Spec}\, k$, define its Poincar\'e polynomial
to be
$$P(K)=\sum_{w\in {\bf Z}}\sum_{i\in {\bf Z}}(-1)^ie_{iw}T^w,
$$ where $e_{iw}$ is the number of eigenvalues with weight $w$
counted with multiplicities of the geometric Frobenius element $F$
in ${\rm Gal}(\bar k/k)$ acting on $H^i(K_{\bar s})$. We often
write $P(K)$ as $P(K_{\bar s})$ by abuse of notations. By the
above discussion, we have
$$P(R\Gamma_c(Y\otimes_k\bar k, {j_Y}_{!\ast}({\cal K}_{\chi}[n])
\otimes g^\ast{\cal L}_{\psi}))= bT^n,$$ where $b={\rm
dim}(H^0_c(Y\otimes_k \bar k, {j_Y}_{!\ast}({\cal K}_{\chi}[n])
\otimes g^\ast{\cal L}_{\psi}))$. As $Y$ is the disjoint union of
$Y\cap O_\sigma$ $(\sigma\in \Sigma)$, and the cone $\sigma=0$
corresponds to $Y\cap {\bf T}_A^n\cong {\bf T}_k^n$, we have
\begin{eqnarray*}
&& P(R\Gamma_c({\bf T}_{\bar k}^n, {\cal K}_{\chi}\otimes
f^\ast{\cal L}_{\psi}))\\
&=& (-1)^n P(R\Gamma_c({\bf T}_{\bar k}^n,  {j_Y}_{!\ast}({\cal
K}_{\chi}[n]) \otimes
g^\ast{\cal L}_{\psi}))\\
&=& (-1)^n \biggl(P(R\Gamma_c(Y\otimes_k\bar k,
{j_Y}_{!\ast}({\cal K}_{\chi}[n]) \otimes g^\ast{\cal
L}_{\psi}))-\sum_{\sigma\not=0} P(R\Gamma_c((Y\cap
O_\sigma)\otimes_k\bar k, {j_Y}_{!\ast}({\cal
K}_{\chi}[n]) \otimes g^\ast{\cal L}_{\psi}))\biggr)\\
&=& (-1)^n \biggl(bT^n-\sum_{\sigma\not=0} P(R\Gamma_c((Y\cap
O_\sigma)\otimes_k\bar k, {j_Y}_{!\ast}({\cal K}_{\chi}[n]) \otimes
g^\ast{\cal L}_{\psi}))\biggr).
\end{eqnarray*}
We claim that if $0\not\in F_{\Delta_\infty(f)}(\sigma)$, then
$$H_c^i((Y\cap O_\sigma)\otimes_k \bar k, {j_Y}_{!\ast}({\cal
K}_\chi[n])\otimes g^\ast {\cal L}_\psi)=0$$ for all $i$. Indeed,
by Lemma 3.7, $R^j(g|_{Y\cap O_\sigma})_!({j_Y}_{!\ast}({\cal
K}_\chi[n])\otimes {\cal L}_\psi)$ are constant sheaves on ${\bf
A}_k^1$ for all $j$. It follows that
$$H_c^i({\bf A}_{\bar k}^1, R^j(g|_{Y\cap O_\sigma})_!({j_Y}_{!\ast}
({\cal K}_\chi[n]))\otimes {\cal L}_\psi)=0$$ for all $i,j$. Our
claim then follows from the spectral sequence
\begin{eqnarray*}
E_2^{ij}&=& H_c^i({\bf A}_{\bar k}^1, R^j(g|_{Y\cap
O_\sigma})_!({j_Y}_{!\ast}({\cal K}_\chi[n]))\otimes {\cal
L}_\psi)\cong  H_c^i({\bf A}_{\bar k}^1, R^j(g|_{Y\cap
O_\sigma})_!({j_Y}_{!\ast}({\cal K}_\chi[n])\otimes g^\ast {\cal
L}_\psi)) \\
&\Rightarrow& H_c^{i+j} ((Y\cap O_\sigma)\otimes_k\bar
k,{j_Y}_{!\ast}({\cal K}_\chi[n])\otimes g^\ast{\cal L}_\psi).
\end{eqnarray*}
So we have $$P(R\Gamma_c({\bf T}_{\bar k}^n, {\cal
K}_{\chi}\otimes f^\ast{\cal L}_{\psi}))=(-1)^n
\biggl(bT^n-\sum_{\sigma\not=0, \;0\in
F_{\Delta_\infty(f)}(\sigma)} P(R\Gamma_c((Y\cap
O_\sigma)\otimes_k\bar k, {j_Y}_{!\ast}({\cal K}_{\chi}[n])
\otimes g^\ast{\cal L}_{\psi}))\biggr).$$ For those $\sigma$ with
$0\in F_{\Delta_\infty(f)}(\sigma)$, by Lemma 3.2 (iii), $Y\cap
O_\sigma$ can be identified with ${O_\sigma}_k$ and $g|_{Y\cap
O_\sigma}:Y\cap O_\sigma\to {\bf A}_k^1$ can be identified with
$f_{\tau_\sigma}:{O_\sigma}_k\to {\bf A}_k^1$, where
${\tau_\sigma}= F_{\Delta_\infty(f)}(\sigma)$. Using Lemma 3.3,
one can check $({j_Y}_{!\ast}({\cal K}_\chi[n]))|_{Y\cap
O_\sigma}$ is identified with $(j_{!\ast}({\cal
K}_\chi[n]))|_{{O_\sigma}_k}$. So we have
$$P(R\Gamma_c({\bf T}_{\bar k}^n, {\cal
K}_{\chi}\otimes f^\ast{\cal L}_{\psi}))=(-1)^n
\biggl(bT^n-\sum_{\sigma\not=0, \;0\in
F_{\Delta_\infty(f)}(\sigma)} P(R\Gamma_c({O_\sigma}_{\bar k},
(j_{!\ast}({\cal K}_\chi[n]))|_{{O_\sigma}_k} \otimes
f_{\tau_\sigma}^\ast{\cal L}_{\psi}))\biggr).$$ By Lemma 2.3, if
${\cal K}_\chi$ is not the inverse image of a Kummer sheaf on
${O_\sigma}_k$ under the projection
$$p_\sigma: {\bf T}_k^n={\rm Spec}\, k[{\bf Z}^n]
\to {O_\sigma}_k={\rm Spec}\, k[{\bf Z}^n\cap \sigma^\perp],$$ then
$(j_{!\ast}({\cal K}_\chi[n]))|_{{O_\sigma}_k} $ is acyclic. Let $S$
be the set of those cones $\sigma$ in $\Sigma$ so that
$\sigma\not=0$, $0\in F_{\Delta_\infty(f)}(\sigma)$, and ${\cal
K}_\chi\cong p_\sigma^\ast{\cal K}_{\chi_\sigma}$ for a Kummer sheaf
${\cal K}_{\chi_\sigma}$ on ${O_\sigma}_k$. For each $\sigma\in S$,
let $\delta_\sigma$ be the image of $\check\sigma$ under the
projection ${\bf R}^n\to {\bf R}^n/{\sigma}^\perp$, let
$j_\sigma':{\rm Spec}\, k[{\bf Z}^n/{\bf Z}^n\cap \sigma^\perp]\to
X_k(\Sigma(\delta_\sigma))$ be the immersion of the open dense torus
in $X_k(\Sigma(\delta_\sigma))$, let $x_\sigma$ be the distinguished
point in $X_k(\Sigma(\delta_\sigma))$, and let
$\pi_\sigma:{O_\sigma}_k\to {\rm Spec}\,k$ be the structure
morphism. Then by Lemma 2.3, we have
$$(j_{!\ast}({\cal K}_\chi[n]))|_{{O_\sigma}_k}\cong
({\cal K}_{\chi_\sigma}[n-{\rm dim}(\sigma)])\otimes
\pi_{\sigma}^\ast x_\sigma^\ast( {j'_\sigma}_{!\ast}(\overline
{\bf Q}_l[{\rm dim}(\sigma)])).$$ Therefore we have
\begin{eqnarray*}
&&P(R\Gamma_c({\bf T}_{\bar k}^n, {\cal K}_{\chi}\otimes
f^\ast{\cal L}_{\psi}))\\
&=&(-1)^n \biggl(bT^n-\sum_{\sigma\in S}
P(R\Gamma_c({O_\sigma}_{\bar k}, ({\cal K}_{\chi_\sigma}[n-{\rm
dim}(\sigma)])\otimes \pi_{\sigma}^\ast x_\sigma^\ast(
{j'_\sigma}_{!\ast}(\overline {\bf Q}_l[{\rm
dim}(\sigma)]))\otimes
f_{\tau_\sigma}^\ast{\cal L}_{\psi}))\biggr)\\
&=& (-1)^n \biggl(bT^n-\sum_{\sigma\in S} (-1)^{n-{\rm
dim}(\sigma)} P(R\Gamma_c({O_\sigma}_{\bar k}, {\cal
K}_{\chi_\sigma}\otimes f_{\tau_\sigma}^\ast{\cal L}_{\psi}))
P(x_\sigma^\ast( {j'_\sigma}_{!\ast}(\overline {\bf Q}_l[{\rm
dim}(\sigma)])))\biggr)
\end{eqnarray*}
By Remark 2.2 and the fact that the polynomial
$\alpha(\delta_\sigma)$ involves only even powers of $T$, we have
$$P(x_\sigma^\ast( {j'_\sigma}_{!\ast}(\overline
{\bf Q}_l[{\rm dim}(\sigma)])))=(-1)^{{\rm
dim}(\sigma)}\alpha(\delta_\sigma).$$ So we get
$$P(R\Gamma_c({\bf T}_{\bar k}^n, {\cal K}_{\chi}\otimes
f^\ast{\cal L}_{\psi}))=(-1)^n \biggl(bT^n-\sum_{\sigma\in S}
(-1)^{n} P(R\Gamma_c({O_\sigma}_{\bar k}, {\cal
K}_{\chi_\sigma}\otimes f_{\tau_\sigma}^\ast{\cal L}_{\psi}))
\alpha(\delta_\sigma)\biggr).$$ By (i), we have
\begin{eqnarray*}
E({\bf T}_k^n, f,\chi)&=&(-1)^{n}P(R\Gamma_c({\bf T}_{\bar k}^n,
{\cal K}_{\chi}\otimes f^\ast{\cal L}_{\psi}))\\
E({O_{\sigma}}_k, f_{\tau_\sigma},\chi_\sigma)&=&(-1)^{n-{\rm
dim}(\sigma)}P(R\Gamma_c({O_\sigma}_{\bar k}, {\cal
K}_{\chi_{\sigma}}\otimes f_{\tau_\sigma}^\ast{\cal L}_{\psi})).
\end{eqnarray*}
So we finally get
$$E({\bf T}_k^n, f,\chi)=bT^n-\sum_{\sigma\in S}
(-1)^{{\rm dim}(\sigma)}E({O_{\sigma}}_k,
f_{\tau_\sigma},\chi_\sigma)\alpha(\delta_\sigma).$$ To get an
explicit formula for $E({\bf T}_k^n, f,\chi)$, we now take
$\Sigma=\Sigma(\Delta_\infty(f))$. By Proposition 1.2,
$$\tau\mapsto ({\rm cone}_{\Delta_{\infty}(f)}(\tau))^\vee$$
defines a one-to-one correspondence between faces of
$\Delta_\infty(f)$ and cones in $\Sigma(\Delta_\infty(f))$.
Moreover, for $\sigma_\tau=({\rm
cone}_{\Delta_{\infty}(f)}(\tau))^\vee$, we have
\begin{eqnarray*}
&&\tau=F_{\Delta_\infty(f)}(\sigma_\tau),\\
&&{\rm dim}(\sigma_\tau)=n-{\rm dim}(\tau),\\
&&\sigma_\tau^\perp=\check\sigma_\tau\cap (-\check \sigma_\tau)=
({\rm cone}_{\Delta_{\infty}(f)}(\tau))\cap (-{\rm
cone}_{\Delta_{\infty}(f)}(\tau))={\rm span}(\tau-\tau),
\end{eqnarray*}
and the image of $\check \sigma_\tau={\rm
cone}_{\Delta_{\infty}(f)}(\tau)$ in ${\bf R}^n/{\sigma_\tau}^\perp$
is just ${\rm cone}_{\Delta_{\infty}(f)}^\circ(\tau)$. Let $T$ be
the set of faces $\tau$ of $\Delta_\infty(f)$ so that
$\tau\not=\Delta_\infty(f)$, $0\in\tau$, and ${\cal K}_\chi\cong
p_\tau^\ast{\cal K}_{\tau}$ for a Kummer sheaf ${\cal K}_{\tau}$ on
${\bf T}_\tau={\rm Spec}\, k[{\bf Z}^n\cap {\rm span}(\tau-\tau)]$,
where
$$p_\tau:{\bf T}_k^n={\rm Spec}\, k[{\bf Z}^n]\to {\bf T}_\tau
={\rm Spec}\, k[{\bf Z}^n\cap {\rm span}(\tau-\tau)]$$ is the
projection. Note that $T$ corresponds to the set $S$ under the
one-to-one correspondence $\tau\mapsto ({\rm
cone}_{\Delta_{\infty}(f)}(\tau))^\vee$. We then have
$$E({\bf T}_k^n, f,\chi)=bT^n-\sum_{\tau\in T}
(-1)^{n-{\rm dim}(\tau)}E({\bf T}_\tau,
f_{\tau},\chi_\tau)\alpha({\rm
cone}_{\Delta_\infty(f)}^\circ(\tau)).$$ By (ii), we have
\begin{eqnarray*}
E({\bf T}_k^n, f,\chi)(1)&=&n!{\rm vol}(\Delta_\infty(f)),\\
E({\bf T}_\tau, f_\tau,\chi_\tau)(1)&=&({\rm dim}(\tau))!{\rm
vol}(\tau).
\end{eqnarray*}
Evaluating the above equality at 1, we get
$$b=n!{\rm vol}(\Delta_\infty(f))+\sum_{\tau\in T}
(-1)^{n-{\rm dim}(\tau)}({\rm dim}(\tau))!{\rm vol}(\tau)\alpha({\rm
cone}_{\Delta_\infty(f)}^\circ(\tau))(1),$$ that is,
$b=e(\Delta_\infty(f),\chi).$ So we have
$$E({\bf T}_k^n, f,\chi)=e(\Delta_\infty(f),\chi)T^n-\sum_{\tau\in T}
(-1)^{n-{\rm dim}(\tau)}E({\bf T}_\tau,
f_{\tau},\chi_\tau)\alpha({\rm
cone}_{\Delta_\infty(f)}^\circ(\tau)).$$ Using this expression,
the definition of $E(\Delta_{\infty}(f),\chi)$, and induction on
${\rm dim}(\Delta_\infty(f))$, we get
$$E({\bf T}_k^n, f,\chi)=E(\Delta_\infty(f),\chi).$$
By [D] 3.3.1 and 3.3.3, $H_c^n({\bf T}_{\bar k}^n, {\cal
K}_{\chi}\otimes f^\ast{\cal L}_{\psi}))$ is mixed with weights
between $0$ and $n$. So $E({\bf T}_k^n, f,\chi)$ is a polynomial
of degree $\leq n$. By the definition of $\alpha$, we have
$${\rm deg}(\alpha({\rm
cone}_{\Delta_\infty(f)}^\circ(\tau)))\leq {\rm dim}({\rm
cone}_{\Delta_\infty(f)}^\circ(\tau))-1=n-{\rm dim}(\tau)-1.$$
Moreover, we have
$${\rm deg}(E({\bf T}_\tau, f_\tau,\chi_\tau))\leq {\rm dim}(\tau).$$ It now
follows from the expression
$$E({\bf T}_k^n, f,\chi)=e(\Delta_\infty(f),\chi)T^n-\sum_{\tau\in T}
(-1)^{n-{\rm dim}(\tau)}E({\bf T}_\tau,
f_{\tau},\chi_\tau)\alpha({\rm
cone}_{\Delta_\infty(f)}^\circ(\tau))$$ that
$$e_n=e(\Delta_\infty(f),\chi).$$

(iv) Since $0$ is an interior point of $\Delta_\infty(f)$, for any
nonzero $\sigma\in\Sigma$, we have $0\not\in
F_{\Delta_\infty(f)}(\sigma)$. We have seen in the proof of (iii)
that this implies
$$H_c^i((Y\cap O_\sigma)\otimes_k \bar k, {j_Y}_{!\ast}({\cal
K}_\chi[n])\otimes g^\ast {\cal L}_\psi)=0$$ for all $i$. Since
$Y-Y\cap {\bf T}_A^n$ is the disjoint union of $Y\cap O_\sigma$
for nonzero $\sigma$, we have
$$H_c^i((Y-Y\cap {\bf T}_A^n)\otimes_k
\bar k,{j_Y}_{!\ast}({\cal K}_\chi[n]) \otimes g^\ast{\cal
L}_\psi)=0$$ for all $i$. (We can also apply Lemma 3.5 (iii) to
$Y-Y\cap {\bf T}_A^n$, $J=\{\sigma\in \Sigma|{\rm dim}(\sigma)=1\}$,
and the closed subschemes $(Y-Y\cap {\bf T}_A^n)\cap V(\sigma)=Y\cap
V(\sigma)$ $(\sigma\in J)$.) So we have
$$H_c^i((Y\cap {\bf T}_A^n)\otimes_k \bar k, {j_Y}_{!\ast}({\cal
K}_\chi[n]) \otimes g^\ast{\cal L}_\psi)\cong H_c^i
(Y\otimes_k\bar k, {j_Y}_{!\ast}({\cal K}_\chi[n])\otimes
g^\ast{\cal L}_\psi)$$ for all $i$, that is,
$$H_c^{i+n}({\bf T}_{\bar k}^n, {\cal K}_\chi
\otimes f^\ast{\cal L}_\psi)\cong H_c^i (Y\otimes_k\bar k,
{j_Y}_{!\ast}({\cal K}_\chi[n])\otimes g^\ast{\cal L}_\psi)$$ for
$i$. In the proof of (iii), we have seen $H_c^0(Y\otimes_k\bar k,
{j_Y}_{!\ast}({\cal K}_\chi[n])\otimes g^\ast{\cal L}_\psi)$ is
pure of weight $n$. So $H_c^{n}({\bf T}_{\bar k}^n, {\cal K}_\chi
\otimes f^\ast{\cal L}_\psi)$ is pure of weight $n$.

\bigskip
\bigskip
\noindent {\bf References.}

\bigskip
\bigskip
\noindent [AS] A. Adolphson and S. Sperber, {\it Twisted exponential
sums and Newton polyhedra}, J. Reine Angew. Math. 443 (1993),
151-177.

\bigskip
\noindent [BBD] A. Beilinson, J. Bernstein and P. Deligne, {\it
Faisceaux pervers}, in {\it Analyse et Topologie sur les Espace
Singuliers (I)}, Ast\'erique 100 (1980).

\bigskip
\noindent [D] P. Deligne, {\it La conjecture de Weil II}, Publ.
Math. IHES 52 (1980), 137-252.

\bigskip
\noindent [DL] J. Denef and F. Loeser, {\it Weights of exponential
sums, intersection cohomology, and Newton polyhedra}, Invent. Math.
106 (1991), 275-294.

\bigskip
\noindent [F] W. Fulton, {\it Introduction to toric varieties},
Annals of Math. Studies (131), Princeton University Press 1993.

\bigskip
\noindent [I1] L. Illusie, {\it Th\'eorie de Brauer et
Caract\'eristique d'Euler-Poincar\'e}, in {\it Caract\'eristique
d'Euler-Poincar\'e}, Ast\'erique 82-83 (1981), 161-172.

\bigskip
\noindent [I2] L. Illusie, {\it Autour du th\'eor\`eme de monodromie
locale}, in {\it P\'eriods p-adiques}, Ast\'erique 223 (1994), 9-57.

\bigskip
\noindent [SGA] S\'eminaire de G\'eom\'etrie Alg\'ebrique du
Bois-Marie.

\bigskip
\noindent [S] R. Stanley, {\it Generalized H-vectors, intersection
cohomology of toric varieties, and related results,} in {\it
Commutative Algebra and Combinatorics} edited by N. Nagata and H.
Matsumura, Adv. Stud. Pure Math., vol 11, 187-213, Amsterdam, New
York, North-Holland 1987.

\bigskip
\noindent [SGA 1] Rev\^{e}tements \'etales et groupe fondemental,
by Grothendieck, {\it Lecture Notes in Mathematics} 224,
Springer-Verlag (1971).

\bigskip
\noindent [SGA 4] Th\'eorie des topos et cohomologie \'etale des
sch\'emas, by M. Artin, A. Grothendieck and J.-L. Verdier, {\it
Lecture Notes in Mathematics} 269, 270, 305, Springer-Verlag
(1972-1973).

\bigskip
\noindent [SGA 4${1\over 2}$] Cohomologie \'etale, by P. Deligne,
{\it Lecture Notes in Mathematics} 569, Springer-Verlag (1977).

\bigskip
\noindent [SGA 5] Cohomologie $l$-adique et fonctions L, by
Grothendieck, {\it Lecture Notes in Mathematics} 589,
Springer-Verlag (1977).

\bigskip
\noindent [SGA 7] Groupes de monodromie en g\'eom\'etrie
alg\'ebrique, I by A. Grothendieck, II by P. Deligne and N. Katz,
{\it Lecture Notes in Mathematics} 288, 340, Springer-Verlag
(1972-1973).
\end{document}